%% file: expl_self_avoid.tex
\documentclass[11pt]{amsart}
\usepackage{amssymb}
\usepackage{pgfplots}
\usepackage{csvsimple}
\usepackage{hyperref}
\usepackage{graphicx,color}
\usepackage[normalem]{ulem}
\usepackage{mathrsfs}
\include{sb_macros}
\include{pr_macros}
\newcommand{\SB}[1]{#1}
\newcommand{\sbb}[1]{#1} 
\newcommand{\sbbb}[1]{#1}

\def\TP{{\rm TP}}
\newcommand{\radius}{R}
\newcommand{\korr}[1]{#1} 
\allowdisplaybreaks

\begin{document}
\title[Gradient flows with tangent-point functionals]{Stability
of a simple scheme \\ for the approximation of elastic knots \\ and self-avoiding inextensible curves}
\author[S. Bartels]{S\"oren~Bartels}
 \address{Abteilung f\"ur Angewandte Mathematik,  
 Albert-Ludwigs-Universit\"at Freiburg, Hermann-Herder-Str.~10, 
 79104 Freiburg i.~Br., Germany}
 \email{bartels@mathematik.uni-freiburg.de}
\author[Ph. Reiter]{Philipp~Reiter}
 \address{Department of Mathematics,
University of Georgia,
Athens GA 30602, USA}
 \email{reiter@uga.edu}
\date{\today}
\renewcommand{\subjclassname}{%
  \textup{2010} Mathematics Subject Classification}
 \subjclass[2010]{65N12 (57M25 
 65N15 65N30)}
 \begin{abstract}
 We discuss a semi-implicit numerical scheme that allows for
 minimizing the bending energy of curves within
 certain isotopy classes.
 To this end we consider a \sbb{weighted sum} of
 the bending energy and the tangent-point functional.
 
 Based on estimates for the second derivative
 of the latter and a uniform bi-Lipschitz radius, we prove a stability result implying
 energy decay during the evolution as well as
 maintenance of arclength parametrization.
 
 Finally we present some numerical experiments
 exploring the energy landscape,
 targeted to the question how to obtain global
 minimizers of the bending energy in knot classes,
 so-called elastic knots.
 \end{abstract}
 \keywords{Self-avoidance, curves, stability,
 bending energy, knot energy, elastic knots, tangent-point energies}
\maketitle

\setcounter{tocdepth}1
\tableofcontents

\section{Introduction}
We aim at numerically detecting configurations of embedded curves
with low bending energy within certain knot classes. For this we 
define an energy functional as a weighted sum of an elastic bending 
energy term and a tangent-point functional that prevents curves from
self-intersecting and pulling-tight of small knotted arcs,
i.e., from leaving the ambient isotopy class.
To this end we define
\begin{equation}\label{eq:total}
 E(u) = \frac{\k}{2} \int_I |u''(x)|^2 \dv{x} + \vrho \TP(u)
\end{equation}
on the class of embedded and arclength parametrized curves 
$u:I\to \R^3$
where $\k$ and $\vrho$ are positive parameters.
We mainly consider the case of periodic unit-length intervals
$I=\R/\Z$, however, the setting can also be extended to arbitrary
intervals~$I\subset\R$ and suitable boundary conditions.
The regularization ansatz~\eqref{eq:total} has already been
discussed by von der Mosel~\cite{vdM:meek} with
O'Hara's energies~\cite{oha:en-fam} in place of the tangent-point functional.
In fact, one might conjecture that any self-avoiding functional
will qualitatively produce the same results.

We are particularly interested in the case $\k\gg\vrho$
which corresponds to the idea of a very thin knotted springy wire.
Disregarding twist and other physical quantities we assume that its 
behavior is only driven by the bending
energy of its centerline.
This case is computationally challenging since strong
forces related to bending effects have to be compensated by repulsive
forces related to the tangent-point functional to avoid self-intersections.
\sbb{We rigorously show in a semi-discrete setting 
that the energy of our time-stepping scheme is monotonically decreasing 
under a moderate condition on the step size which we call 
{\em energy stability}.}

\sbb{Our scheme and numerical analysis apply to arbitrary choices of
parameters $\k>0$ and $\vrho \ge 0$ but the stability conditions
become more restrictive as $\k/\vrho \to 0$. The
experiments provided in this article give rise to the following conclusions:
\begin{itemize}
\item[(1)] Our numerical scheme is energy stable if the time-step size
satisfies a condition $\tau \le F(\k,\vrho)$  with $F(\k,\vrho) = \cO(1)$ 
if $\k \gg \vrho$ and $F(\k,\vrho) \to 0$ as $\k/\vrho \to 0$. A mesh
dependence is observed if $\vrho \gg \k$ which appears to be related to
errors introduced by quadrature for $\TP$. 
\item[(2)] Our scheme preserves the isotopy class of a closed curve if 
it is energy stable and if for fixed $\k>0$ the spacial mesh size $h>0$ 
satisfies $h \le G(e_0,\vrho)$ with the initial energy $e_0 = E(u^0)$
where $G(e_0,\vrho)\to 0$ as $e_0 \to \infty$ or $\vrho \to 0$. 
\end{itemize}
The condition $\tau \le F(\k,\vrho)$ is related to our explicit treatment
of the tangent-point functional while the condition $h\le G(e_0,\vrho)$ ensures
that the discretized tangent-point functional defines a sufficiently large
energy barrier relative to the initial energy to prevent isotopy changes.} 

We believe that our scheme can be helpful in investigating properties of 
knots and elastic curves. We have therefore developed a flexible online tool 
which is available under the following adress:  
\begin{center}
\url{https://aam.uni-freiburg.de/agba/forschung/knotevolve/} \\
\end{center}
A password is required to work with the current beta version which can be 
obtained via email from the authors. Further details about the tool 
will be provided in~\cite{BaFaRe18-pre}.

While our numerical scheme can help understanding 
configurations with low bending energy, it remains a challenging task to rigorously
predict the equilibrium shapes depending on isotopy classes.
Langer and Singer~\cite{LS} have proven that the circle is the
only local minimizer in the unknot class.
For the trefoil (and all other two bridge torus knot classes)
global minimizers tend to the doubly covered circle as $\vrho\searrow0$~\cite{GRvdM}.
These limit objects are referred to as \emph{elastic knots}.

The characterization of elastic knots for general knot classes
is wide open. Langer and Singer~\cite{LS},
Gallotti and Pierre-Louis~\cite{gallotti},
as well as Sossinsky in different collaborations (see~\cite{sossinsky}
and references therein) have carried out several experiments in that direction.

\subsection*{Bending energy}

The bending energy of a curve is the integral of its squared curvature.
It has been proposed as a simple model for the energy of a thin springy
wire almost three hundred years ago by Daniel Bernoulli.
It can be seen as one of the most elementary examples of nonlinear
functionals and plays a fundamental r\^ole in elasticity theory.
One-dimensional bending theories have been rigorously derived by
Mora and M\"uller~\cite{MorMul03}.

We find applications in different fields such as
the modeling of cell filaments (Manhart et al.~\cite{MOSS15}),
textile fabrication processes (Grothaus and Marheineke~\cite{GroMar16}),
and computer graphics (Wardetzky et al.~\cite{WBHZG07}).

In recent time, gradient flows have received much attention,
with respect to rigorous analysis, see Dziuk et al.~\cite{DzKuSc02}
as well as regarding discretization aspects,
see Deckelnick and Dziuk~\cite{DeDz09},
Barrett et al.~\cite{BaGaNu12},
Bartels~\cite{Bar13},
Dall'Acqua et al.~\cite{DaLiPo14},
Pozzi and Stinner~\cite{PozSti17}.

In this article, we rely on the scheme introduced
by the first author in~\cite{Bar13}. Implementing a constraint
ensuring that the curves stay close to arclength parametrization
if the initial curve is arclength parametrized,
the bending energy is replaced by the squared $L^{2}$ norm
of the second derivative of the curve in~\eqref{eq:total}.
It is obvious that both functionals agree in case of arclength 
parametrization, and the same applies for their derivatives
in normal directions.

\subsection*{Tangent-point functional}
The \emph{tangent-point radius} $r(u(y), u(x))$ is the radius of the circle that is tangent to $u$ at the point $u(y)$ and that intersects with $u$ in $u(x)$.
The tangent-point functional is essentially the $q$-th power of the corresponding
$L^{q}$ norm \sbbb{of its reciprocal value}, 
i.e., for arclength parametrized curves we have
\begin{equation}\label{eq:tp}
 u\mapsto\frac1{2^{q}q}\iint_{I\times I}\frac{\d x\d y}{r(u(y),u(x))^{q}},
 \qquad q>2.
\end{equation}

The tangent-point energies
have been proposed as a family of self-avoiding functionals by
Gonzalez and Maddocks~\cite{GM99};
the scale invariant case $q=2$ (which we will disregard in this paper) already appears in a paper by Buck and Orloff~\cite{BuckOrloff}.
They are defined on
(smooth) embedded curves $u:I\to\R^{n}$
and take values in~$[0,+\infty]$,
see Strzelecki and von der Mosel~\cite{StrMos12} and references therein.
Blatt~\cite{Blatt13} has characterized the energy spaces
in terms of Sobolev--Slobodecki{\u\i} spaces;
regularity aspects are discussed in~\cite{BR15}.

The main feature of the tangent-point energies and many further
so-called knot energies~\cite{oha:book} is that they provide
a monotonic uniform bound on the bi-Lipschitz constant.
This implies in particular that the energy values of a sequence of
embedded curves converging
to a curve with a self-intersection will necessarily blow up.

Compared to other knot energies,
the tangent-point energies seem to be particularly well suited for numerics.
In contrast to the \emph{ropelength} functional, i.e., 
the quotient of length over thickness~\cite{GM99}, they are
smooth. Its variations have integrable integrands and do not 
contain intrinsic terms which can be an issue in the case of 
O'Hara's energy family~\cite{oha:en-fam}.
Furthermore, the evaluation of the tangent-point energies only
requires the evaluation of a double integral, \sbbb{while}
the integral Menger curvature is defined by a triple integral.
We refer to Scholtes~\cite{Schol16} for an outline of
the discretization of several self-avoiding energies.
A scheme for the gradient flow of the integral Menger curvature
has been devised by
Hermes~\cite{Hermes12}.

We will denote by $\TP$ a variant, see~\eqref{eq:vartp} below,
which agrees with the functional~\eqref{eq:tp} up to
a multiplicative factor which amounts to one on arclength parametrized
curves. So, prescribing an arclength constraint, we may replace~\eqref{eq:tp} by $\TP$,
thereby avoiding additional multiplicative terms which
would just extend the already quite involved derivative formulae.
Most results on~\eqref{eq:tp} carry over to $\TP$ by minor modifications
of the corresponding arguments.
In contrast to~\eqref{eq:tp}, the functional $\TP$ is \emph{not} invariant under
reparametrization.
Several aspects including
error estimates for a spacial discretization of $\TP$
and its first variation have been derived in~\cite{BRR}.

Note that, besides its self-avoiding feature, the tangent-point
functional is a curvature energy (as well as the other functionals
mentioned above).
For $\k\approx\vrho$ we may expect a similar qualitative behavior
compared to the case $\k=0$.
In the latter situation we expect nonplanar minimizers
(in nontrivial knot classes) as the self-avoiding potential clearly dominates.
The advantage of $\k>0$ is that we can employ a semi-implicit scheme
which does not involve implicit terms in a fractional operator.
\sbb{We refer the reader to Lin and Schwetlick~\cite{LinSchwe10} for an 
experimental investigation of the case $\kappa=\rho$ using the M\"obius 
energy instead of $\TP$.}

\subsection*{Gradient flow}
We consider a gradient flow to model certain dynamics and to define
\sbb{a family of arclength parametrized curves} $(u_t)_{t\in [0,T]}$ that converges to a stationary 
configuration for $E$. Given an inner product $(\cdot,\cdot)_{\korr X}$ on 
$H^2(I;\R^3)$ the evolution is specified by the parabolic system
\begin{equation}\label{eq:gradflow}
 (\p_t u, w)_{\korr X} = -\delta E(u)[w] = -\k (u'',w'')_{\korr{L^{2}}} - \vrho \delta \TP(u)[w],
\end{equation}
for all admissible $w$,
subject to initial and boundary conditions 
\[
u(0,\cdot) = u_0, \quad L_{bc}(u(t,\cdot)) = \ell,
\]
and subject to the linearized arclength condition
\[
\p_t u'(t,\cdot) \cdot u'(t,\cdot) = 0, \quad w' \cdot u'(t,\cdot) = 0.
\]
Here $L_{bc}$ denotes a suitable linear operator which imposes, e.g.,
periodic or clamped boundary conditions.

\subsection*{Time-stepping}

Time-stepping schemes provide the basis for numerical methods
and the generally most stable approach uses a fully implicit treatment
of the nonlinearities. This, however, is of limited practical use
since the variation of the tangent-point functional defines a nonlocal
\SB{and nonlinear} operator. We therefore aim at analyzing schemes that treat this term
explicitly, e.g., in the form
\begin{equation}\label{eq:discgradflow}
 (d_t u^k,w)_{X} + \k ([u^k]'',w'')_{L^{2}} = - \vrho \delta \TP(u^{k-1})[w].
\end{equation}
This avoids inverting a fully populated matrix related to 
a spacially discretized tangent-point functional and only requires solving 
sparse linear systems of equations in the time steps. Here, $d_t$ is
the backward difference quotient operator, defined with a step-size
$\tau>0$ via
\[
d_t u^k = \frac{1}{\tau} (u^k-u^{k-1}).
\]
For ease of presentation we restrict our stability analysis to the discretized 
$H^2$ flow. In some cases, e.g., 
when $\k$ is small, a weaker flow may be of interest. By making appropriate use
of inverse estimates in a spacially discrete setting, our arguments can be carried
over leading to a stability result under more restrictive step size conditions. 
Finally, we note that the $H^2$ flow may serve as an iterative solver 
for the minimization problems in the time steps of a fully implicit discretization
of an $L^2$ flow. Higher order gradient flows have recently also been 
used in~\cite{Schu18} for the simulation of elastic knots. 

\subsection*{Outline}

The paper is organized as follows.
In Section~\ref{sec:stab} we prove the stability result
(Proposition~\ref{prop:stab})
which, among other techniques, is based on estimates on the tangent-point functional $\TP$
that are derived in Section~\ref{sec:tp}.
Here we first discuss some general remarks on $\TP$,
including a characterization of the energy spaces
and an approximation result on the relation of~\eqref{eq:tp}
and~$\TP$.
Then we derive \SB{certain} estimates on the bi-Lipschitz constant
which leads to a uniform bi-Lipschitz radius in Corollary~\ref{cor:bil}.
Furthermore, we derive estimates for $\TP$ and its derivatives
in Propositions~\ref{prop:0th}, \ref{prop:1st}, and~\ref{prop:2nd}
as well as~\eqref{eq:kth}.
Numerical experiments in Section~\ref{sec:experiments}
confirm the good stability properties of the proposed scheme and its
suitability to maintain the isotopy class.
We experimentally explore the complex energy landscape 
defined by the functional~$E$ in Section~\ref{sec:elastic} by using our
scheme to relax the energy of various knot configurations.

\section{Semi-discrete stability analysis}\label{sec:stab}


We start with an auxiliary statement providing uniform bounds
on the derivatives of the tangent-point functional.
Its proof relies on statements which are derived in Section~\ref{sec:tp}.

Recall that, for practical purposes, we replace the tangent-point energy~\eqref{eq:tp} by a variant~\eqref{eq:vartp}
which is analytically less involved and coincides with
the original functional on arclength parametrized curves.

Throughout what follows we write $W^{s,q}$ instead of $W^{s,q}(\R/\Z,\R^n)$
\sbbb{and} if $q=2$ we use the symbol $H^s$.

\begin{lemma}\label{lem:2estimates}\label{la:vars_tp}
 There are constants $c_{1},c_{2},\radius>0$
 only depending on $0<\lambda\le\Lambda$, $M>0$, $n$, and $q$
 such that any embedded and regular curve $u\in W^{2-1/q,q}$ with
 \[ \lambda\le\abs{u'}\le\Lambda \qquad\text{and}\qquad \TP(u)\le M \]
 satisfies
 \[ \abs{\delta\TP(u)[w]} \le c_{1}\norm{u'}_{W^{1-1/q,q}}^{q}\norm{w'}_{W^{1-1/q,q}} \]
 and, for any $z\in W^{2-1/q,q}$ with $\norm{z'}_{\korr{W^{1-1/q,q}}}\le\radius$,
 \[ \abs{\delta^{2}\TP(u+z)[v,w]} \le c_{2}\br{\norm{u'}_{W^{1-1/q,q}}^{2q+2}+1}\norm{v'}_{W^{1-1/q,q}}\norm{w'}_{W^{1-1/q,q}}. \]
\end{lemma}

In particular, the second estimate applies to a Taylor expansion
\[ \TP(\tilde u) = \TP(u) + \delta\TP(u)[\tilde u-u]
+ \int_{0}^{1}(1-s)\delta^{2}\TP\br{u+s(\tilde u-u)}[\tilde u-u,\tilde u-u]\d s \]
provided
\[ \norm{(\tilde u-u)'}_{\korr{W^{1-1/q,q}}}\le\radius. \]

\begin{proof}[Proof of Lemma~\ref{lem:2estimates}]
 The first statement is an immediate consequence
 of Propositions~\ref{prop:1st} and~\ref{prop:unibi}.
 To show the second one, we first infer
 \begin{align*}
  &\abs{\delta^{2}\TP(u+z)[v,w]} \\
  &\le C\br{\bil\br{u+z}}^{2q+4}\norm{\br{u+z}'}_{W^{1-1/q,q}}^{2q+2}\cdot\norm{v'}_{W^{1-1/q,q}}\norm{w'}_{W^{1-1/q,q}}
 \end{align*}
 from Proposition~\ref{prop:2nd}. We have
 \[ \norm{\br{u+z}'}_{W^{1-1/q,q}}^{2q+2}
 \le C\br{\norm{u'}_{W^{1-1/q,q}}^{2q+2} + \radius^{2q+2}}
 \le C'\br{\norm{u'}_{W^{1-1/q,q}}^{2q+2} + 1}. \]
 \korr{%
 From Corollary~\ref{cor:bil} we obtain a radius $\widetilde\radius
 =\widetilde\radius(\lambda,\Lambda,M,n,q)>0$
 such that $\norm{z'}_{L^{\infty}}\le\widetilde\radius$ implies
 $\bil\br{u+z} \le 2\bil u$.
 Due to the embedding $W^{1-1/q,q}\hookrightarrow L^{\infty}$
 we have $\norm{z'}_{L^\infty}\le\tilde c_{q} \norm{z'}_{W^{1-1/q,q}}$
 for some $\tilde c_{q}>0$,
 so $\radius=\tilde c_{q}^{-1}\widetilde\radius$ yields the desired.}
\end{proof}

\SB{We consider a periodic setting and an $H^2$ flow which allows us to 
make use of certain Sobolev and Poincar\'e inequalities.}

\korr{%
\begin{remark}[Poincar\'e inequalities]
 For periodic functions $w$ we have $\int_{\R/\Z}w'=0$
 which implies $\norm{w'}_{L^2}\le\norm{w''}_{L^2}$.
 Owing to Sobolev embeddings $H^{1}\hookrightarrow W^{1-1/q,q}$
 and $H^{1}\hookrightarrow L^{\infty}$ respectively, we infer
 \begin{equation}\label{eq:bc}
 {\|w'\|_{W^{1-1/q,q}}} \le c_q \|w''\|_{L^{2}}, \qquad \|w'\|_{L^\infty} \le c_\infty {\|w''\|_{L^2}}
 \end{equation}
 where $c_{q}$ depends on $n$ and $q\in(2,4]$ and $c_{\infty}$ on~$n$. 
\end{remark}%
}%

\SB{The following proposition shows that our numerical scheme for
approximating the $H^2$ flow of the energy functional $E$ is energy stable
under a moderate step size condition.}

\begin{proposition}\label{prop:stab}
For given $\kappa>0$, $\rho\ge0$, and $q\in(2,4]$, \SB{let $(u^k)_{k=0,\dots,K} \subset H^2$}
be the uniquely defined sequence 
generated for an initial $u^0$ with $|[u^0]'|^2 =1$ and $E(u^0)< \infty$ by the scheme
\[
(d_t u^k,w)_{H^{2}} + \k ([u^k]'',w'')_{L^{2}} = -\vrho \delta \TP(u^{k-1})[w],
\]
subject to the linearized arclength conditions
\[
[d_t u^k]'\cdot [u^{k-1}]' = 0, \qquad w'\cdot [u^{k-1}]' = 0.
\]
There exists $c'=c'(\k,\vrho,e_0,n,q)>0$ with $e_0 = E(u^0)$ which 
is independent of $\tau>0$ such that if 
$\tau c' \le 1/2$ then we have  the energy stability property
\begin{equation}\label{eq:stability}
 E(u^L) + (1-c'\tau) \tau \sum_{\ell=1}^L \|d_t u^\ell\|_{H^2}^2 \le E(u^0)
\end{equation}
for all $0\le L \le K$. Moreover, we have that
\begin{equation}\label{eq:arclength}
\max_{k=0,\dots,K} \big\| |[u^k]'|^2 - 1 \big\|_{L^\infty} \le \tau c_\infty^2 4 e_0. 
\end{equation}
\end{proposition}

\begin{proof}
It follows from the Lax--Milgram lemma that the iterates 
$(u^k)_{k=0,\dots,K}$ are well defined. To prove the asserted
energy law we argue by induction over $K$, and note
that for $K=0$ it is trivially satisfied. Assume that the estimate
holds with $K$ replaced by $K-1$ for some $K\ge 1$. For every
$1 \le k\le K$ we then have, assuming $c' \tau \le \tfrac12$, that
\[
\frac{\k}{2} \|[u^{k-1}]''\|_{L^{2}}^2 + \frac{\tau}{2} \sum_{\ell=1}^{k-1} \|d_t u^\ell\|_{H^2}^2
\le  e_0. 
\]
Choosing $w=d_t u^k$ in the scheme and using the first estimate
of Lemma~\ref{la:vars_tp} with~\eqref{eq:bc} yields that 
\begin{align*}
&\|d_t u^k\|_{H^2}^2 + \k \frac{d_t}{2} \|[u^k]''\|_{L^{2}}^2 
+ \k \frac{\tau}{2} \|[d_t u^k]''\|_{L^{2}}^2 \\
&\le \frac{\vrho^2}{2} c_1^2 c_{q}^{2} \korr{\|[u^{k-1}]'\|_{W^{1-1/q,q}}^{2q}}
+ \frac1{2c_{q}^{2}} \underbrace{\korr{\|[d_t u^k]'\|_{W^{1-1/q,q}}^2}}_{\le c_{q}^{2}\norm{d_{t}u^{k}}_{H^{2}}^{2}}.
\end{align*}
Using
\begin{equation}\label{eq:uW}
 \korr{\|[u^{k-1}]'\|_{W^{1-1/q,q}}} \le c_q \|[u^{k-1}]''\|_{L^{2}} \le c_q \k^{-1/2} e_0^{1/2}\sqrt2,
\end{equation}
\korr{owing to}~\eqref{eq:bc} this implies that 
\begin{equation}\label{eq:ca}
 \frac12 \|d_t u^k\|_{H^2}^2 + \k \frac{d_t}{2} \|[u^k]''\|_{L^{2}}^2
 \le \frac{\vrho^2}{2} c_1^2 c_{q}^{2} \big( c_{q}^{2q}\kappa^{-q}e_{0}^{q}2^{q} \big) = c_a.
\end{equation}
Imposing the condition $\tau c_a \le e_0$ shows that
\begin{equation}\label{eq:ener_subopt}
\frac{\tau}{2} \sum_{\ell=1}^k \|d_t u^\ell \|_{H^{2}}^2 + \frac{\k}{2} \|[u^k]''\|_{L^{2}}^2 
\le e_0 + \tau c_a \le 2 e_0.
\end{equation}
With this (suboptimal) auxiliary bound we aim at deriving the asserted energy
law up to level~$K$. We first note that the iterates $(u^k)_{k=0,\dots,L}$
satisfy, because of the linearized arclength condition, 
\[
|[u^k]'|^2 = |[u^{k-1}]'|^2 + \tau^2 |[d_t u^k]'|^2 
= \dots = |[u^0]'|^2 + \tau^2 \sum_{\ell=1}^k |[d_t u^\ell]'|^2.
\]
Since $|[u^0]'|^2 = 1$ we find that $|[u^k]'|^2 \ge 1$ and \SB{verify~\eqref{eq:arclength}}.
 To deduce the energy law we again choose $w=d_t u^k$ in the scheme 
but this time use a Taylor expansion of the potential, i.e., 
\begin{equation}\label{eq:ener_taylor}
\begin{split}
&\|d_t u^k\|_{H^2}^2 + \k \frac{d_t}{2} \|[u^k]''\|_{L^{2}}^2 
+ \k \frac{\tau}{2} \|[d_t u^k]''\|_{L^{2}}^2 \\
& = -\vrho \delta \TP(u^{k-1})[d_t u^k]  \\
& = -\vrho d_t \TP(u^k) + \tau \vrho
\int_{0}^{1}(1-s)\delta^{2}\TP\br{u^{k-1}+s\tau d_{t}u^{k}}[d_{t}u^{k},d_{t}u^{k}]\d s.
\end{split}
\end{equation}
\sbbb{Requiring that
\begin{equation}\label{eq:radius}
\tau 4 e_0 \le \radius^2 c_q^{-2}
\end{equation}
we infer $\norm{[\tau d_{t}u^{k}]'}_{W^{1-1/q,q}}\le\radius$
from~\eqref{eq:ener_subopt}.}
The second estimate of Lemma~\ref{la:vars_tp} and~\eqref{eq:uW} imply that
\[\begin{split}
 \big|\delta^2 \TP(\cdots)[d_t u^k,d_t u^k]\big| 
&\le  c_2 \big(\|[u^{k-1}]'\|_{W^{1-1/q,q}}^{2q+2} + 1\big) \|[d_t u^k]'\|_{W^{1-1/q,q}}^2 \\
&\le  c_2 \big( \br{2c_{q}^{2}\kappa^{-1}e_{0}}^{q+1} +1\big) c_{q}^{2} \|d_t u^k\|_{H^2}^2.
\end{split}\]
Setting $c_b = \vrho c_2 \big( \br{2c_{q}^{2}\kappa^{-1}e_{0}}^{q+1} +1\big) c_{q}^{2}$ and
assuming that $\tau c_b \le\tfrac12$ we absorb the last term in the right-hand
side of~\eqref{eq:ener_taylor} and deduce that
\[
(1-c_b \tau) \|d_t u^k\|_{H^2}^2 + d_t E(u^k) \le 0. 
\]
This implies that the energy decay law also holds with $K$ and thus
proves the assertion with $c'=\korr{\max}\set{c_{a}\korr{e_{0}^{-1}}/2,
2 e_0 c_q^2 \radius^{-2}, c_b }$.  
\end{proof}

\begin{remarks}
\textup{(i)} Note that the radius $\radius$ can be chosen independently of~$k$.
According to Lemma~\ref{lem:2estimates} it only depends on
$\lambda$, $\Lambda$, $M$, $n$, and $q$, so we merely have to
show that there is a uniform bound on the deviation from
arclength parametrization and on the tangent-point energy.
The former can be achieved
from~\eqref{eq:arclength} and~\eqref{eq:radius} 
by additionally claiming $R\le\tfrac12$.
The latter follows from $\TP(u^{k})\le e_{0}\vrho^{-1}$ according to~\eqref{eq:stability}.

\textup{(ii)} 
From $c'\simeq\max\set{e_{0}^{q-1}\k^{-q}\rho^{2},
\sbbb{e_0 R^{-2},} 
\br{e_{0}^{q+1}\k^{-q-1}+1}\rho}$
we infer that we have to choose smaller time steps in case
either $\k^{-1}$ or $\rho$ increases.

\textup{(iii)} 
In case of a metric related to a norm $\|\cdot \|_X$ that defines the gradient 
flow an inequality $\|w'\|_{W^{1-1/q,q}} \le c_q \|w\|_X$ is required which may 
involve an inverse estimate in a fully discrete setting and then implies a more
restrictive step size condition. Discrete norms such as
\[
\|w_h\|_{h,r}^2 = \|w_h\|_{L^2}^2 + h^r \|w_h''\|_{L^2}^2 
\]
with a mesh-size parameter $h>0$ that mimic the scaling properties of fractional 
Sobolev spaces may lead to practical and efficient numerical schemes.

\textup{(iv)} 
Minimal modifications of the stability result are required if other boundary
conditions are considered. Note that $d_t u^k$ satisfies corresponding homogeneous 
boundary conditions which, e.g., in the case of a clamped boundary condition
imply the Poincar\'e inequalities. The inhomogeneous conditions then only enter 
in~\eqref{eq:uW}.

\end{remarks}%

\input{tp-2var}

\input{experiments}

\subsection*{Acknowledgments}

Philipp Reiter was partially supported by DFG-Grant RE~3930/1--1. 
\sbbb{The work on this manuscript was initiiated during the workshop
``Geometric curvature functionals and discretizations'' organized
by Heiko von~der~Mosel which took place in Kloster Steinfeld in September 2017. }

\bibliographystyle{abbrvhref}
\bibliography{bib_project}

\end{document}

%% file: sb_macros.tex
\def\R{\mathbb{R}}

\def\N{\mathbb{N}}

\def\Z{\mathbb{Z}}

\def\cA{\mathcal{A}}

\def\cM{\mathcal{M}}

\def\cO{\mathcal{O}}

\def\d{\delta}

\def\k{\kappa}

\def\p{\partial}

\def\veps{\varepsilon}
\def\vrho{\varrho}

\def\NN{{\rm N}} 
 
\newcommand{\dv}[1]{\,{\mathrm d}#1}

\newcommand{\wcheck}[1]{#1\hspace{-.8ex}\mbox{\huge {\lower.45ex \hbox{$\textstyle \check{}$}}} \hspace{.5ex}}

\let\oldmarginpar\marginpar
\renewcommand\marginpar[1]{
  \oldmarginpar[\raggedleft\footnotesize #1]
  {\raggedright\footnotesize #1}}
\newtheorem{definition}{Definition}
\newtheorem{lemma}[definition]{Lemma}
\newtheorem{proposition}[definition]{Proposition}

\newtheorem{corollary}[definition]{Corollary}
\newtheorem{remark}[definition]{Remark}
\newtheorem{remarks}[definition]{Remarks}

\numberwithin{definition}{section}

%% file: pr_macros.tex
\usepackage{esint}

\hypersetup{
breaklinks=true,
colorlinks=true,
linkcolor=blue,
citecolor=blue,
urlcolor=blue,
}

\newcommand{\abs}[1]{\ensuremath{\left\lvert#1\right\rvert}}
\newcommand{\br}[1]{\ensuremath{\left(#1\right)}}

\newcommand{\linorm}[1]{\ensuremath{\left\|#1'\right\|_{L^{{\infty}}}}}
\newcommand{\norm}[1]{\ensuremath{\left\lVert#1\right\rVert}}

\newcommand{\seqn}[2][k]{\ensuremath{\br{{#2}_{#1}}_{#1\in\N}}}
\newcommand{\set}[1]{\ensuremath{\left\{#1\right\}}}

\newcommand{\sett}[2]{\left\{#1\;\middle|\;#2\right\}}
\renewcommand{\sp}[1]{\ensuremath{\left\langle#1\right\rangle}}

\newcommand{\sq}[1]{\ensuremath{\left[#1\right]}}

\DeclareMathOperator{\bil}{biL}
\newcommand{\CC}[1][]{\ensuremath{\mathscr{C}_{#1}}}

\newcommand\D[1]{{#1(x)-#1(y)}}
\renewcommand{\d}{\ensuremath{\,\mathrm{d}}}

\newcommand{\eps}{\ensuremath{\varepsilon}}

\newcommand{\iintrz}{\iint_{\R/\Z\times\R/\Z}}

\renewcommand{\kappa}{\varkappa}

\renewcommand{\rho}{\ensuremath{\varrho}}
\newcommand\rzd{(\R/\Z,\R^{n})}

\newcommand{\TP}{\mathrm{TP}}



%% file: tp-2var.tex

\section{Tangent-point energies}\label{sec:tp}

We consider the functional
\begin{equation}\label{eq:vartp}
 \TP(u)=\tfrac1q\iintrz\frac{\abs{u'(y)\wedge\br{\D u}}^{q}}
 {\abs{\D u}^{2q}} \d x\d y,
 \qquad q>2,
\end{equation}
where $u:\R/\Z\to\R^{n}$ is a continuously differentiable curve
which is embedded (injective) and {regular} ($\min_{\R/\Z}\abs{u'}>0$).
If we multiply the integrand of~\eqref{eq:vartp} by
$\abs{u'(x)}^{1-q}\abs{u'(y)}$, we obtain the
``classical'' tangent-point functional~\eqref{eq:tp} mentioned
in the introduction, so their difference 
linearly tends to zero as $u$ approaches arclength parametrization, see Proposition~\ref{prop:consist} below.
Corresponding statements for the derivatives of these functionals
can be obtained in a straightforward manner, so we omit them here.

Formula~\eqref{eq:vartp} applies to curves in euclidean spaces of
arbitrary dimension; to this end we point out that only the norm
of the vector product $\wedge$ appears which regardless of dimension can be evaluated
in terms of scalar products via
\begin{equation}\label{eq:wedge}
 \sp{a\wedge b,c\wedge d} = \sp{a,c}\sp{b,d} - \sp{a,d}\sp{b,c}
 \qquad\text{for all }a,b,c,d\in\R^{n}.
\end{equation}


As mentioned in the introduction,
the characteristic feature of the tangent-point energies (and many other knot energies)
is a uniform bound on the \emph{bi-Lipschitz constant} of a curve~$u$, namely
$\sup_{x,y\in\R/\Z}\frac{d_{u}(x,y)}{\abs{u(x)-u(y)}}$
where $d_{u}(x,y)$ denotes the \emph{intrinsic} distance of $u(x)$
and $u(y)$ on the curve~$u$.
In case of arclength parametrization, this term is equivalent
to Gromov's distortion~\cite{gromov},
and we also have $d_{u}(x,y)=\abs{x-y}$ for all $x,y\in\R/\Z$
where $\abs z=\min_{k\in\Z}\abs{z-k}$ for $z\in\R/\Z$.
As we are interested in settings close to that condition, we
replace $d_{u}(x,y)$ by $\abs{x-y}$ and define
\begin{equation}\label{eq:bil}
 \bil(u) = \sup_{\substack{x,y\in\R/\Z\\x\ne y}} \frac{\abs{x-y}}{\abs{\D u}}\ge \max_{x\in\R/\Z}\frac1{\abs{u'(x)}}
 = \frac1{\min_{x\in\R/\Z}\abs{u'(x)}}.
\end{equation}

In order to discuss some fundamental facts about $\TP$,
we define
\[ \CC = \sett{u\in W^{2-1/q,q}\rzd}{u\text{ is embedded and regular}}. \]
Here, for $s\in(0,1)$,
\[ \sq f_{W^{s,q}} = \br{\iint_{\R/\Z\times\R/\Z}{\frac{\abs{{\D f}}^{q}}{\abs{x-y}^{1+sq}}}\d x\d y}^{1/q} \]
denotes the Sobolev--Slobodecki\u\i\ seminorm.
The corresponding spaces are given by
\[ W^{k+s,q}=\sett{f\in L^{q}}{\sq f_{W^{k+s,q}}<\infty}, \qquad k\in\N\cup\set0. \]
We will see in Section~\ref{sec:prelim} below that
$\CC$ is contained in the energy space for $\TP$,
i.e., $u\in\CC\implies\TP(u)<\infty$.
Note that, as the domain $\R/\Z$ is one-dimensional,
we have the embedding $W^{2-1/q,q}\hookrightarrow C^{1,1-2/q}$
for $q>2$ which shows that the initial assumption of $C^{1}$ curves
does not imply any additional restriction.

The case $q=2$ is excluded for several reasons.
It corresponds to the analytically challenging scale invariant case.
While the image of a finite-energy curve
is a topological manifold~\cite[Thm.~1.1]{StrMos12}
it is not clear to what extent the other results for $q>2$ also apply to $q=2$.
However the discrete counterpart
of $\TP$ is just too weak for modeling self-avoidance if $q=2$~\cite[Sect.~4.2]{BRR},
so we will rather choose large values of~$q$.
On the other hand, discretization estimates impose
$q<4$~\cite[Lem.~3.1]{BRR}.

\subsection{Preliminaries}\label{sec:prelim}

We start with some general facts on $\TP$ and its relation
to the ``classical'' tangent-point functional~\eqref{eq:tp}.
However, the stability result in Section~\ref{sec:stab}
does not rely on the following statements.

The crucial observation for analytical investigation of $\TP$
is that, using $u'(y)\wedge u'(y)=0$, we may write
\begin{equation}\tag{\ref{eq:vartp}*}\label{eq:vartp-}
 \TP(u)=\tfrac1q\iintrz\frac{\abs{u'(y)\wedge\br{\D u-(x-y)u'(y)}}^{q}}
 {\abs{\D u}^{2q}} \d x\d y.
\end{equation}

\begin{proposition}[Energy estimate]\label{prop:0th}
 The functional $\TP$ is continuous on~$\CC$
 and satisfies
 \begin{equation}\label{eq:0th}
 \abs{\TP(u)} \le \tfrac1q\br{\bil u}^{2q}\norm{u'}_{L^\infty}^{q}\sq{u'}_{W^{1-1/q,q}}^{q} \qquad\text{for any }u\in\CC.
 \end{equation}
\end{proposition}

Continuity can be shown by
bounding the difference $\TP(\tilde u)-\TP(u)$
for $\tilde u,u\in\CC$ in terms of $\sq{\br{\tilde u-u}'}_{W^{2-1/q,q}}$
using Corollary~\ref{cor:bil} below.
However, it is sufficient to do this for the first or second derivative
and then conclude by integrating
as in~\cite[Rem.~3.1]{BR15} for the first variation of~\eqref{eq:tp}.
We omit the details.
As its arguments will be repeatedly applied throughout
this section, we provide a proof of Blatt's energy estimate~\cite{Blatt13} for the reader's convenience.

\begin{proof}[Proof of~\eqref{eq:0th}]
 We start by noting that \sbb{due to~\eqref{eq:vartp-} and~\eqref{eq:wedge} we have}
 \begin{align*}
 &0\le\frac{q\TP(u)}{\br{\bil u}^{2q}}
 \le\iintrz\frac{\abs{u'(y)\wedge\br{\D u}}^{q}}
 {\abs{x-y}^{2q}} \d x\d y \\
 &=\iintrz\frac{\abs{u'(y)\wedge(x-y)\br{\int_{0}^{1}u'(y+\theta(x-y))\d\theta-u'(y)}}^{q}}
 {\abs{x-y}^{2q}} \d x\d y \\
 &\le\norm{u'}_{L^\infty}^{q}\iintrz\frac{\abs{\int_{0}^{1}u'(y+\theta(x-y))\d\theta-u'(y)}^{q}}{\abs{x-y}^{q}} \d x\d y \\
 &\le\norm{u'}_{L^\infty}^{q}\iintrz\int_{0}^{1}\abs{\frac{{u'(y+\theta(x-y))-u'(y)}}{{x-y}}}^{q}\d\theta \d x\d y, 
 \end{align*}
 where we invoked Jensen's inequality in the last step.
 Applying Fubini's theorem, we arrive at
 \begin{align*}
 \frac{q\TP(u)}{\br{\bil u}^{2q}\norm{u'}_{L^\infty}^{q}}
 &= \int_{0}^{1}\iintrz\abs{\frac{{u'(y+\theta(x-y))-u'(y)}}{{x-y}}}^{q}\d x\d y\d\theta \\
 &= \int_{0}^{1}\int_{\R/\Z}\int_{-1/2}^{1/2}\abs{\frac{{u'(y+\theta z)-u'(y)}}{z}}^{q}\d z\d y\d\theta \\
 &= \int_{0}^{1}\int_{\R/\Z}\int_{-\theta/2}^{\theta/2}\abs{\frac{{u'(y+\tilde z)-u'(y)}}{\tilde z}}^{q}\theta^{q-1}\d\tilde z\d y\d\theta \\
 &\le\int_{\R/\Z}\int_{-1/2}^{1/2}\abs{\frac{{u'(y+\tilde z)-u'(y)}}{\tilde z}}^{q}\d\tilde z\d y
 = \sq{u'}_{W^{1-1/q,q}}^{q}
 \end{align*}
 which is the desired estimate.
\end{proof}


Working in the class of $H^{2}$ curves (which is contained in
$W^{2-1/q,q}$ provided $q\le4$) in Section~\ref{sec:stab},
a characterization of the energy spaces of $\TP$ is not required.
However, we briefly state the following.

\begin{remark}[Energy spaces]
 Revisiting the proof of~\eqref{eq:0th} and replacing the term
 $\br{\int_{0}^{1}u'(y+\theta(x-y))\d\theta-u'(y)}$ by
 \[ {\int_{0}^{1}\br{u'(y+\theta(x-y))-\frac{u'(y)}{\abs{u'(y)}}\abs{u'(y+\theta(x-y))}}\d\theta}, \]
 we arrive at
 \begin{equation}\label{eq:0th*}\tag{\ref{eq:0th}*}
  \abs{\TP(u)} \le \tfrac1q\br{\bil u}^{2q}\norm{u'}_{L^\infty}^{2q}\sq{\frac{u'}{\abs{u'}}}_{W^{1-1/q,q}}^{q} \qquad\text{for any }u\in\CC.
 \end{equation}
 Together with the subsequent estimate~\eqref{eq:energ-bound} obtain
 the following characterization of the energy spaces.
 \emph{If $u$ is an embedded regular $C^{1}$ curve then
 \[ \TP(u)<\infty\iff\frac{u'}{\abs{u'}}\in W^{1-1/q,q}\rzd. \]}%
 Some further comments are in order.
 \setlength{\leftmargini}{0em}
 \begin{enumerate}
 \setlength{\itemindent}{3em}
 \item Strzelecki and von der Mosel~\cite[Thm.~1.3]{StrMos12} have
 shown that any curve for which~\eqref{eq:tp} is finite is
 a $C^{1,1-2/q}$ manifold.
 Adapting the arguments in~\cite{StrMos12}
 one may relax the condition $u\in\CC$
 to $u\in W^{1,1}$ with positive lower and finite upper bound
 on $\abs{u'}$.
 \item As pointed out in~\cite[p.~2]{StrMos12},
 a finite $\TP$ value does not imply injectivity.
 The curve $u$ could cover a path several times or even
 form a manifold with boundary.
 We have to exclude these phenomena by claiming embeddedness
 in order to state bi-Lipschitz estimates.
 \item The condition $u\in W^{2-1/q,q}$ implies $\frac{u'}{\abs{u'}}\in W^{1-1/q,q}$,
 but the converse is not true.
 Note that an open embedded regular $C^{1}$ curve $u$ lying on a straight line has tangent-point
 energy zero while we can certainly produce a parametrization
 where $\sq{u'}_{W^{1-1/q,q}}>0$.
 However, as its unit tangent is constant, we infer
 $\sq{\frac{u'}{\abs{u'}}}_{W^{1-1/q,q}}=0$ in accordance with~\eqref{eq:energ-bound} below.
 \end{enumerate}
 For the sake of simplicity, we will keep on restricting to $W^{2-1/q,q}$ curves.
\end{remark}

Now we proceed to the converse estimate to~\eqref{eq:0th*}
which is in fact used for the proof of Proposition~\ref{prop:unibi} later on.

\begin{proposition}[Necessary regularity for finite energy]
 Let $u\in C^{1}(\R/\Z,\R^n)$ be embedded and 
 regular with $\TP(u)<\infty$.
 Then $\frac{u'}{\abs{u'}}\in W^{1-1/q,q}$ and
 \begin{equation}\label{eq:energ-bound}
  \sq{\frac{u'}{\abs{u'}}}_{W^{1-1/q,q}}^q 
  \le C\br{\TP(u)+ \TP(u)^{\frac{q-1}{q-2}}}
 \end{equation}
 where $C$ only depends on $\norm{u'}_{L^{\infty}}$, $\min_{\R/\Z}\abs{u'}$, and $q$.
\end{proposition}

The proof follows by minor modifications from the arguments given in~\cite[Prop.~2.5]{BR15}.

Although we will not further rely on it, we briefly note that the functionals~\eqref{eq:tp}
and~\eqref{eq:vartp} are closely related
and in fact agree as $u$ approaches arc\-length parametrization.

\begin{proposition}[Approximation]\label{prop:consist}
 For any $u\in\CC$ and $\Lambda\ge1$ with
 \begin{equation*}
 \frac1\Lambda\le\abs{u'(x)}\le\Lambda\qquad\text{for all }x\in\R/\Z,
 \end{equation*}
 the difference between the ``classical''
 tangent-point functional~\eqref{eq:tp}
 and $\TP$ defined in~\eqref{eq:vartp} is bounded above by
 \[ \br{\bil u}^{2q}\Lambda^{2q}\sq{u'}_{W^{1-1/q,q}}^{q}(\Lambda-1). \]
\end{proposition}

Corresponding estimates for the derivatives of~$\TP$ can be derived accordingly.

\begin{proof}
 The ``classical'' tangent-point functional is given by
 \begin{equation}\label{eq:tp-expl}\tag{\ref{eq:tp}*}
 u\mapsto\tfrac1q\iintrz
 \frac{\abs{P_{u'(y)}^{\perp}\br{{\D u}}}^{q}}
 {\abs{\D u}^{2q}}\abs{u'(x)}\abs{u'(y)}\d x\d y
 \end{equation}
 where $P_{u'(y)}^{\perp}$ denotes the projection onto the normal
space to $\R u'(y)$, cf.~\cite{BR15}.
 From
\[ \begin{split}
  \abs{P_{u'(y)}^{\perp}\br{{\D u}}}^{2}
  &=\abs{\D u - \sp{\D u,\frac{u'(y)}{\abs{u'(y)}}}\frac{u'(y)}{\abs{u'(y)}}}^{2} \\
  &=\abs{\D u}^{2} - \abs{\sp{\D u,\frac{u'(y)}{\abs{u'(y)}}}}^{2} \\
  &=\abs{u'(y)}^{-2}\abs{{u'(y)}\wedge{\D u}}^{2}
 \end{split}\]
 we infer that the integrands of~\eqref{eq:tp-expl} and~\eqref{eq:vartp}
 just differ by a factor of $\abs{u'(x)}^{1-q}\abs{u'(y)}$.
 Applying the estimate $\abs{a^{p}-b^{p}}\le p\max(a,b)^{p-1}\abs{a-b}$
 which holds for all $a,b\ge0$, $p\ge1$,
 as well as the techniques employed for the proof of
 Proposition~\ref{prop:0th}, we arrive at the assertion.
%
\end{proof}

For discretization aspects, involving an estimate
on the contribution of the diagonal of $\R/\Z\times\R/\Z$
we refer to~\cite{BRR}.

\subsection{Uniform bi-Lipschitz continuity}

As repeatedly pointed out, a uniform bi-Lipschitz bound
in terms of the energy is the essential property of~$\TP$.
\sbbb{Using~Proposition~\ref{prop:consist},}
the following statement follows by the same lines as
the corresponding proof given in~\cite[Prop.~2.7]{BR15}.

\begin{proposition}[Uniform bi-Lipschitz estimate]\label{prop:unibi}
 For any $\lambda,\Lambda,M>0$ there is a constant $C=C(\lambda,\Lambda,M,n,q)>0$
 such that any curve $u\in\CC$ with
 \[ \lambda\le\abs{u'}\le\Lambda \qquad\text{and}\qquad \TP(u)\le M \]
 satisfies $\bil u\le C$, more precisely,
  \[ \abs{x-y}_{\R/\Z}\le C\abs{u(x)-u(y)}_{\R^{n}}
 \qquad\text{for all }x,y\in\R/\Z. \]
\end{proposition}

Before discussing the first variation of $\TP$,
we will identify a radius~$\radius$ about $u$ in \korr{$W^{1,\infty}$} in terms of $\TP(u)$ ensuring a bi-Lipschitz constant
of $2\bil u$. To this end, we first need the following characterization of bi-Lipschitz curves.

\begin{lemma}[Bi-Lipschitz radius]\label{lem:bil}
 A curve $u\in C^{1}(\R/\Z,\R^{n})$ is embedded and regular if and only if
 \[ \bil u < \infty. \]
 Moreover, the inequality
 \[ \bil\tilde u \le 2\bil u \]
 holds for all $\tilde u\in C^{1}(\R/\Z,\R^{n})$ provided
 \begin{equation}\label{eq:u-u}
 \norm{\br{\tilde u-u}'}_{L^{\infty}} \le \frac1{2\bil u}.
 \end{equation}
\end{lemma}

Consequently, if $u\in C^{1}\rzd$ is bi-Lipschitz continuous,
all curves in \sbbb{a ball of radius $(2\bil u)^{-1}$ around~$u$}
 are bi-Lipschitz continuous as well.

\begin{proof}
 From~\eqref{eq:bil} we read off that $\bil u<\infty$
 implies injectivity of $u$ (which gives embeddedness)
 as well as $\min_{\R/\Z}\abs{u'}>0$, i.e., $u$ is regular.
 If $\bil u=\infty$ there are sequences $\seqn x, \seqn y\subset\R/\Z$ such that
 \[ \frac{x_{k}-y_{k}}{\abs{u(x_{k})-u(y_{k})}} \xrightarrow{k\to\infty}\infty. \]
 By compactness we may assume $x_{k}\to x_{0}$, $y_{k}\to y_{0}$.
 As the nominator is bounded, we infer $u(x_{k})-u(y_{k})\to0$
 as $k\to\infty$, so $u(x_{0})=u(y_{0})$ by continuity.
 Either $u$ is not embedded or $x_{0}= y_{0}$.
 The latter gives
 \[ 0\xleftarrow{k\to\infty}\abs{\frac{u(x_{k})-u(y_{k})}{x_{k}-y_{k}}}
 =\abs{\fint_{x_{k}}^{y_{k}}u'(\xi)\d\xi}\xrightarrow{k\to\infty}\abs{u'(x_{0})}, \]
 so $u$ is not regular.
 
 Now let $\bil u<\infty$ and assume that~\eqref{eq:u-u} holds.
 Abbreviating $w=\tilde u-u$ we have
 for all $x,y\in\R/\Z$, $x\ne y$,
 \begin{align*}
 &\abs{\frac{\tilde u(x)-\tilde u(y)}{x-y}}^{2}
 = \abs{\fint_{x}^{y}\tilde u'(\xi)\d\xi}^{2} \\
 &=\fint_{x}^{y}\fint_{x}^{y}\br{u+w}'(\xi)\br{u+w}'(\eta)\d\xi\d\eta \\
 &=\fint_{x}^{y}\fint_{x}^{y}\sq{u'(\xi)u'(\eta)+w'(\xi)w'(\eta)+u'(\xi)w'(\eta)+u'(\eta)w'(\xi)}\d\xi\d\eta \\
 &\ge\abs{\fint_{x}^{y}u'(\xi)\d\xi}^{2} + \abs{\fint_{x}^{y}w'(\xi)\d\xi}^{2} - 2\abs{\fint_{x}^{y}u'(\xi)\d\xi\fint_{x}^{y}w'(\eta)\d\eta} \\
 &\ge\tfrac12\abs{\fint_{x}^{y}u'(\xi)\d\xi}^{2} - \abs{\fint_{x}^{y}w'(\xi)\d\xi}^{2} \\
 &\ge\frac1{2\br{\bil u}^{2}} - \norm{w'}_{L^{\infty}}^{2}
 \ge \frac1{\br{2\bil u}^{2}},
 \end{align*}
\sbbb{which proves the assertion.}
\end{proof}

\begin{corollary}\label{cor:bil}
 For any $\lambda,\Lambda,M>0$ there is a radius $\radius=\radius(\lambda,\Lambda,M,n,q)>0$
 such that any curve $u\in\CC$ with
 \[ \lambda\le\abs{u'}\le\Lambda \qquad\text{and}\qquad \TP(u)\le M \]
 and any $\tilde u\in W^{2-1/q,q}(\R/\Z,\R^{n})$ with
 \[ \norm{\br{\tilde u-u}'}_{L^{\infty}} \le \radius \]
 satisfy
 \[ \bil\tilde u\le2\bil u<\infty. \]
\end{corollary}

\begin{proof}
 The assumption on $u$ implies $\TP(u)<\infty$. According to
 Proposition~\ref{prop:unibi} there is some $C=C\br{\min_{\R/\Z}\abs{u'}, \norm{u'}_{L^{\infty}},\TP(u),n,q}>0$
 with $\bil u\le C$.
 Choosing $\radius = \frac1{2C} \le \frac1{2\bil u}$, Lemma~\ref{lem:bil}
 yields the desired \sbbb{estimate.}
\end{proof}

\subsection{First derivative}

The first variation of $\TP$ as well as its discretization
\sbbb{have} already been derived in~\cite{BRR}. Its formula reads
\setlength{\multlinegap}{0pt}
\newcommand{\MM}[3][u]{\mathcal M(#1;#2,#3)}
\renewcommand{\AA}[3][u]{\mathcal A(#1;#2,#3)}
\renewcommand{\NN}[7][u]{\mathcal N(#1;#2,#3,#4,#5,#6,#7)}
\newcommand{\PP}[4][u]{\mathcal P(#1;#2,#3,#4)}
\newcommand{\BB}[7][u]{\mathcal B(#1;#2,#3,#4,#5,#6,#7)}
\newcommand{\XX}[9]{\mathcal X(u;#1,#2,#3,#4,#5,#6,#7,#8,#9)}
\newcommand{\W}{W^{2-1/q,q}}
\begin{equation}\label{eq:firstvar}
 \delta\TP(u)[\varphi] = \MM u\varphi + \MM \varphi u - 2\AA u\varphi
\end{equation}
where
\begin{align*}
\MM vw &= \iintrz
\frac{\abs{u'(y)\wedge\br{\D u}}^{q-2}}{\abs{\D u}^{2q}} \cdot \\
&\qquad\qquad{}\cdot\Big\langle u'(y)\wedge\br{\D u},\cdot \\
&\qquad\qquad\qquad{}\cdot v'(y)\wedge\br{\D w-(x-y)w'(y)}\Big\rangle\d x\d y, \\
\AA vw &= \iintrz
\frac{\abs{u'(y)\wedge\br{\D u}}^{q}}{\abs{\D u}^{2q+2}}\cdot \\
&\qquad\qquad\qquad\qquad{}\cdot\sp{\D v,\D w}\d x\d y.
\end{align*}
Note that for the implementation of the algorithm we may omit
the $-(x-y)w'(y)$ term in $\MM vw$ for it cancels in symmetric expressions
$\MM vw + \MM wv$ due to $v'\wedge w'+w'\wedge v'=0$.

The first variation formula is considerably simpler than the corresponding one in~\cite[(1.11), Rem.~3.1]{BR15}
that has been derived for the parametrization invariant functional~\eqref{eq:tp}.
If $u$ is parametrized by arclength and $u'\perp \varphi'$
both formulae agree.

\begin{proposition}\label{prop:1st}
 For any $u\in\CC$
 the functional $\TP$ is continuously differentiable.
 In particular, its first variation defines a bounded linear form,
 and we have for any $w\in\W\rzd$
 \begin{equation}\label{eq:1th}
  \abs{\delta\TP(u)[w]} \le C\br{\bil u}^{2q+2}\norm {u'}_{L^{\infty}}^{q+1}\norm{u'}_{W^{1-1/q,q}}^{q}\cdot\norm{w'}_{W^{1-1/q,q}}.
 \end{equation}
 The constant $C>0$ only depends on $n$ and $q$.
\end{proposition}

The existence of the first variation~\eqref{eq:firstvar} can be
rigorously derived by
establishing a uniform $L^{1}$ bound for the difference quotient
of the integrand of $\tfrac1\tau\br{\TP(u+\tau w)-\TP(u)}$,
cf.~\cite[Sect.~3]{BR15}.
As for Proposition~\ref{prop:0th}, we will omit the proof of
continuity.

\begin{proof}[Proof of \eqref{eq:1th}]
 We obtain
 \begin{align*}
 &\abs{\mathcal M(u;v,w)} \\
 &\le\iint_{\R/\Z\times\R/\Z}\frac{\abs{u'(y)\wedge\br{\D u-(x-y)u'(y)}}^{q-1}}{\abs{\D u}^{2q}}\cdot{} \\
 &\qquad\qquad\qquad\qquad{}\cdot\abs{v'(y)\wedge\br{\D w-(x-y)w'(y)}}\d x\d y \\
 &\le\br{\bil u}^{2q}\norm{u'}_{L^{\infty}}^{q-1}\norm{v'}_{L^{\infty}}
 \iint_{\R/\Z\times\R/\Z}\br{\frac{\abs{{\D u-(x-y)u'(y)}}}{\br{x-y}^{2}}}^{q-1}\cdot{}\\
 &\qquad\qquad\qquad\qquad\qquad\qquad\qquad\qquad{}\cdot\frac{\abs{\D w-(x-y)w'(y)}}{(x-y)^{2}}\d x\d y \\
 &\le\br{\bil u}^{2q}\norm{u'}_{L^{\infty}}^{q-1}\norm{v'}_{L^{\infty}} \cdot \\
 &\qquad\qquad\qquad{}\cdot\br{\iint\limits_{\R/\Z\times\R/\Z}{\frac{\abs{{\D u-(x-y)u'(y)}}^{q}}{\abs{x-y}^{2q}}}\d x\d y}^{\frac{q-1}q}\cdot{}\\
 &\qquad\qquad\qquad{}\cdot
 \br{\iint\limits_{\R/\Z\times\R/\Z}{\frac{\abs{{\D w-(x-y)w'(y)}}^{q}}{\abs{x-y}^{2q}}}\d x\d y}^{\frac1q} \\
 &\le \br{\bil u}^{2q}\norm{u'}_{L^{\infty}}^{q-1}
 \sq{u'}_{W^{1-1/q,q}}^{q-1} \norm{v'}_{L^{\infty}} \sq{w'}_{W^{1-1/q,q}} \\
 &\le \br{\bil u}^{2q+2}\norm{u'}_{L^{\infty}}^{q+1}
 \sq{u'}_{W^{1-1/q,q}}^{q-1} \norm{v'}_{L^{\infty}} \sq{w'}_{W^{1-1/q,q}}
 \end{align*}
 where we used the estimates employed for the proof of~\eqref{eq:0th}.
 In the same way, we derive
 \[ \abs{\mathcal A(u;v,w)} \le \br{\bil u}^{2q+2}\norm{u'}_{L^{\infty}}^{q}\sq{u'}_{W^{1-1/q,q}}^{q}\norm{v'}_{L^{\infty}}\norm{w'}_{L^{\infty}}. \]
 The assertion follows by the embedding $W^{1-1/q,q}\hookrightarrow L^{\infty}$.
\end{proof}

\subsection{Second derivative}

In order to compute the second variation, we introduce
\begin{align*}
&\NN vw\phi\psi\xi\eta\\
&= \iintrz
\frac{\abs{u'(y)\wedge\br{\D u}}^{q-4}}{\abs{\D u}^{2q}} \cdot \\
&\qquad\qquad{}\cdot\sp{u'(y)\wedge\br{\D u},v'(y)\wedge\br{\D w}}\cdot\\
&\qquad\qquad{}\cdot\sp{\phi'(y)\wedge\br{\D\psi},\xi'(y)\wedge\br{\D\eta}}\d x\d y, \\
&\PP vw\phi\\
&= \iintrz
\frac{\abs{u'(y)\wedge\br{\D u}}^{q-2}}{\abs{\D u}^{2q+2}} \cdot \\
&\qquad\qquad{}\cdot\sp{u'(y)\wedge\br{\D u},v'(y)\wedge\br{\D w}}\cdot\\
&\qquad\qquad{}\cdot\sp{\D u,\D\phi}\d x\d y, \\
&\BB vw\phi\psi\xi\eta\\
&= \iintrz
\frac{\abs{u'(y)\wedge\D u}^{q-2}}{\abs{\D u}^{2q+4}}\sp{\D v,\D w}\cdot\\
&\qquad\qquad{}\cdot\sp{u'(y)\wedge\D u,\phi'(y)\wedge\D\psi}\cdot\\
&\qquad\qquad{}\cdot\sp{\D\xi,\D\eta}\d x\d y.
\end{align*}

Computing the first variation of $\MM\cdot\cdot$ and $\AA\cdot\cdot$
we arrive at
\begin{align*}
&\delta\mathcal M(u;v,w;h) = \left.\tfrac\d{\d\tau}\right|_{\tau=0}\MM[u+\tau h]vw \\
&=(q-2)\iintrz
\frac{\abs{u'(y)\wedge\br{\D u}}^{q-4}}{\abs{\D u}^{2q}} \cdot \\
&\qquad\qquad{}\cdot\sp{u'(y)\wedge\br{\D u},v'(y)\wedge\br{\D w}} \cdot\\
&\qquad\qquad{}\cdot\Big\langle u'(y)\wedge\br{\D u}, \cdots{} \\
&\qquad\qquad\qquad{}\cdots h'(y)\wedge\br{\D u}+u'(y)\wedge\br{\D h}\Big\rangle \d x\d y \\
&\quad{}-2q\iintrz
\frac{\abs{u'(y)\wedge\br{\D u}}^{q-2}}{\abs{\D u}^{2q+2}} \cdot \\
&\qquad\qquad{}\cdot\sp{u'(y)\wedge\br{\D u},v'(y)\wedge\br{\D w}}\cdot\\
&\qquad\qquad{}\cdot\sp{\D u,\D h}\d x\d y \\
&\quad{}+
\iintrz
\frac{\abs{u'(y)\wedge\br{\D u}}^{q-2}}{\abs{\D u}^{2q}} \cdot \\
&\qquad\qquad{}\cdot\Big\langle h'(y)\wedge\br{\D u}+u'(y)\wedge\br{\D h}, \cdots{} \\
&\qquad\qquad{}\cdots v'(y)\wedge\br{\D w}\Big\rangle\d x\d y \\
& = (q-2)\br{\NN vwuuhu + \NN vwuuuh} - 2q\PP vwh + {}\\
&\quad{}+ \NN uuhuvw + \NN uuuhvw
\end{align*}
and
\begin{align*}
&\delta\mathcal A(u;v,w;h) = \left.\tfrac\d{\d\tau}\right|_{\tau=0}\AA[u+\tau h]vw \\
&=q\iintrz
\frac{\abs{u'(y)\wedge\br{\D u}}^{q-2}}{\abs{\D u}^{2q+2}} \cdot \\
&\qquad\qquad{}\cdot\sp{\D v,\D w} \cdot\\
&\qquad\qquad{}\cdot\Big\langle u'(y)\wedge\br{\D u},\cdots{} \\
&\qquad\qquad\qquad{}\cdots h'(y)\wedge\br{\D u}+u'(y)\wedge\br{\D h}\Big\rangle \d x\d y \\
&\quad{}-(2q+2)\iintrz
\frac{\abs{u'(y)\wedge\br{\D u}}^{q}}{\abs{\D u}^{2q+4}} \cdot \\
&\qquad\qquad{}\cdot\sp{\D v,\D w} \cdot\\
&\qquad\qquad{}\cdot\sp{\D u,\D h} \d x\d y \\
&=q\br{\BB vwhuuu+\BB vwuhuu}-2(q+1)\BB vwuuuh.
\end{align*}

\korr{%
Using the relation
\begin{equation}\tag{\ref{eq:firstvar}}
 \delta\TP(u)[\varphi] = \MM u\varphi + \MM \varphi u - 2\AA u\varphi
\end{equation}
and that $\cM$ and $\cA$ are linear in
their second and third arguments, we arrive at}%
\begin{align*}
&\delta^{2}\TP(u)[v,w] \\
&=\delta\mathcal M(u;u,v;w) + \MM vw +\delta\mathcal M(u;v,u;w) + \MM wv \\
&\quad{}-2\delta\mathcal A(u;u,v;w) - 2\AA vw \\
&= (q-2)\Big( \NN uvuuwu + \NN uvuuuw + {} \\
&\quad{}+ \NN vuuuwu + \NN vuuuuw \Big) + \NN uuwuuv + {} \\
&\quad{}+ \NN uuuwuv + \NN uuwuvu + \NN uuuwvu - {} \\
&\quad{}- 2q\Big(\PP uvw + \PP vuw\Big) + \MM vw + \MM wv - {} \\
&\quad{}- 2q\Big( \BB uvwuuu + \BB uvuwuu \Big) + {} \\
&\quad{}+ 4(q+1) \BB uvuuuw - 2\AA vw.
\end{align*}

\begin{proposition}\label{prop:2nd}
 For any embedded curve $u\in\W\rzd$
 the functional $\TP$ is twice continuously differentiable.
 In particular, its second variation defines a bounded bilinear form,
 and we have for all $v,w\in\W\rzd$
 \begin{equation}\label{eq:2nd}
  \abs{\delta^{2}\TP(u)[v,w]} \le C\br{\bil u}^{2q+4}\norm {u'}_{L^{\infty}}^{q+2}\norm{u'}_{W^{1-1/q,q}}^{q}\cdot\norm{v'}_{W^{1-1/q,q}}\norm{w'}_{W^{1-1/q,q}}.
 \end{equation}
\end{proposition}

A similar statement can be derived for the parametrization invariant
tangent-point functional.

\begin{proof}[Proof of \eqref{eq:2nd}]
It turns out that all functionals $\mathcal{M,A,N,P,B}$ can be
brought into a common form, namely
\begin{align*}
&\XX vw\phi\psi\xi\eta\zeta\theta\iota \\
&= \iintrz
\frac{\abs{u'(y)\wedge\br{\D u}}^{q-4}}{\abs{\D u}^{2q+4}} \cdot \\
&\qquad\qquad{}\cdot\sp{u'(y)\wedge\br{\D u},v'(y)\wedge\br{\D w}}\cdot\\
&\qquad\qquad{}\cdot\sp{\phi'(y)\wedge\br{\D\psi},\xi'(y)\wedge\br{\D\eta}} \cdot \\
&\qquad\qquad{}\cdot\sp{\D u,\D \zeta}\sp{\D \theta,\D \iota}\d x\d y.
\end{align*}
To be more precise, we have
\begin{align*}
 \MM vw &= \XX vwuuuuuuu, \\
 \AA vw &= \XX uuuuuuuvw, \\
 \NN vw\phi\psi\xi\eta &= \XX vw\phi\psi\xi\eta uuu, \\
 \PP vw\phi &= \XX vw uuu u\phi uu, \\
 \BB vw\phi\psi\xi\eta &= \XX \phi\psi uuuu u\xi\eta.
\end{align*}
We may replace all factors of the form
$v'(y)\wedge\br{\D w}$ by
\[ v'(y)\wedge\br{\D w-(x-y)w'(y)}. \]
To this end, we infer
\[ u'(y)\wedge\br{\D u} = u'(y)\wedge\br{\D u-(x-y)u'(y)} \]
from $a\wedge a=0$ in the sense of~\eqref{eq:wedge}.
Recalling
\[ \sp{u'(y),v'(y)} = \sp{u'(y),w'(y)} = 0 \]
we obtain
\begin{align*}
 &\sp{u'(y)\wedge\br{\D u},v'(y)\wedge\br{\D w}} \\
 &=\sp{u'(y)\wedge\br{\D u-(x-y)u'(y)},v'(y)\wedge\br{\D w}} \\
 &=\Big\langle u'(y)\wedge\br{\D u-(x-y)u'(y)},\cdots \\
 &\qquad{}\cdots v'(y)\wedge\br{\D w-(x-y)w'(y)}\Big\rangle.
\end{align*}
There are no further cases appearing in the formula for $\delta^{2}\TP$
except for the $\mathcal N$ terms.
Here we observe that the additional terms for
\[ \NN uvuuwu + \NN uvuuuw \]
cancel due to
\[ \sp{u'(y)\wedge\br{\D u},w'(y)\wedge\br{\D u} + u'(y)\wedge\br{\D w}} = 0. \]
A similar reasoning applies to
$\NN uuwuuv + \NN uuuwuv + \NN uuwuvu + \NN uuuwvu$.
Here one may expand each term, e.g.,
\begin{align*}
 &\sp{w'(y)\wedge\br{\D u-(x-y)u'(y)},u'(y)\wedge\br{\D v-(x-y)v'(y)}} \\
 &\quad{}-\sp{w'(y)\wedge\br{\D u},u'(y)\wedge\br{\D v}} \\
 &=-(x-y)\sp{w'(y)\wedge\br{\D u},u'(y)\wedge\br{v'(y)}} \\
 &\quad{}-(x-y)\sp{w'(y)\wedge\br{u'(y)},u'(y)\wedge\br{\D v}} \\
 &\quad{}+(x-y)^{2}\sp{w'(y)\wedge{u'(y)},u'(y)\wedge{v'(y)}}
\end{align*}
and see that the right-hand sides sum up to zero.

Therefore we may consider
\begin{align*}
&\abs{\XX vw\phi\psi\xi\eta\zeta\theta\iota} \\
&= \Bigg|\;\iintrz
\frac{\abs{u'(y)\wedge\br{\D u-(x-y)u'(y)}}^{q-4}}{\abs{\D u}^{2q+4}} \cdot \\
&\qquad\qquad{}\cdot\Big\langle u'(y)\wedge\br{\D u-(x-y)u'(y)},\cdots \\
&\qquad\qquad\qquad{}\cdots v'(y)\wedge\br{\D w-(x-y)w'(y)}\Big\rangle\cdot\\
&\qquad\qquad{}\cdot\Big\langle\phi'(y)\wedge\br{\D\psi-(x-y)\psi'(y)},\cdots \\
&\qquad\qquad\qquad{}\cdots\xi'(y)\wedge\br{\D\eta-(x-y)\eta'(y)}\Big\rangle \cdot \\
&\qquad\qquad{}\cdot\sp{\D u,\D \zeta}\sp{\D \theta,\D \iota}\d x\d y\Bigg| \\
&\le \linorm u^{q-2}\linorm v\linorm\phi\linorm\xi\linorm\zeta\linorm\theta\linorm\iota\cdot \\
&\qquad{}\cdot\iintrz
\frac{\abs{{\D u-(x-y)u'(y)}}^{q-4}}{\abs{\D u}^{2q+4}}\abs{x-y}^{4} \cdot \\
&\qquad\qquad{}\cdot\abs{{\D u-(x-y)u'(y)}}\;\abs{{\D w-(x-y)w'(y)}}\cdot\\
&\qquad\qquad{}\cdot\abs{{\D\psi-(x-y)\psi'(y)}}\;\abs{{\D\eta-(x-y)\eta'(y)}} \d x\d y.
\end{align*}
Using Proposition~\ref{prop:unibi} we arrive at
\begin{align*}
&\abs{\XX vw\phi\psi\xi\eta\zeta\theta\iota} \\
&\le \linorm u^{q-2}\linorm v\linorm\phi\linorm\xi\linorm\zeta\linorm\theta\linorm\iota\bil(u)^{2q+4}\cdot \\
&{}\cdot\iintrz
\frac{\abs{{\D u-(x-y)u'(y)}}^{q-3}}{\abs{x-y}^{2q}} \;\abs{{\D w-(x-y)w'(y)}}\cdot \\
&\qquad\qquad{}\cdot\abs{{\D\psi-(x-y)\psi'(y)}}\;\abs{{\D\eta-(x-y)\eta'(y)}} \d x\d y \\
&\le \linorm u^{q-2}\linorm v\linorm\phi\linorm\xi\linorm\zeta\linorm\theta\linorm\iota\bil(u)^{2q+4}\cdot \\
&\qquad{}\cdot\norm u_{\W}^{q-3}\norm w_{\W}\norm\psi_{\W}\norm\eta_{\W}.
\end{align*}
Note that at least one \korr{of} the variables $w$, $\psi$, and $\eta$ coincides
with $u$ in the formula of $\TP$, so H\"older's inequality in the
last step does not require $q>3$.
The statement follows from the embedding theorems for
Sobolev spaces.
The bi-Lipschitz constant only depends on $\TP(u)$.
\end{proof}

\subsection{Higher derivatives}

In a similar fashion, we derive for any $k\in\N$ the general estimate
 \begin{equation}\label{eq:kth}
  \abs{\delta^{k}\TP(u)[v_{1},\dots,v_{k}]} \le C\br{\bil u}^{2q+2k}\norm {u'}_{L^{\infty}}^{q+k}\norm{u'}_{W^{1-1/q,q}}^{q}\cdot
  \prod_{j=1}^{k}\norm{v_{j}'}_{W^{1-1/q,q}}.
 \end{equation}

%% file: experiments.tex

\section{Stability and isotopy preservation tests}
\label{sec:experiments}\label{sec:kapparho}

\newcommand{\pspace}{\spacefactor 3000{ }\ignorespaces}

\sbb{For the spacial discretization of our numerical scheme we follow~\cite{BRR} 
and approximate curves using piecewise cubic, continuously differentiable functions
on a given fixed partition of the parameter domain.
The time-stepping scheme~\eqref{eq:discgradflow} is different} 
from that introduced in~\cite{BRR} as we consider the $H^{2}$ gradient flow in the former
and the $L^{2}$ gradient flow in the latter. \SB{Moreover, the tangent-point functional
is treated fully explicitly which dramatically improves the numerical efficiency since the
assembly can be easily parallelized and fully populated matrices are avoided. 
The employed quadrature is the same as the one proposed and analyzed in~\cite{BRR}.
Every time step only requires the solution of a linear system of equations with sparse
system matrix that, due to the linearized arclength condition, has the structure of a 
saddle-point problem. \sbbb{We visualize discrete curves using an artificial small
thickness and a coloring encodes their curvature ranging from blue to 
yellow for small to larger curvature.} 

Our stability result guarantees an energy decay under a moderate condition on the
step size. This however does not imply that self-intersections are avoided. For this, 
the spacial discretization has to be sufficiently fine \sbb{relative to the initial
energy and the parameter $\vrho$}, so that the discretized 
tangent-point functional defines a discrete energy barrier that is larger than the
initial energy, cf.~\cite[Sect. 4.2]{BRR} for further details. 
Fortunately, these conditions to not conflict each other.


We test stability properties
of the flow depending on the parameters $\k$, $\vrho$,
the maximal spacial step size $h_{\max}$ (which is inversely proportional 
to the number of nodes), and the time step size $\tau$.
To this end, we use an initial curve
of length $\ell_{\mathrm{ini}}\approx39.9$ which
belongs to the $5_{2}$ knot class
and is given by
\newcommand{\vektor}[1]{\begin{pmatrix}#1\end{pmatrix}}
\begin{align*}
x\mapsto \tfrac1{100}
&\left[\vektor{-33\\-57\\34}\cos(2\pi x)
+\vektor{43\\99\\-21}\sin(2\pi x)
+\vektor{0\\-54\\-100}\cos(4\pi x) \right. \\
&\quad{}+\vektor{214\\-159\\-93}\sin(4\pi x)
+\vektor{101\\-117\\-27}\cos(6\pi x)
+\vektor{-47\\-5\\-16}\sin(6\pi x) \\
&\quad{}+\left.\vektor{0\\-31\\52}\cos(8\pi x)
+\vektor{11\\-45\\84}\sin(8\pi x) \right].
\end{align*}

\begin{table}[p]
%
\newcount\fileno
\newcommand{\exper}[3][0]{&
\global\advance\fileno1\relax\makebox[2ex][r]{}&
\ifnum#2=0no\else yes\ifnum#2>1\va\fi
\fi&\ifnum#3>0yes\else no\fi}
\newcommand{\va}{\makebox[0ex][l]}

\begin{tabular}{cccccccccc}
\multicolumn{10}{c}{(a)\pspace $\kappa=0$, $\rho = 1$}\\\hline
 \# nodes && \multicolumn{2}{c}{$\tau = \tfrac15h_{\max}^{1/2}$} && \multicolumn{2}{c}{$\tau = \tfrac15h_{\max}$} && \multicolumn{2}{c}{$\tau = \tfrac15h_{\max}^{2}$} \\
 && stab. & isot. && stab. & isot. && stab. & isot. \\\hline
 \textcolor{white}{0}%
  50 \exper00 \exper00 \exper00 \\
 100 \exper00 \exper00 \exper21 \\
 200 \exper00 \exper01 \exper11 \\
 400 \exper00 \exper11 \exper11 \\
 \hline
\end{tabular}
\bigskip

\begin{tabular}{cccccccccc}
\multicolumn{10}{c}{(b)\pspace $\kappa=\tfrac1{10}$, $\rho = 1$}\\\hline
 \# nodes && \multicolumn{2}{c}{$\tau = \tfrac15h_{\max}^{1/2}$} && \multicolumn{2}{c}{$\tau = \tfrac15h_{\max}$} && \multicolumn{2}{c}{$\tau = \tfrac15h_{\max}^{2}$} \\
 && stab. & isot. && stab. & isot. && stab. & isot. \\\hline
 \textcolor{white}{0}%
  50 \exper00 \exper00 \exper00 \\
 100 \exper00 \exper00 \exper21 \\
 200 \exper00 \exper21 \exper11 \\
 400 \exper00 \exper11 \exper11 \\
 \hline
\end{tabular}
\bigskip

\begin{tabular}{cccccccccc}
\multicolumn{10}{c}{(c)\pspace $\kappa=1$, $\rho = 1$}\\\hline
 \# nodes && \multicolumn{2}{c}{$\tau = \tfrac15h_{\max}^{1/2}$} && \multicolumn{2}{c}{$\tau = \tfrac15h_{\max}$} && \multicolumn{2}{c}{$\tau = \tfrac15h_{\max}^{2}$} \\
 && stab. & isot. && stab. & isot. && stab. & isot. \\\hline
 \textcolor{white}{0}%
  50 \exper00 \exper00 \exper00 \\
 100 \exper00 \exper00 \exper11 \\
 200 \exper00 \exper11 \exper11 \\
 400 \exper00 \exper11 \exper11 \\
 \hline
\end{tabular}
\bigskip

\begin{tabular}{cccccccccc}
\multicolumn{10}{c}{(d)\pspace $\kappa=1$, $\rho = \tfrac1{10}$}\\\hline
 \# nodes && \multicolumn{2}{c}{$\tau = \tfrac15h_{\max}^{1/2}$} && \multicolumn{2}{c}{$\tau = \tfrac15h_{\max}$} && \multicolumn{2}{c}{$\tau = \tfrac15h_{\max}^{2}$} \\
 && stab. & isot. && stab. & isot. && stab. & isot. \\\hline
 \textcolor{white}{0}%
  50 \exper11 \exper11 \exper11 \\
 100 \exper11 \exper11 \exper11 \\
 200 \exper11 \exper11 \exper11 \\
 400 \exper11 \exper11 \exper11 \\
 \hline
\end{tabular}
\bigskip

\begin{tabular}{cccccccccc}
\multicolumn{10}{c}{(e)\pspace $\kappa=1$, $\rho = \tfrac1{100}$}\\\hline
 \# nodes && \multicolumn{2}{c}{$\tau = \tfrac15h_{\max}^{1/2}$} && \multicolumn{2}{c}{$\tau = \tfrac15h_{\max}$} && \multicolumn{2}{c}{$\tau = \tfrac15h_{\max}^{2}$} \\
 && stab. & isot. && stab. & isot. && stab. & isot. \\\hline
 \textcolor{white}{0}%
  50 \exper10 \exper10 \exper10 \\
 100 \exper11 \exper11 \exper11 \\
 200 \exper11 \exper11 \exper11 \\
 400 \exper11 \exper11 \exper11 \\
 \hline
\end{tabular}
\bigskip


\caption{\label{tab:kapparho} 
Experimental results for the energy stability and isotopy preservation test in 
Section~\ref{sec:kapparho} for different ratios $\k/\vrho$ (cases (a)-(e)),
different spacial discretizations with $h_{\max} \simeq (\#\, \text{nodes})^{-1}$, 
and different choices of the step size $\tau$ relative to $h_{\max}$. Entries 
``yes/no'' indicate whether energy stability and isotopy preservation are observed.}
\end{table}

The parameters \sbbb{and results} of our experiments are listed in Table~\ref{tab:kapparho}.
The entries in the column ``stab.\@'' indicate whether the
discrete energy stability condition
\[  \frac{E(u^{k+1}) - E(u^{k})}\tau \le \tfrac32\tau^{1/2} \]
is satisfied. \sbbb{Allowing for a small tolerance on the right-hand side
accounts for discretization errors related to quadrature.}
The ``isot.\@'' column reports on whether the isotopy type of
the initial curve is maintained during the evolution.
In each case we observed the evolution for about fifty to one-hundred time steps.
As $h_{\max}$ is close to one for the case of $50$ nodes,
the entries in the corresponding row coincide.

Note that we have to choose $q<4$
in order to absorb the error of cutting out an $\eps$-neighborhood
of the diagonal of $\R/\Z\times\R/\Z$ in our discretization of the $\TP$ term, cf.~\cite{BRR}.
Throughout this experiment we use $q=3$ and $\veps = 2 h_{\max}$.

In general, both stability and isotopy maintenance are improved
as the time-step size or the spacial discretization is refined. The 
latter dependence is attributed to large values of the initial energy $E(u_h^0)$
for coarse spacial resolutions and the dependence of the step-size condition
identified in Proposition~\ref{prop:stab} as well as errors related to 
quadrature for $\TP$ in particular when $\vrho \gg \k$. 

If the ratio of $\k$ over $\vrho$ is small,
i.e., the self-avoiding part plays a prominent r\^ole,
stability seems to imply that the knot class is preserved, see Table~\ref{tab:kapparho}~(a)--(c).
In particular, this applies to $\k=0$ in~(a) which corresponds to a fully explicit discretization of the $\TP$-flow.
A typical instability (along with isotopy and significant length preservation
violation) is shown in Figure~\ref{fig:kapparho(b)}.

In the case of larger ratios of $\k$ over $\vrho$,
i.e., the bending energy dominates,
the scheme tends to be more stable, even for
relatively coarse spacial and temporal discretizations,
see Table~\ref{tab:kapparho}~(d)--(e).
However, stability does not guarantee preservation of the isotopy type.
An isotopy violation due to a too coarse spacial discretization
is depicted in Figure~\ref{fig:kapparho(e)}.


\newcommand{\bild}[2]{\fbox{\includegraphics[scale=.1,trim=220 30 170 70,clip]{knot_-3_200/#1.jpg}}\makebox[0ex][r]{\makebox[0ex][r]{\tiny#1}\hspace{15ex}\makebox[0ex][r]{\tiny length #2} }\,\ignorespaces}

\begin{figure}

%
%
%

\renewcommand{\bild}[3]{\fbox{\includegraphics[scale=.1,trim=220 30 170 70,clip]{#1}}\makebox[0ex][r]{\makebox[0ex][r]{\tiny#2}\hspace{15ex}\makebox[0ex][r]{\tiny length #3} }\,\ignorespaces}


\bild{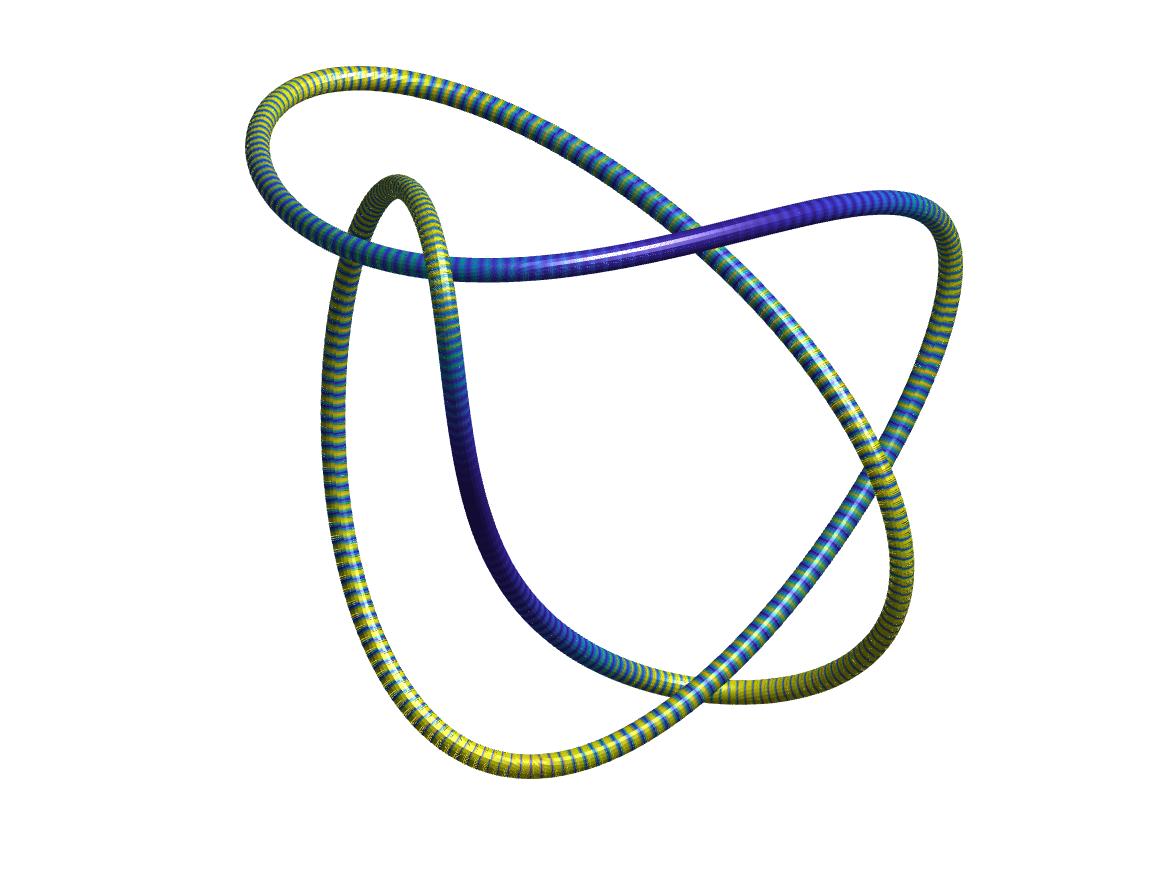}{1}{39.919092}
\bild{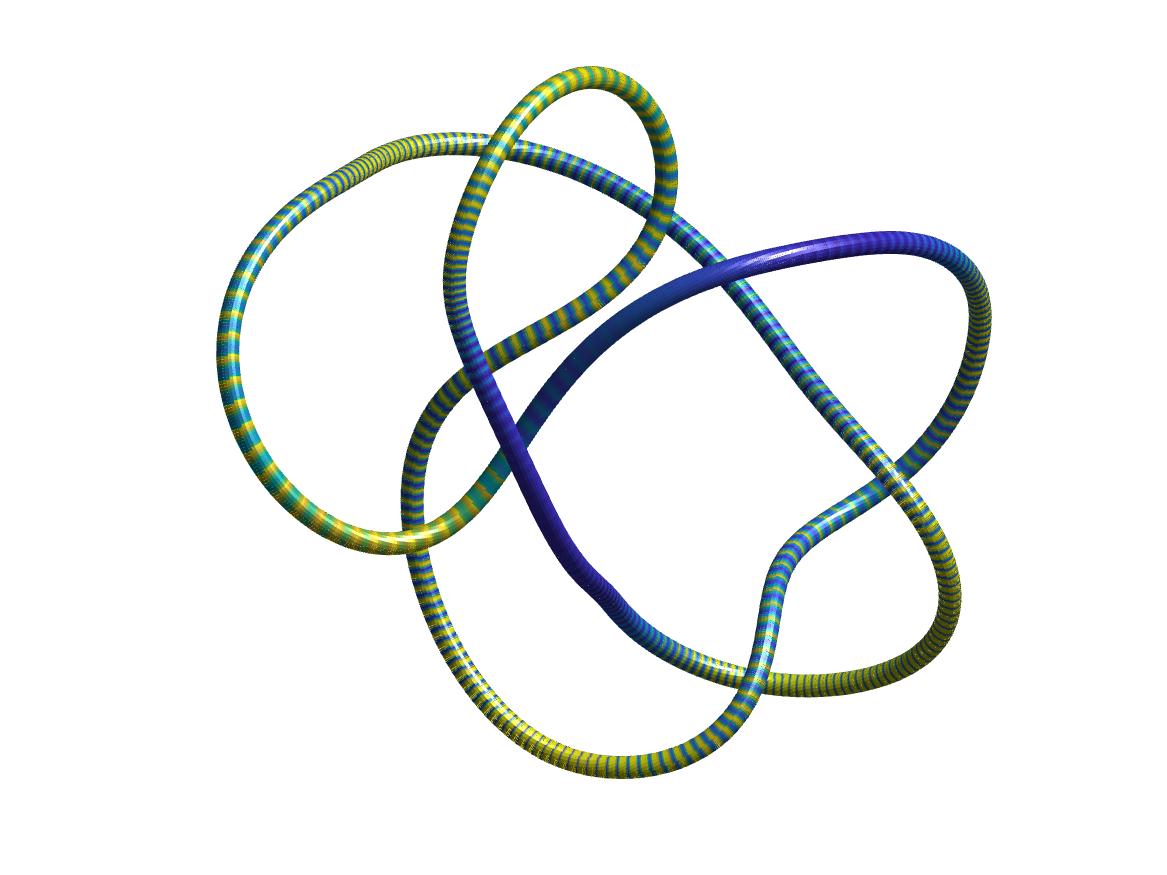}{2}{47.526655}
\bild{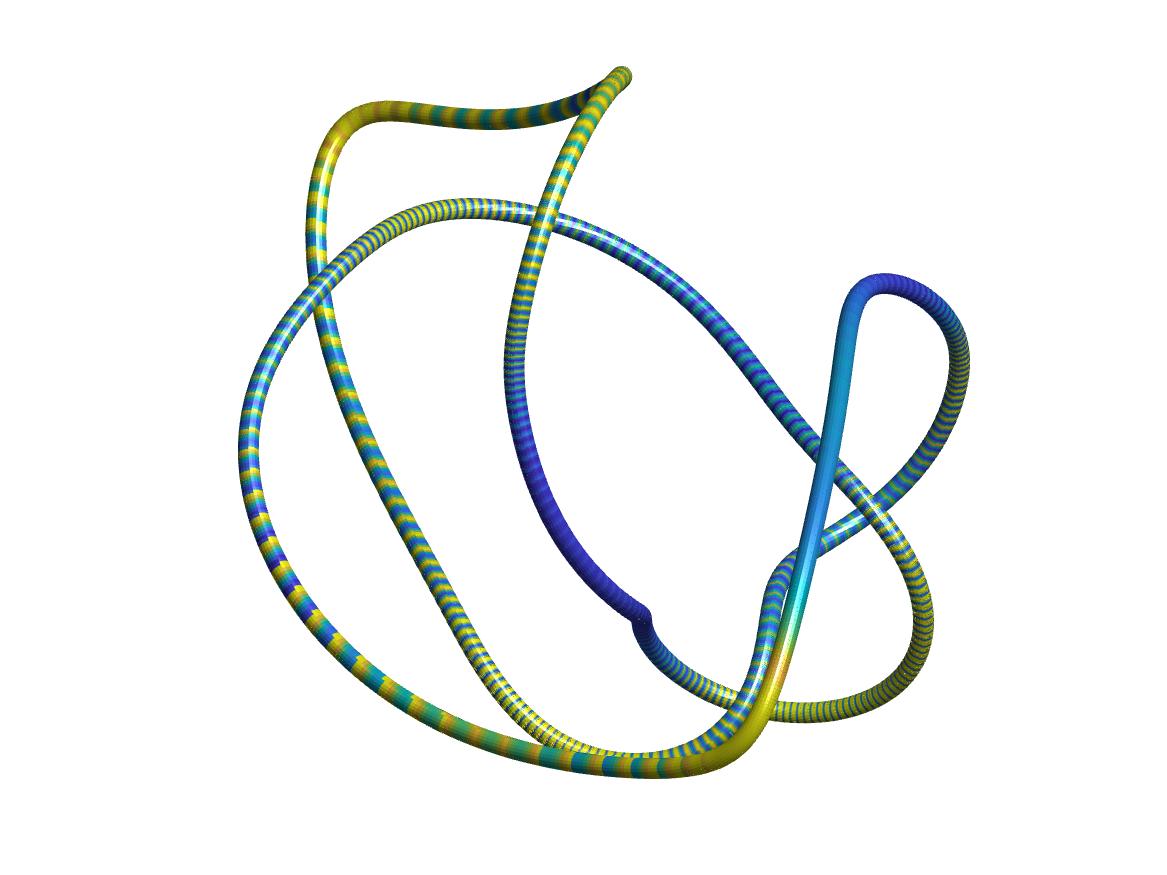}{3}{59.041003}
\bild{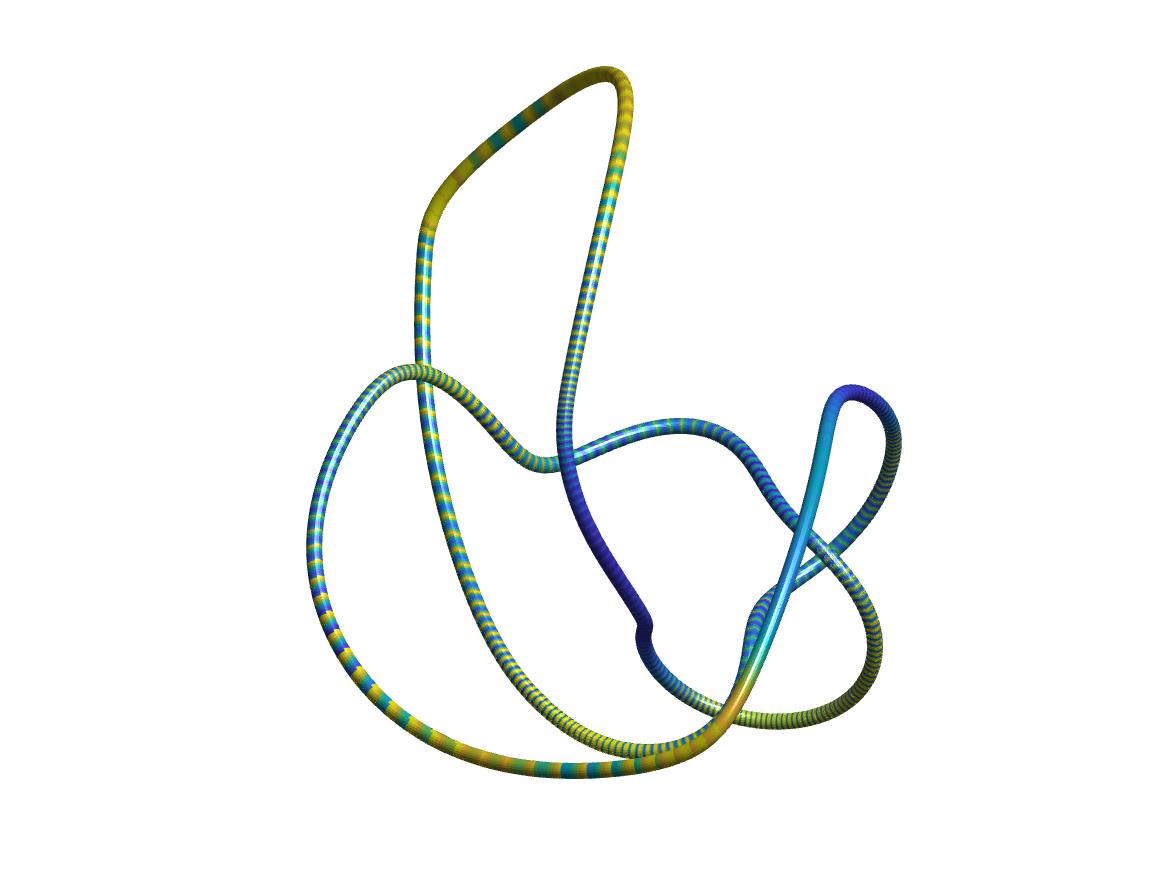}{4}{61.777820}

\bild{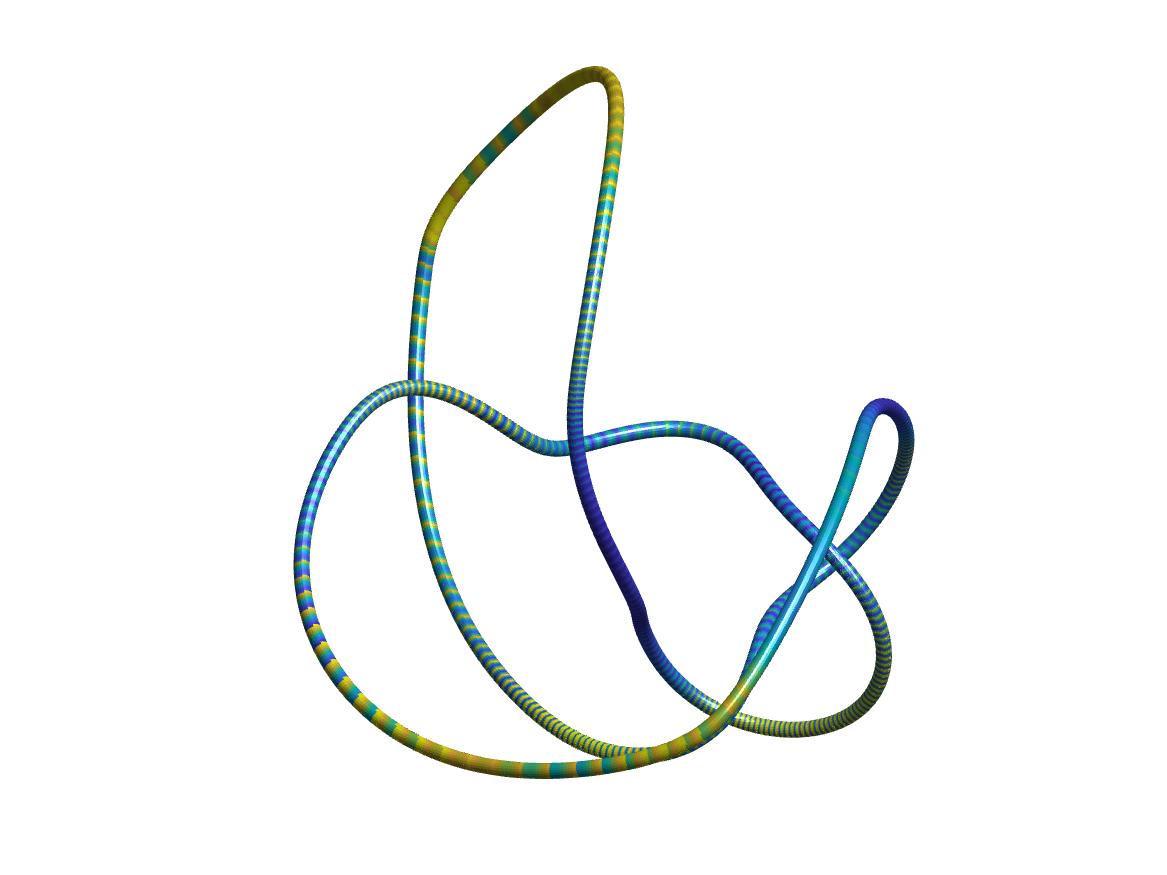}{6}{62.042947}
\bild{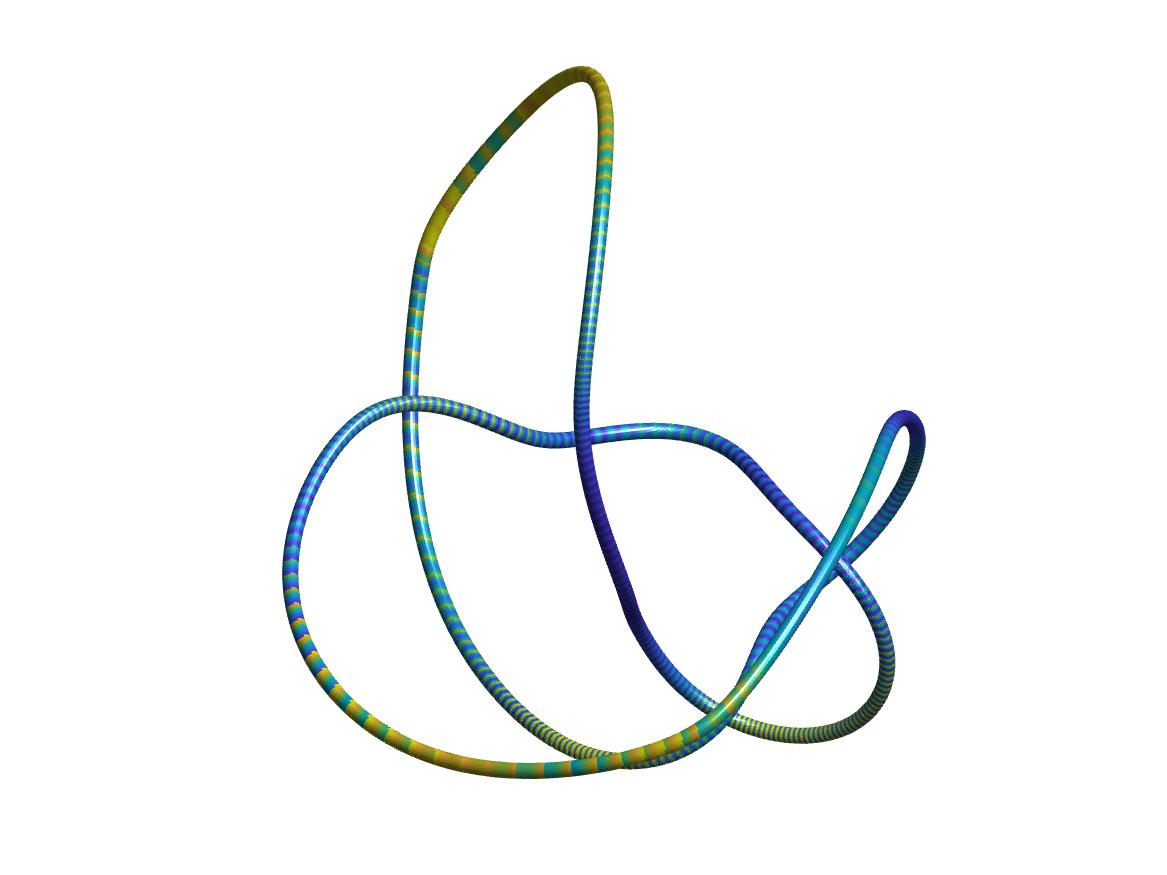}{8}{62.151359}
%
\bild{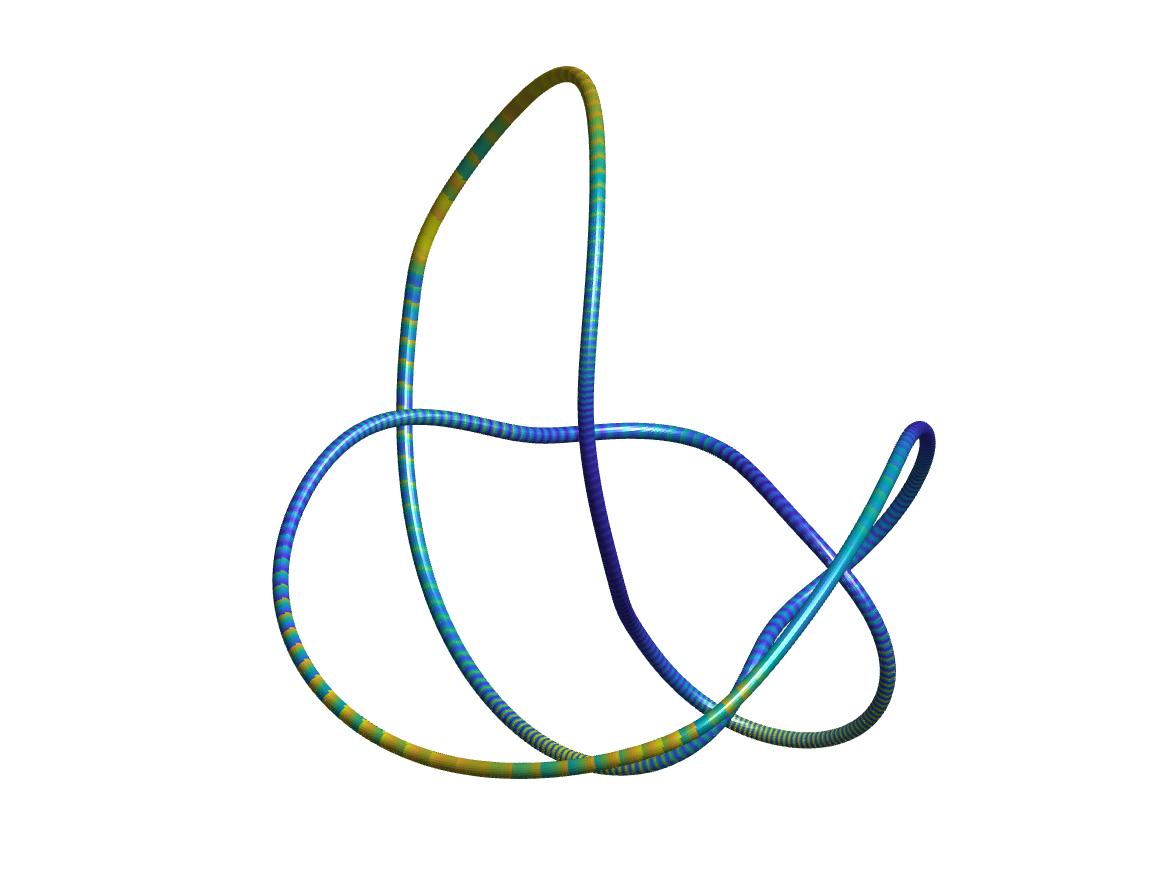}{10}{62.216288}
%
\bild{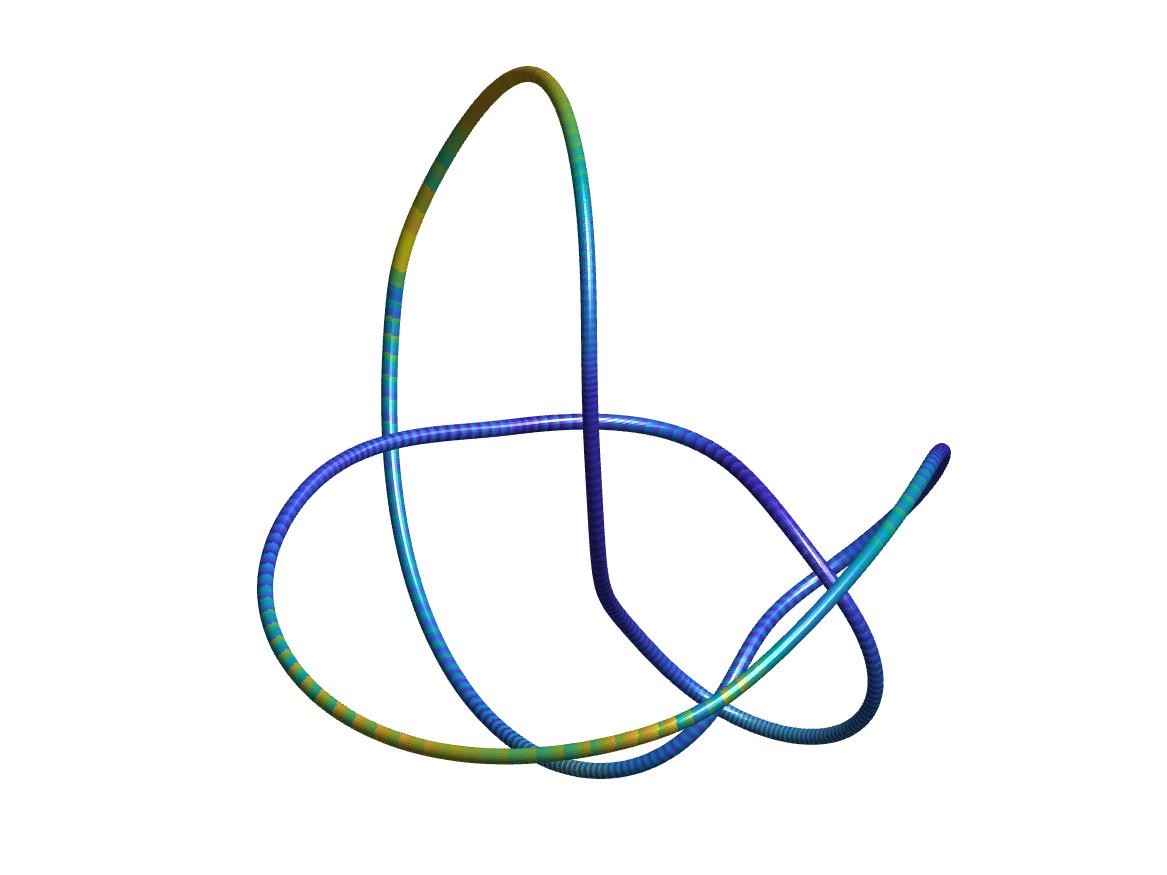}{16}{62.312458}

\vspace*{4mm}

\includegraphics[scale=.475]{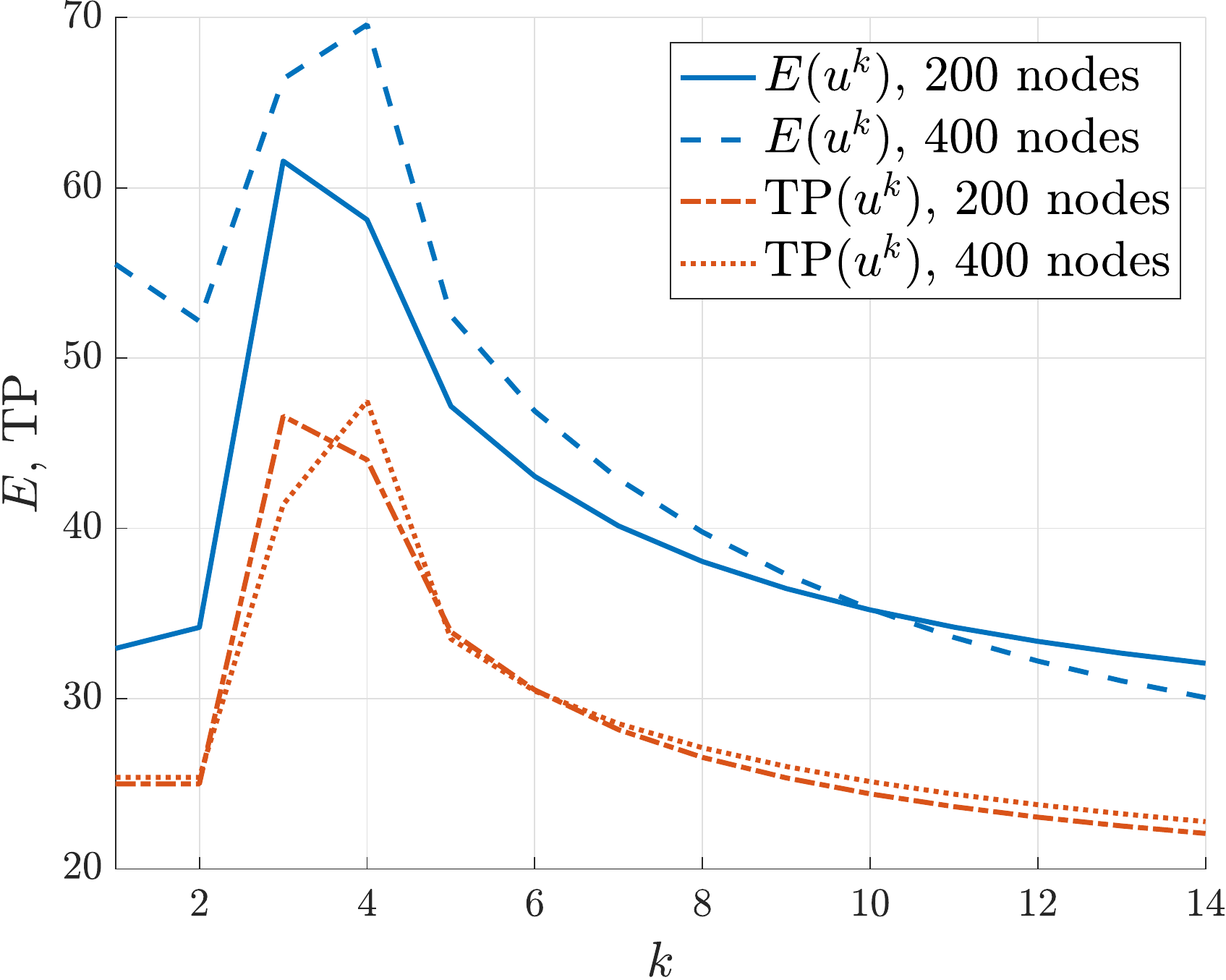}


\caption{\label{fig:kapparho(b)}\label{fig:kapparho(b)plot}
Instability accompanied by isotopy change for the inital curve defined in Section~\ref{sec:kapparho} 
with $\k=1/10$ and $\vrho = 1$ (case~(b) in Table~\ref{tab:kapparho}) and 
$\tau = \tfrac15h_{\max}^{1/2}$. 
Top: snapshots and curve lengths for $h_{\max} \simeq 1/400$
after $k=1,2,3,4,6,8,10,16$ time steps (note different magnification); 
bottom: energy and self-avoidance potentials for $k=1,2,\dots,14$ with
$h_{\max} \simeq 1/200$ and $h_{\max} \simeq 1/400$.} 
\end{figure}

%

\renewcommand{\bild}[2]{\fbox{\includegraphics[scale=.115,trim=280 100 250 270,clip]{#1}}\makebox[0ex][r]{{\tiny#2} }\,\ignorespaces}

\begin{figure}

\bild{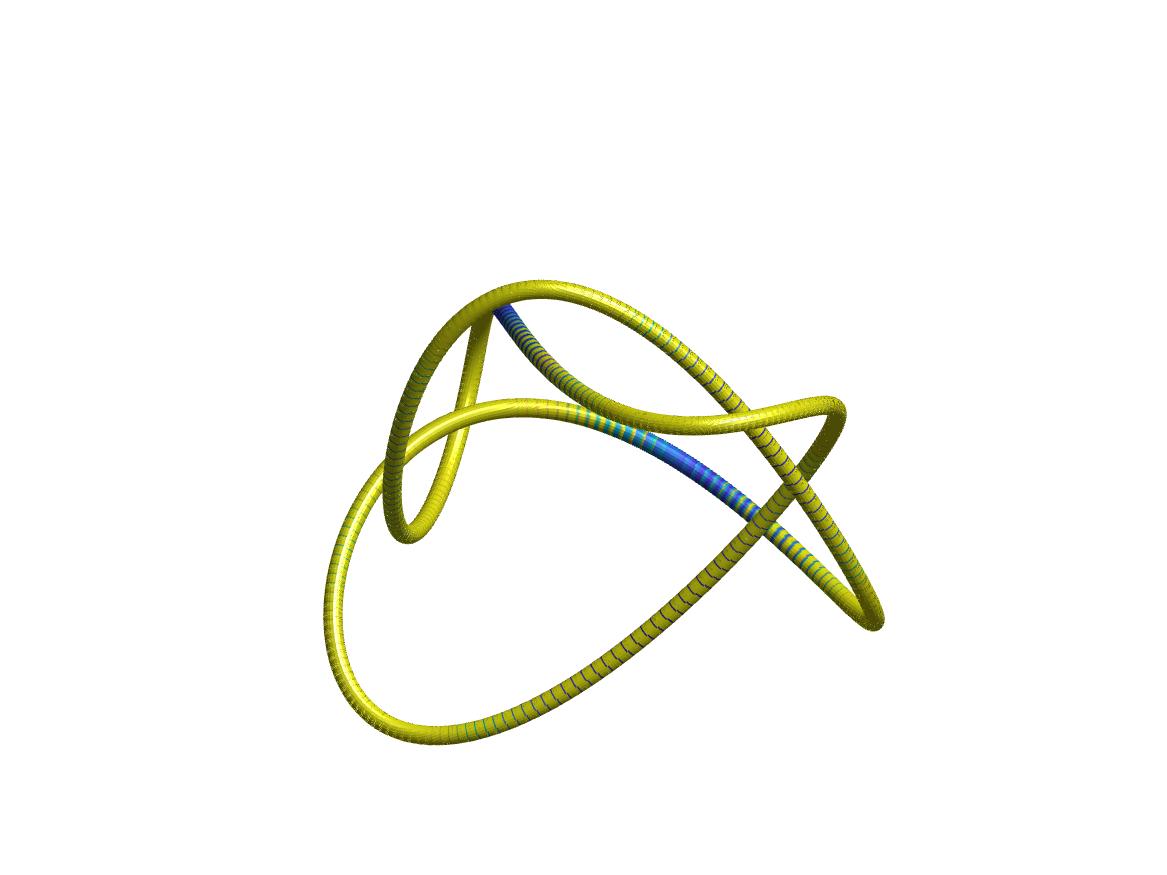}{1}
\bild{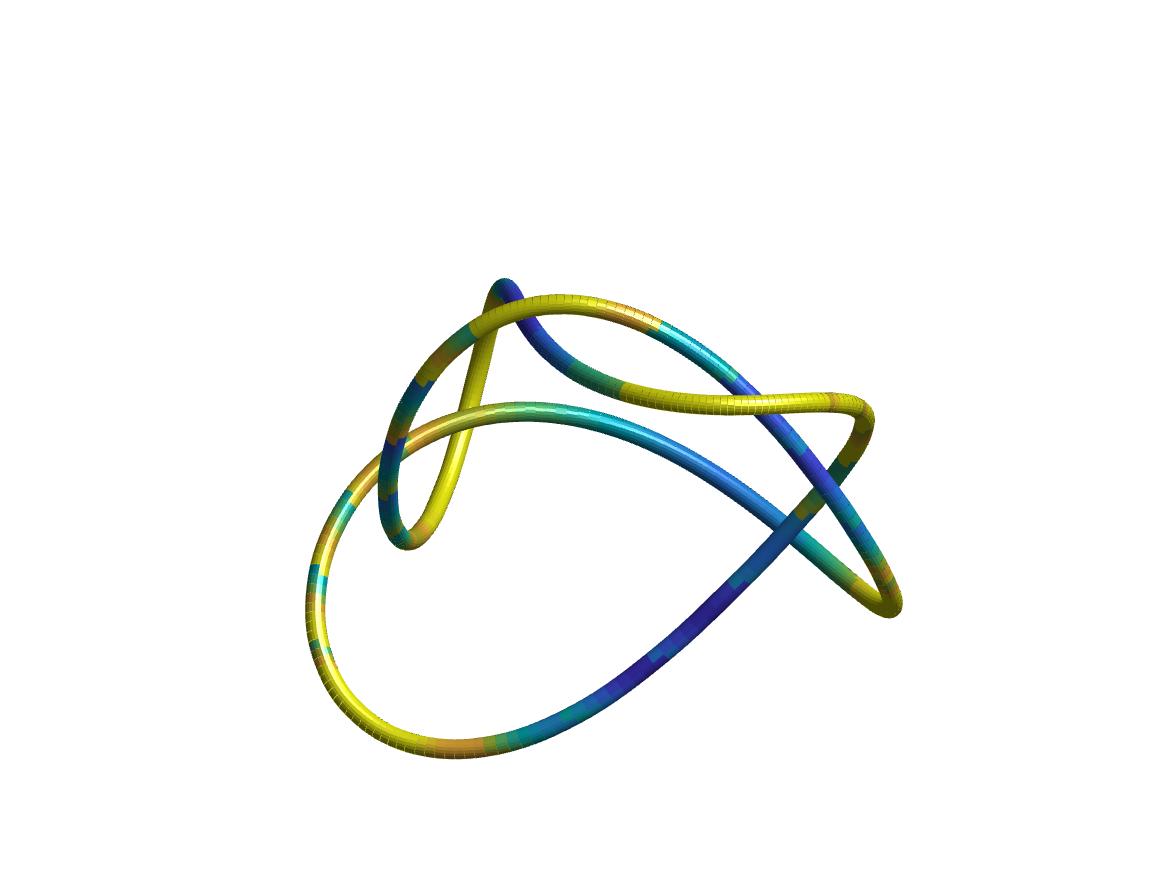}{4}
\bild{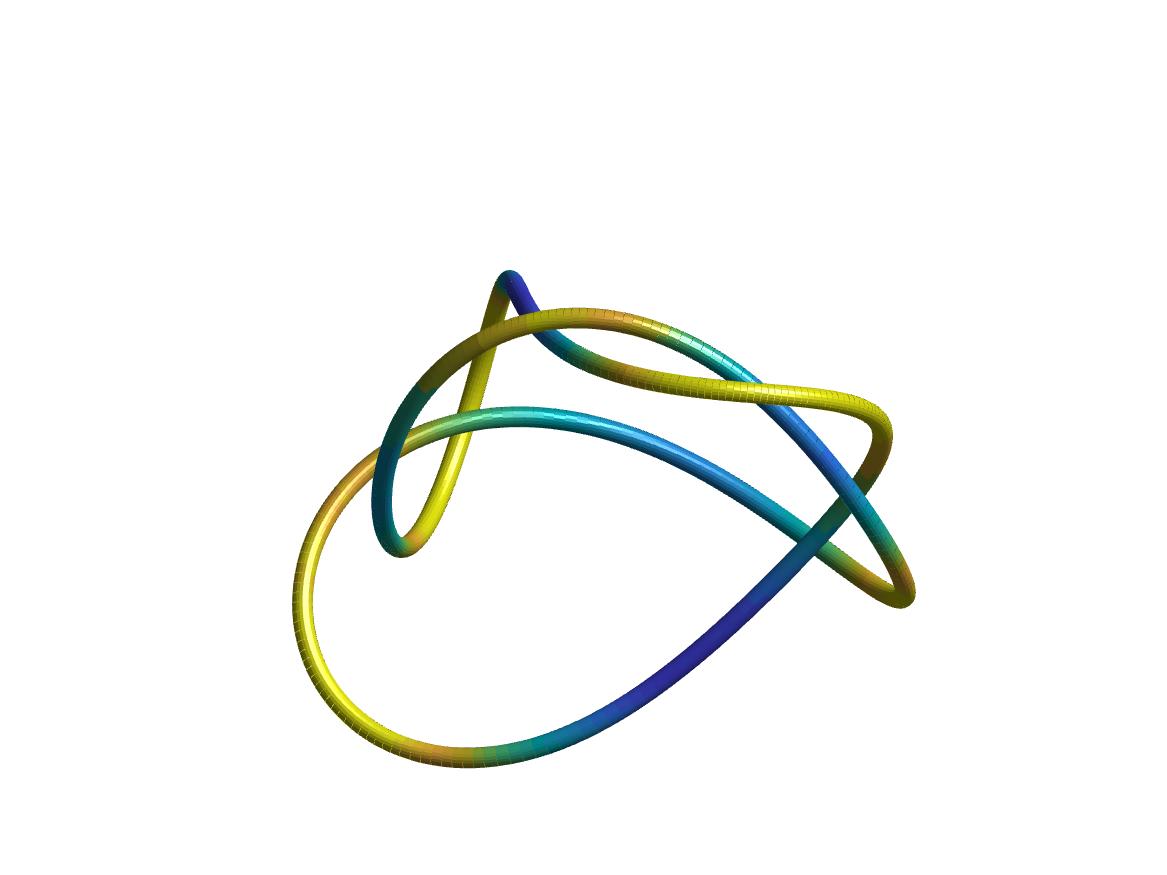}{8}
\bild{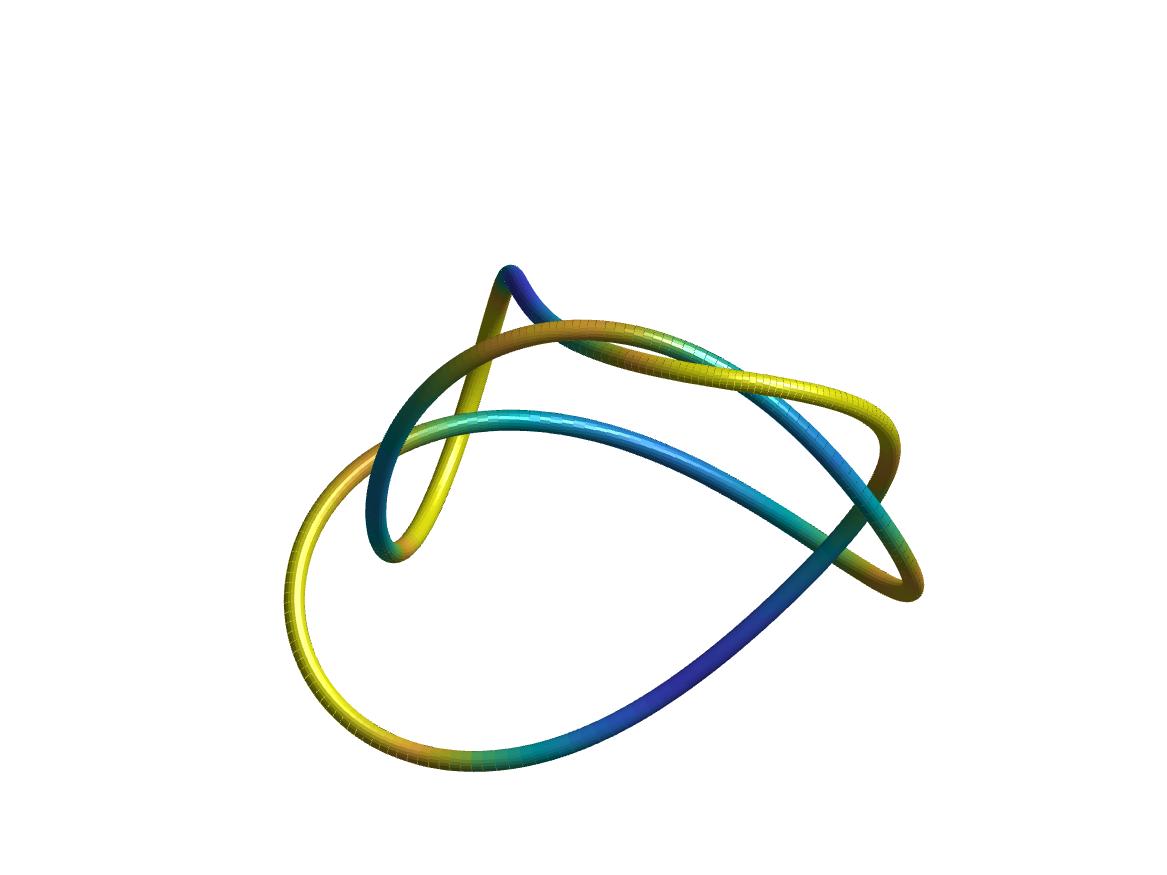}{12}
\makebox[0ex][l]{\hspace{-10.5ex}\raisebox{10ex}{\framebox[6ex]{\rule{0pt}{1.5ex}}}}\ignorespaces

\renewcommand{\bild}[2]{\fbox{\includegraphics[scale=.22,trim=530 470 378 270,clip]{#1}}\makebox[0ex][r]{\makebox[0ex][r]{\tiny#2}\hspace{9.5ex} }\,\ignorespaces}

\bild{knot_-3_50_12.jpg}{12}
\bild{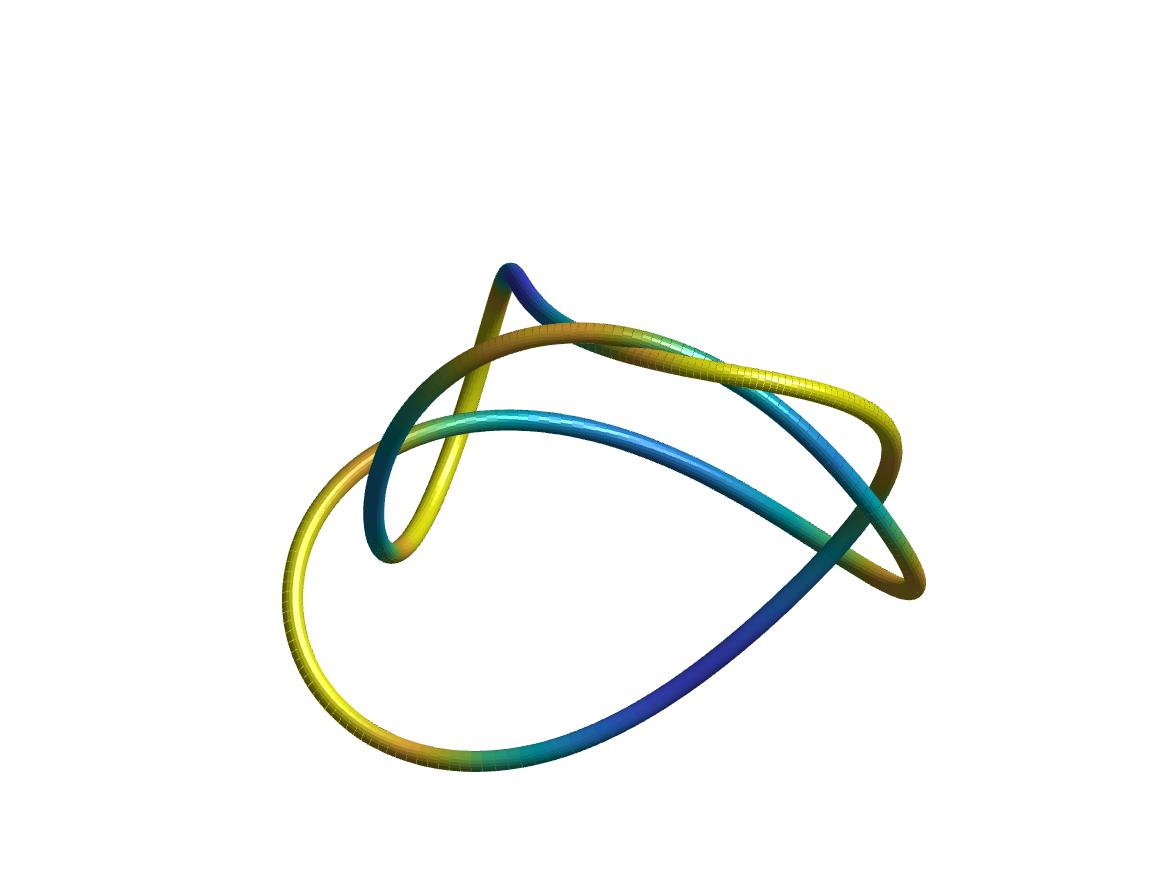}{13}
\bild{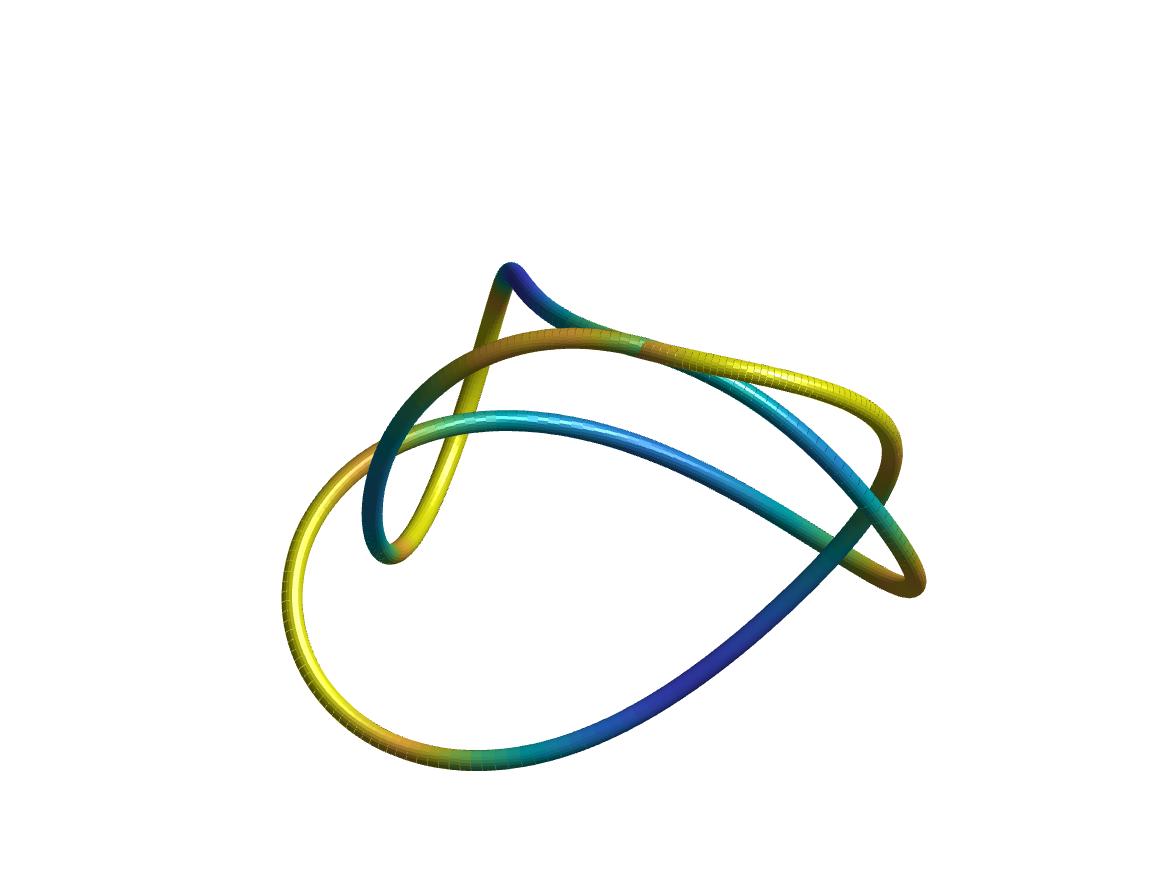}{14}
\bild{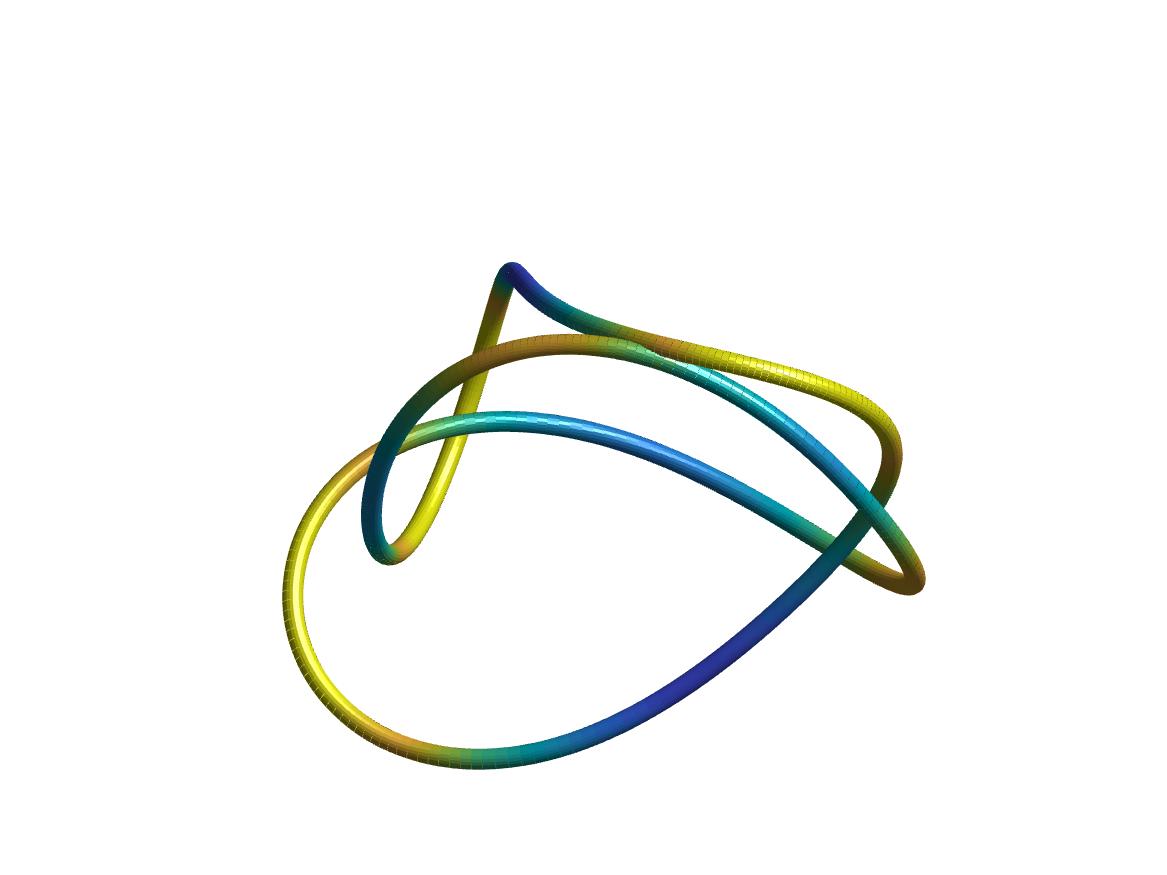}{15}
\bild{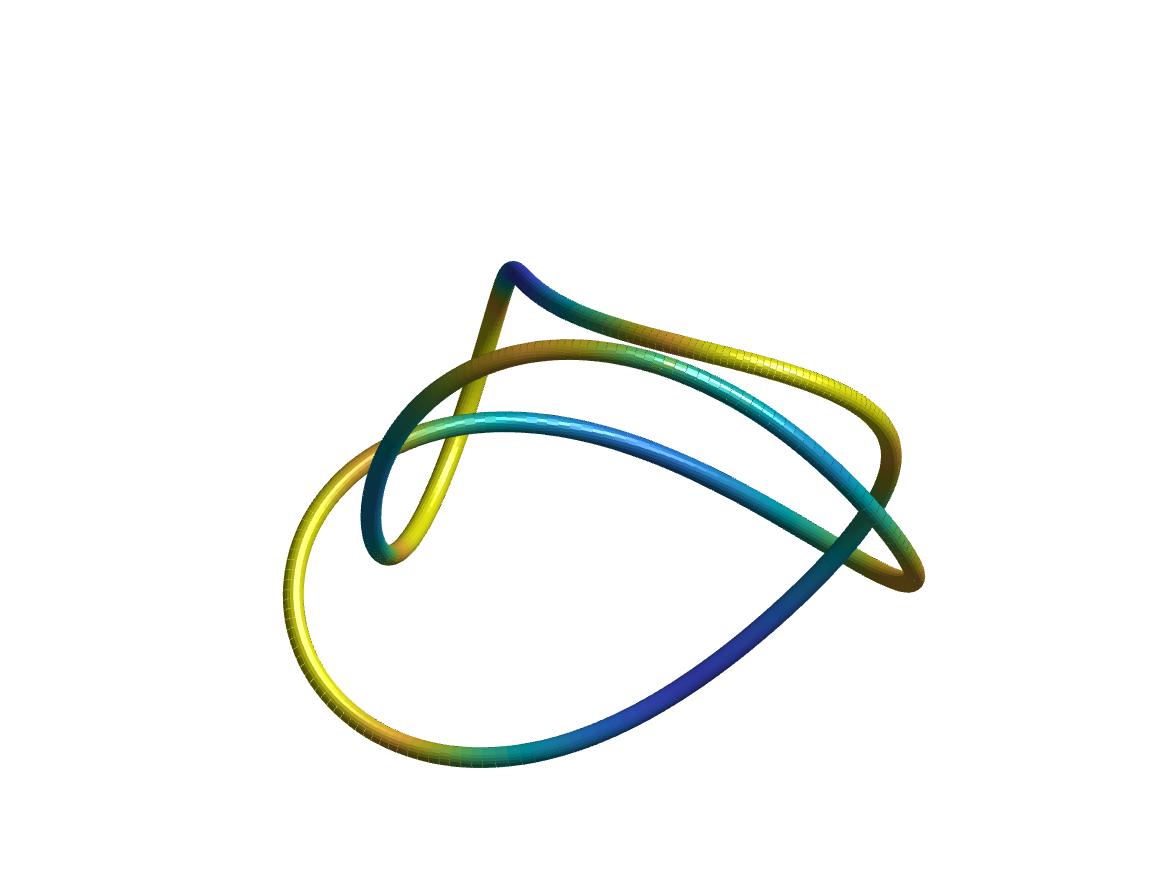}{16}

\renewcommand{\bild}[2]{\fbox{\includegraphics[scale=.115,trim=280 110 250 210,clip]{#1}}\makebox[0ex][r]{\makebox[0ex][r]{\tiny#2}\hspace{0ex} }\,\ignorespaces}

\bild{knot_-3_50_16.jpg}{16}\makebox[0ex][l]{\hspace{-10.5ex}\raisebox{10ex}{\framebox[6ex]{\rule{0pt}{1.5ex}}}}\ignorespaces
\bild{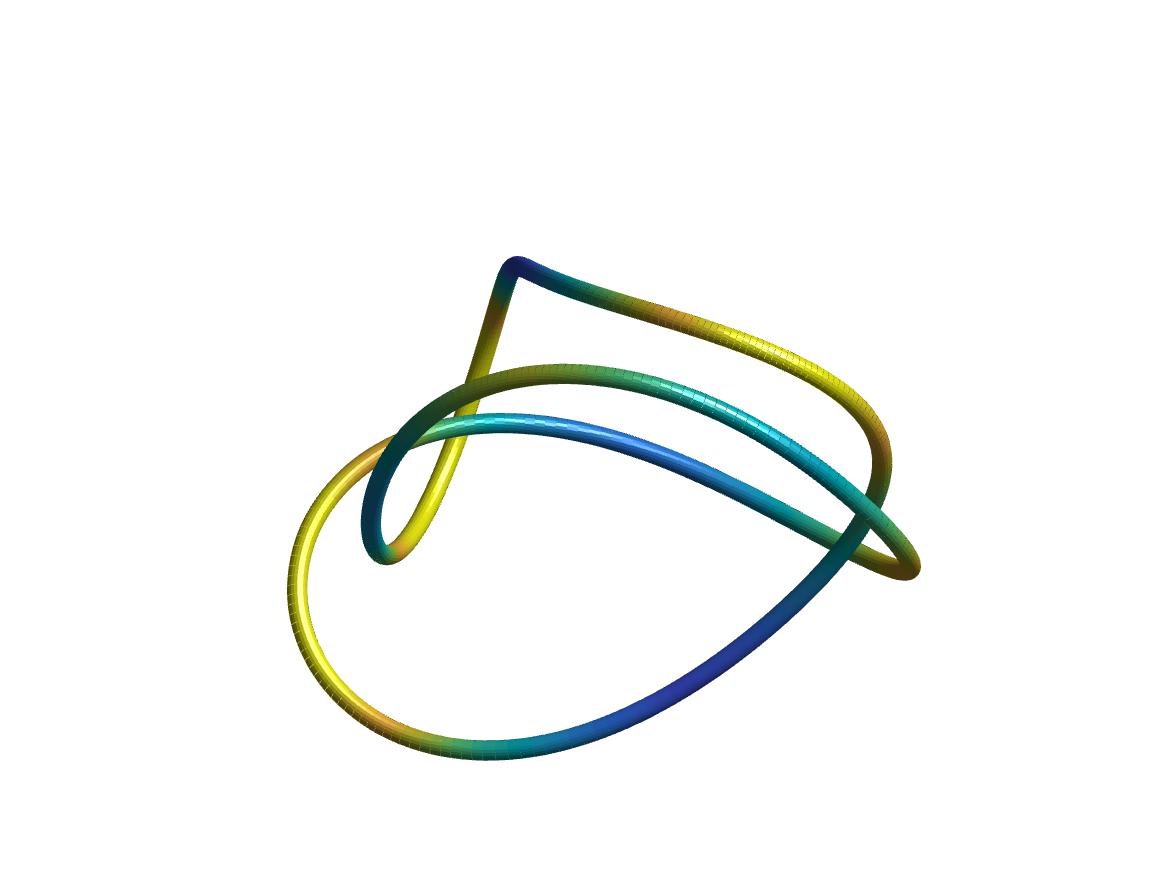}{20}
\bild{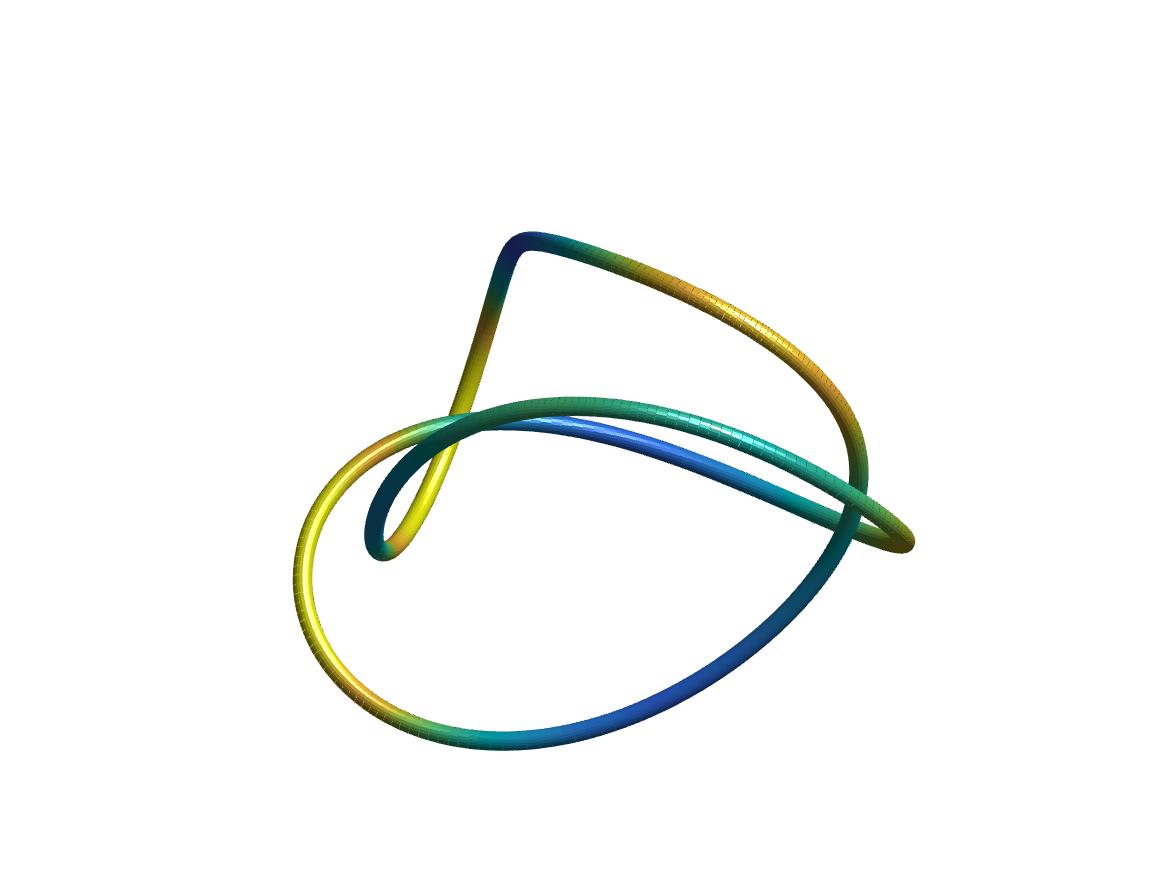}{30}
\bild{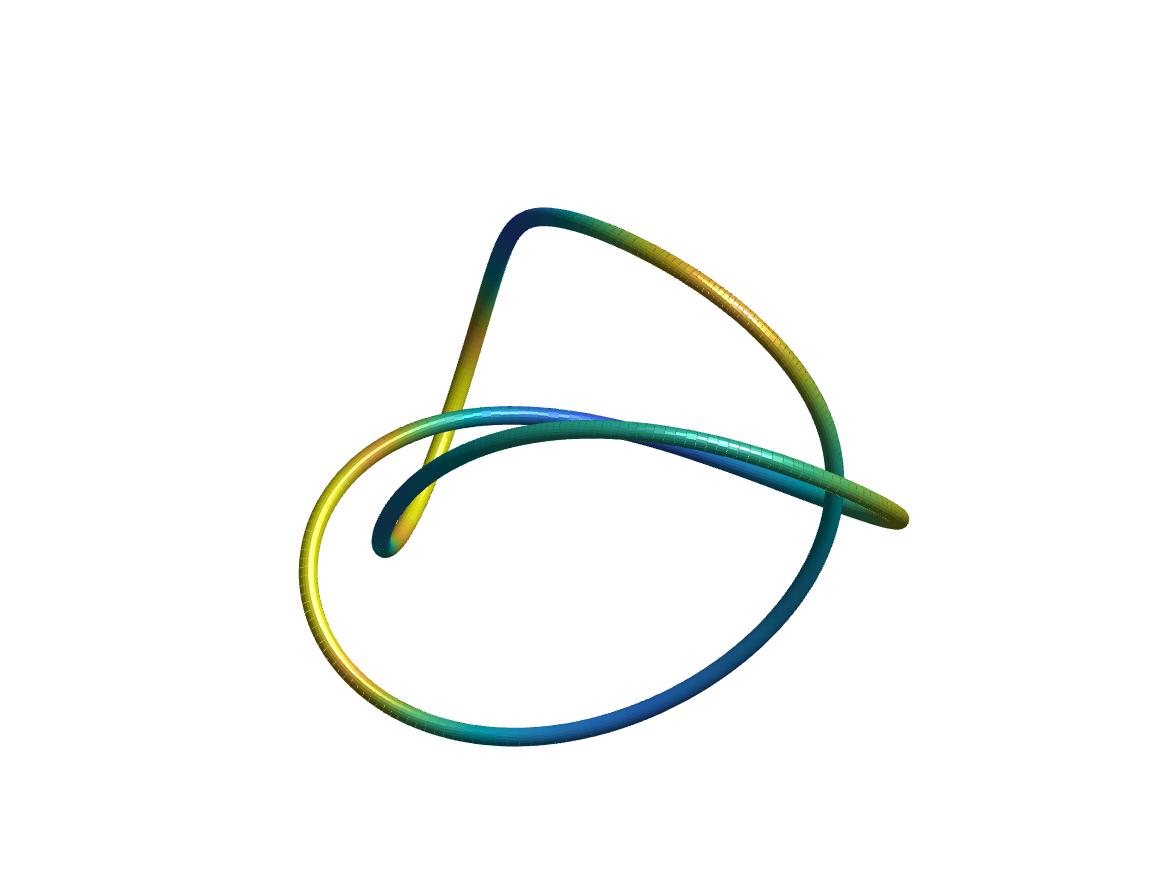}{40}

%

\vspace*{3mm} 

\includegraphics[scale=.475]{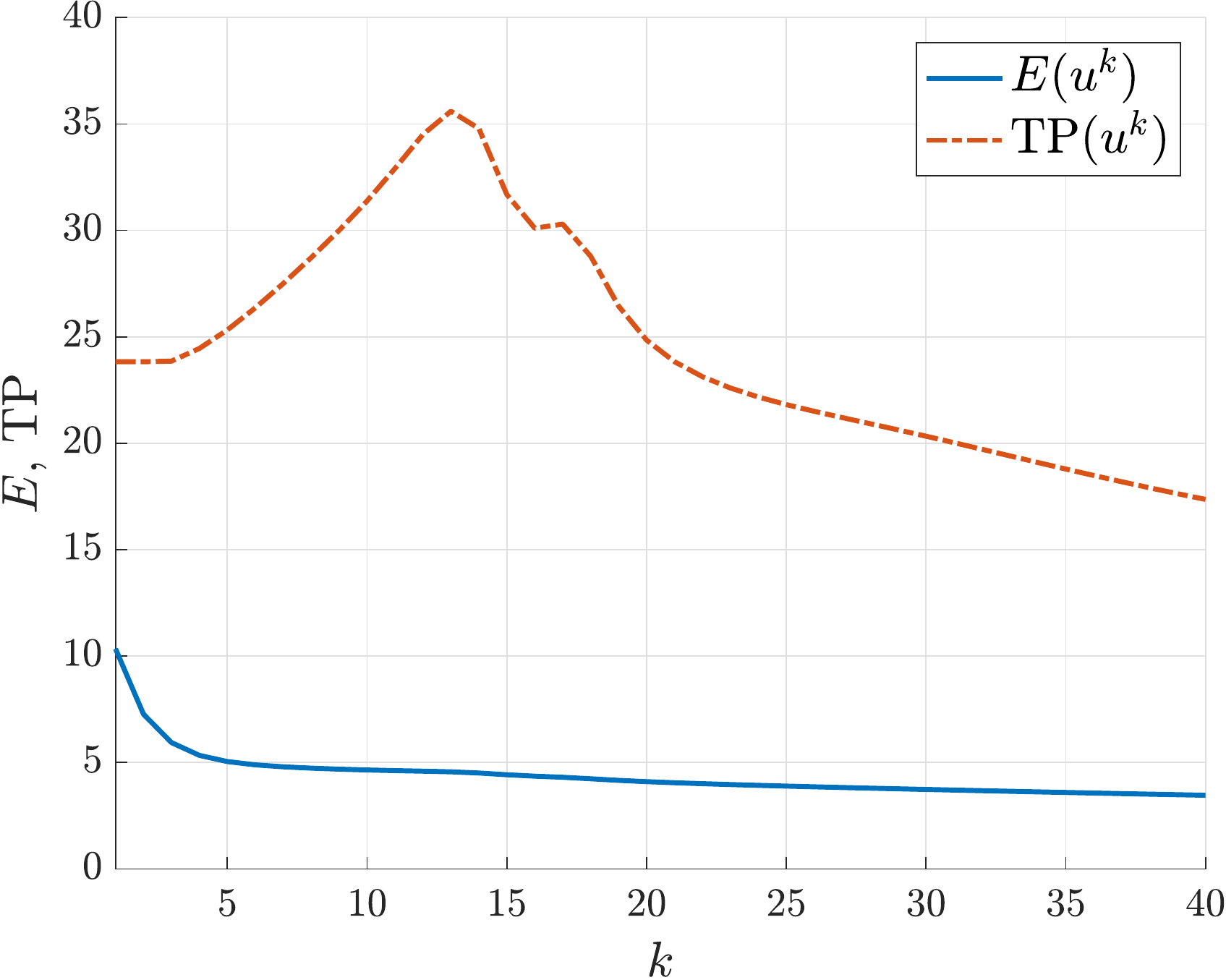}

\caption{\label{fig:kapparho(e)plot}\label{fig:kapparho(e)}
Energy stable evolution accompanied by an isotopy change for the 
inital curve defined in Section~\ref{sec:kapparho} 
with $\k=1$ and $\vrho = 1/100$ (case~(e) in Table~\ref{tab:kapparho}). 
Top: a self-penetration of the curve occurs after~$14$ iterations;
bottom: while the total energy decreases the tangent-point functional has
a peak when the isotopy change takes place.} 
\end{figure}

\section{Simulating elastic knots}\label{sec:elastic}

\SB{A primary goal in the experimental study of elastic knots is the determination 
of global minimizers of the energy functional $E$. As is typical
for gradient flows and in particular for those related to singularly perturbed
functionals, the evolution may become stationary at local minimizers or nearly stationary
at so-called metastable states, see Carr and Pego~\cite{CarPeg89} and 
Otto and Reznikoff~\cite{OtRe07}. In some cases small perturbations of the iterates
can avoid these phenomena.}

In this section we report on experiments targeted at the approximation of
elastic knots~\cite{GRvdM} as global minimizers of $E$ when $\kappa\gg\vrho$.
In all cases, the evolution reaches some ``stationary state'' after finite time
which seems to be a stable configuration.
In general, it is difficult to decide whether
it is in fact
a local minimum, without even being a \emph{global} minimum. 
Experiments with physical wires suggest the existence of
several non-global local minimizers.

\begin{table}
%
\begin{tabular}{lccccccccc}\hline
 & \# nodes & $h_{\max}$ & pert. & $\ell_{\mathrm{ini}}$ & $\ell_{\mathrm{end}}$ \\\hline
 \ref{sec:unkn}
 & $376$ & $0.1255$ & yes & 46.863580 & 46.855587 \\
 \ref{sec:tref} (a) 
 & $401$ & $0.1312$ & no & 49.996110 & 49.999110 \\
 \ref{sec:tref} (b) 
 & $401$ & $0.1312$ & yes & 49.996110 & 49.997392 \\
 \ref{sec:f8} (a) 
 & $400$ & $0.1370$ & yes & 49.871712 & 49.864798 \\
 \ref{sec:f8} (b) 
 & $415$ & $0.1306$ & yes & 49.884779 & 49.878556 \\\hline 
\end{tabular}
\bigskip

\caption{Discretization parameters for experiments defined in Subsections~\ref{sec:unkn}--\ref{sec:f8},
we always use $\tau = \tfrac1{30}{h_{\max}}$ and $\veps = 2 h_{\max}$. }\label{tab:param}
\end{table}

Throughout this section we use $\k= 1$, $\vrho = 10^{-3}$, and $q=3.9$ to define $E$. 
The discretization parameters used for the experiments defined below
are listed in Table~\ref{tab:param}
where we use the same notation as in Section~\ref{sec:kapparho}.
The entries in the ``pert.\@'' column indicate whether a slight randomized perturbation
was performed every hundred steps.
The parameters $\ell_{\mathrm{ini}}$ and $\ell_{\mathrm{end}}$
denote the length of the curve at the initial step and the last
step of the evolution respectively, more precisely, the length
of the polygonal curve defined by the vertices.

Our general observation is that owing to the extreme ratio $\k/\vrho = 10^4$
one has to \sbbb{suitably} choose
the discretization parameters in order to prevent self-intersections
during the evolution. For all experiments reported below the moderate 
number of approximately $400$ nodes and the relation $\tau = \frac{1}{30} h_{\max}$
were found to be sufficient. 

\subsection{Unknot}\label{sec:unkn}

We experiment with an initial configuration proposed by Avvakumov and Sossinsky~\cite{sossinsky},
consisting of a polygon describing a unilateral triangle
with ``twisted vertices''. The discrete evolution is depicted along with an energy 
plot in Figure~\ref{fig:unkn(a)} and we observe that the algorithm gets stuck in a configuration
different from the global minimizer which is the round circle, cf.~\cite{GRvdM}.

It is likely that this is an analytical feature of
the gradient flow~\eqref{eq:gradflow} and not an artifact of
the numerical scheme.
Therefore, at least for small values of $\vrho>0$, the gradient flow~\eqref{eq:gradflow} does not seem} to be a candidate for
a retract of the unknots to the round circle (which exists
due to the Smale conjecture, see Hatcher~\cite{hatcher}).


\renewcommand{\bild}[2]{\fbox{\includegraphics[scale=.09,trim=180 100 160 70,clip]{#1}}\makebox[0ex][r]{\makebox[0ex][c]{\tiny#2}\hspace{7.67ex} }\,\ignorespaces}

\begin{figure}

\bild{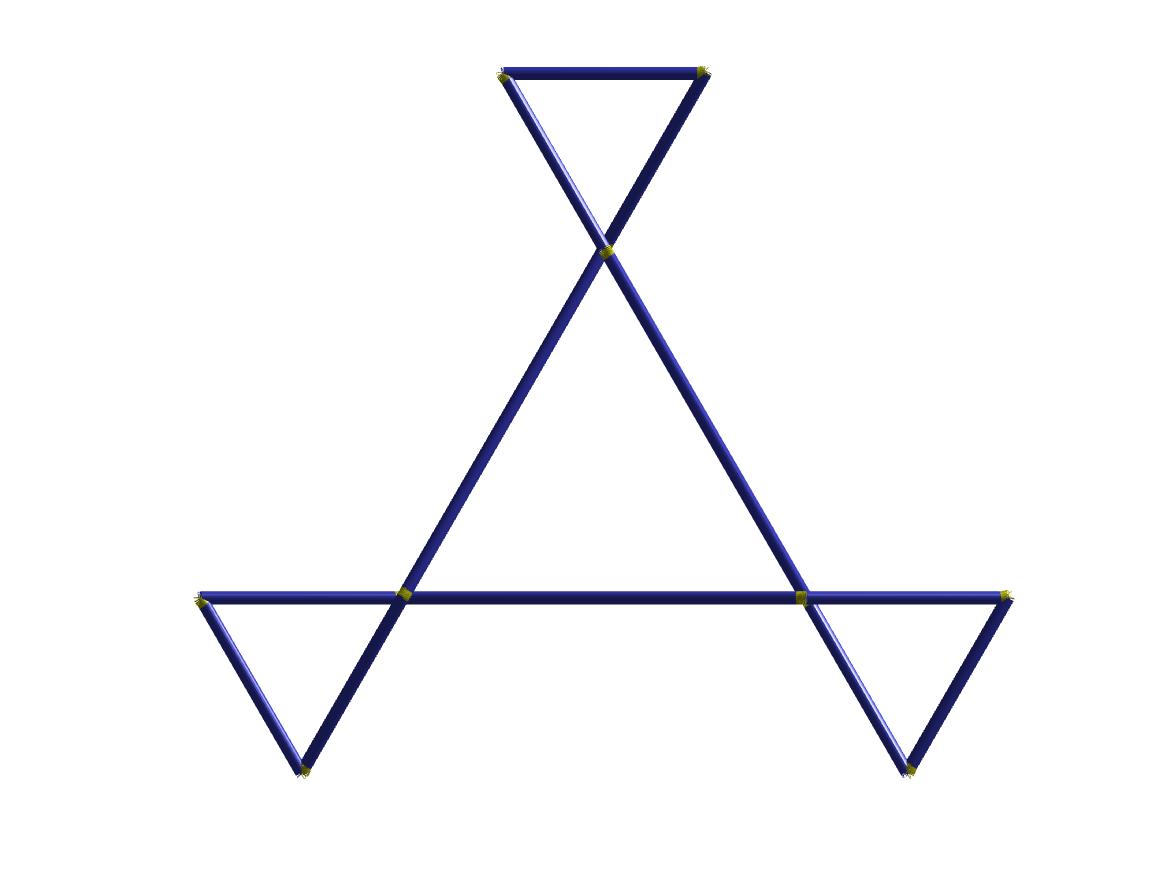}{1}
\bild{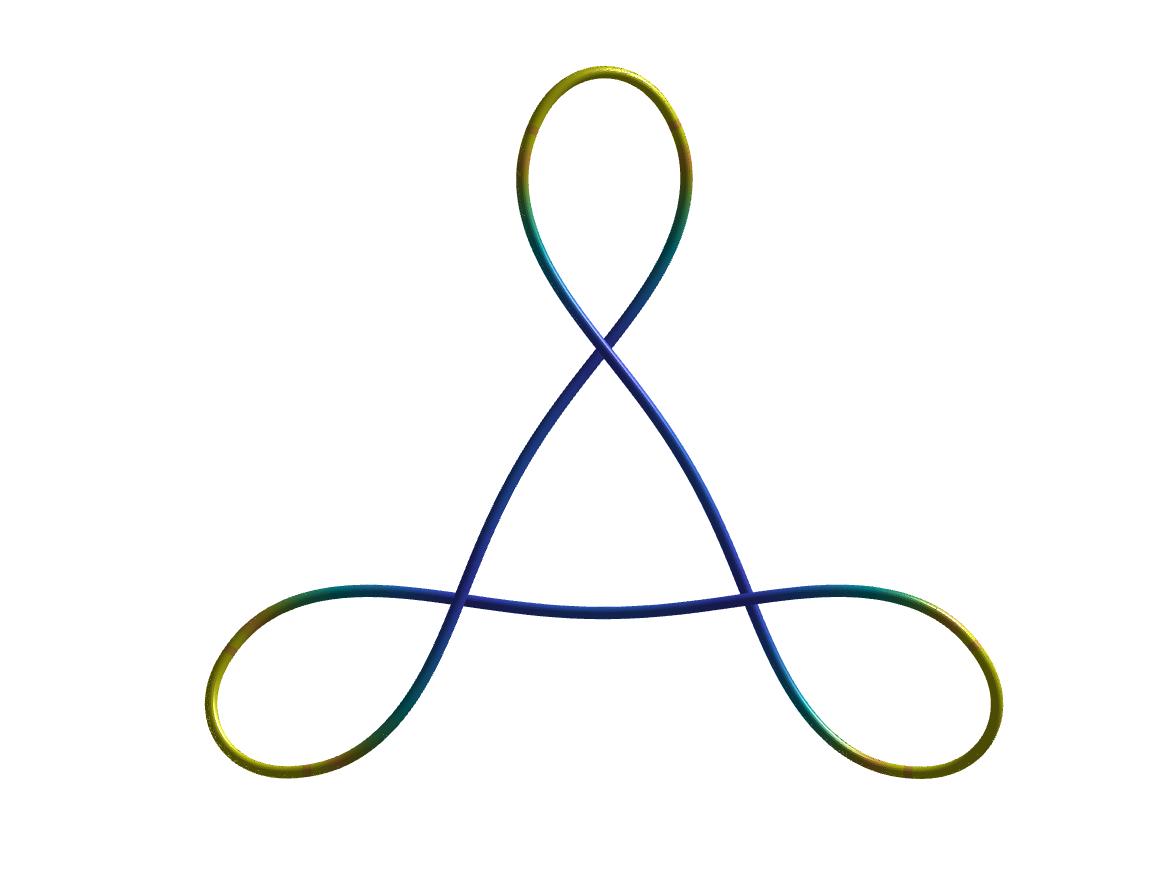}{100}
\bild{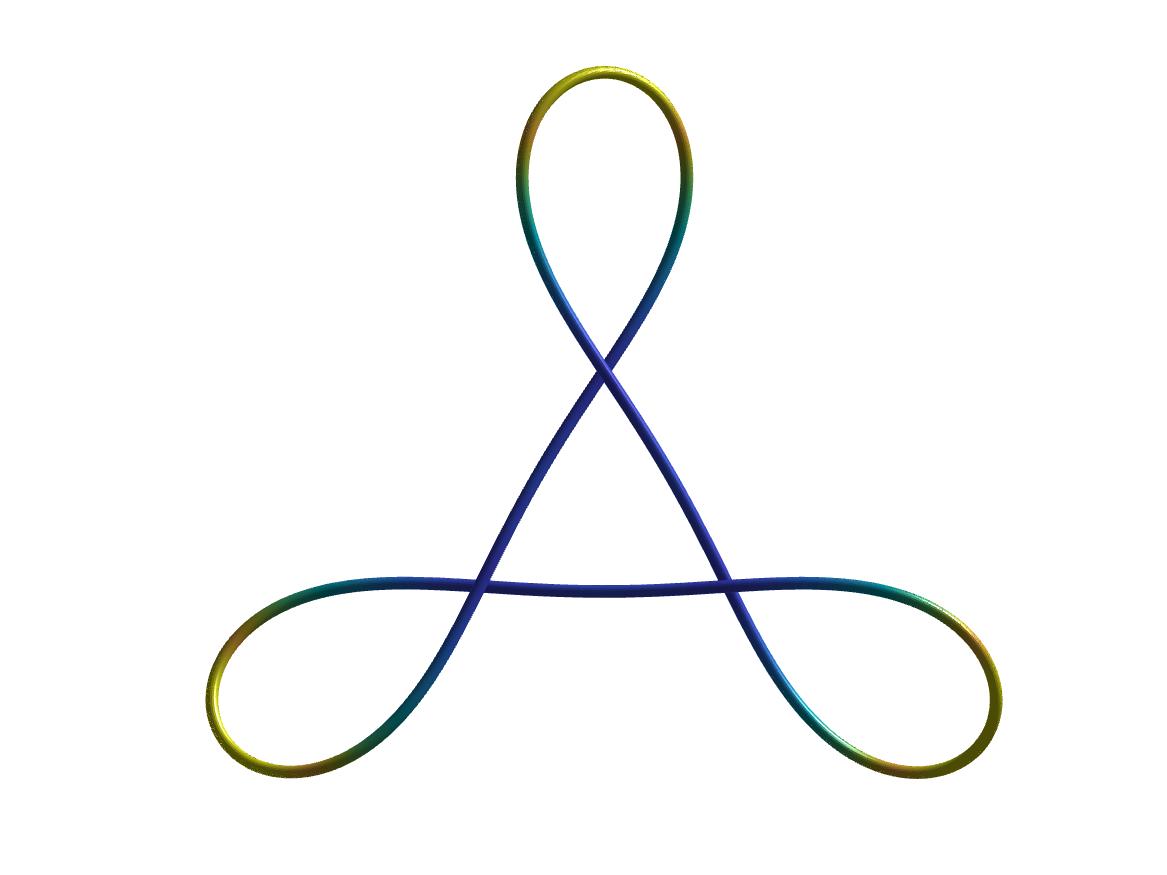}{200}
\bild{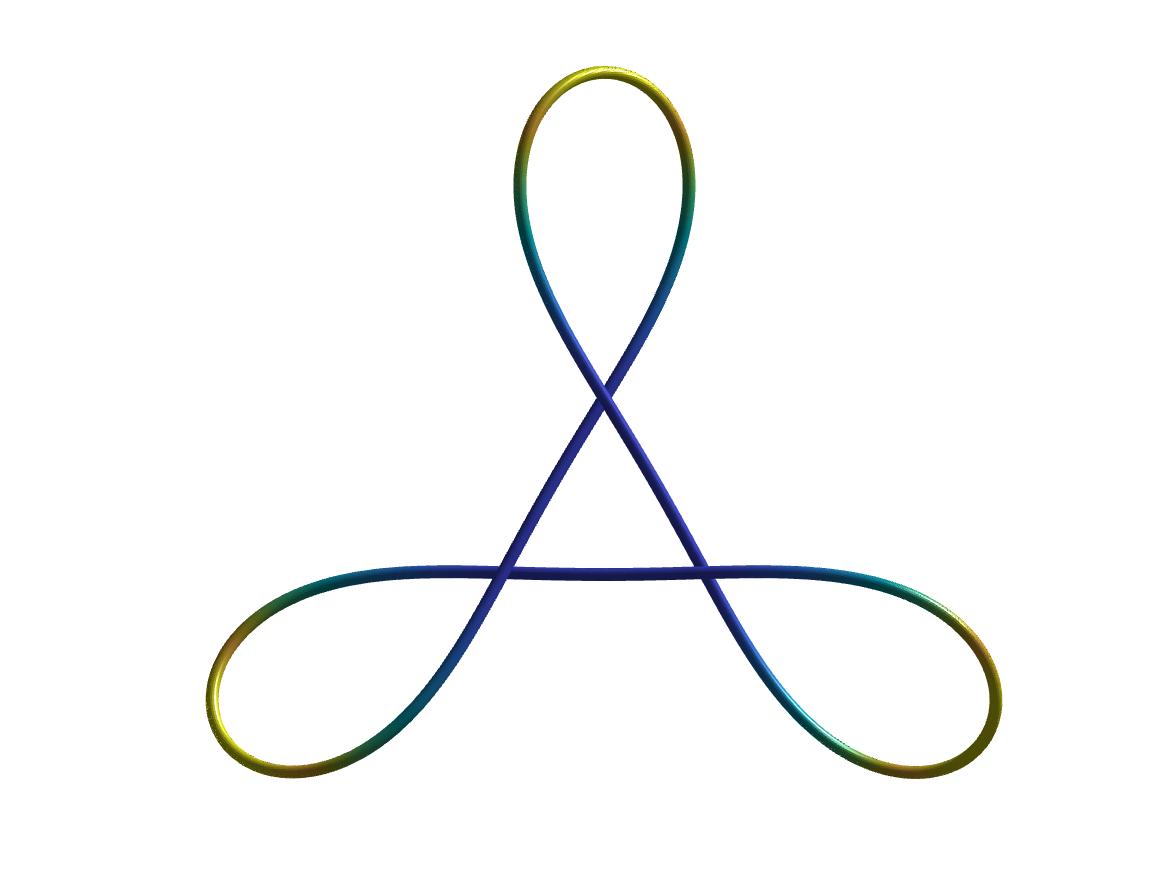}{300}

\bild{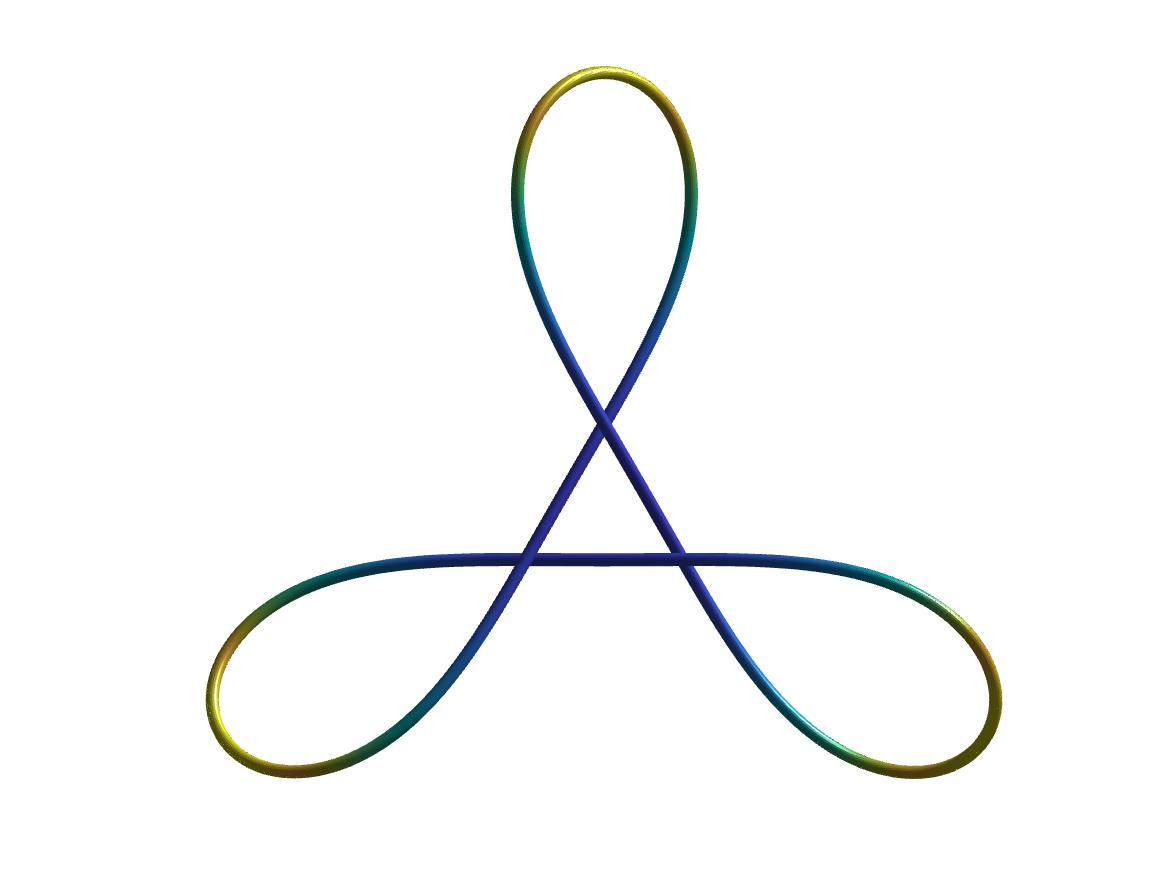}{400}
\bild{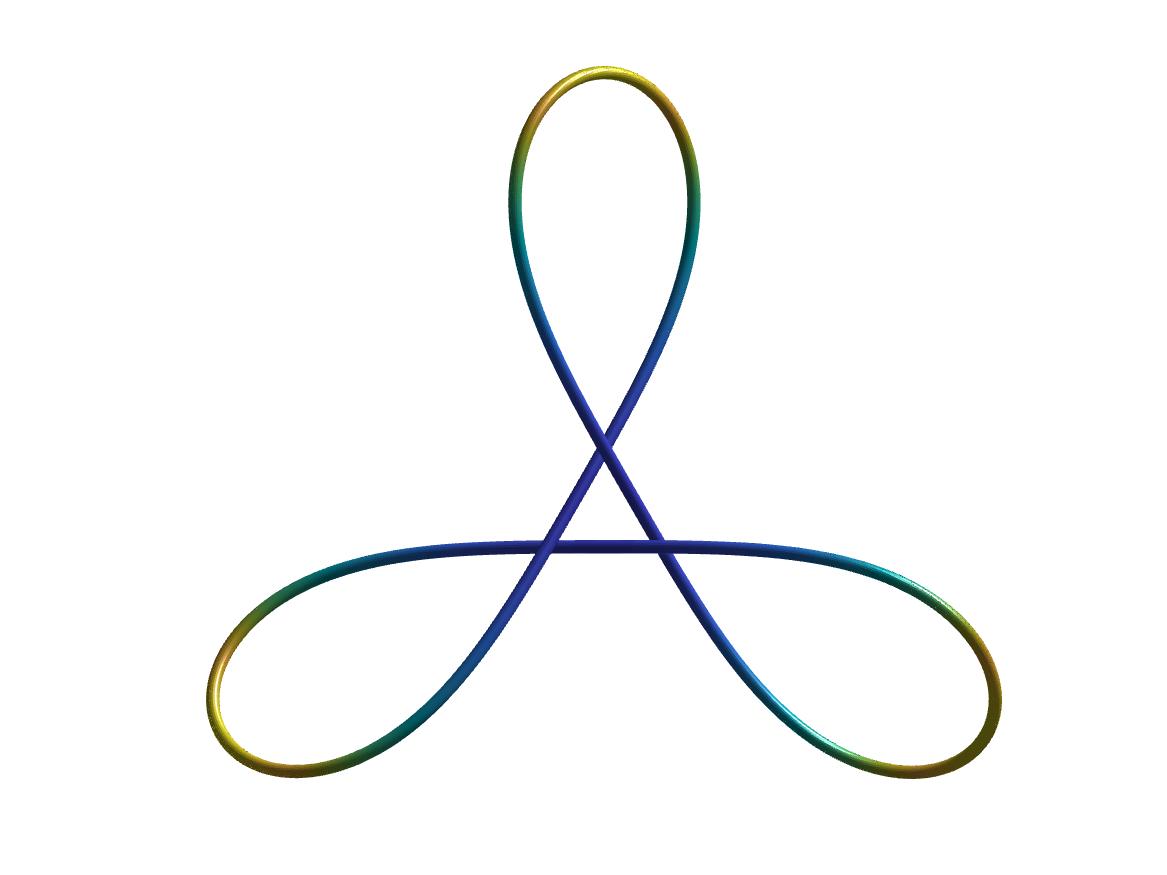}{500}
\bild{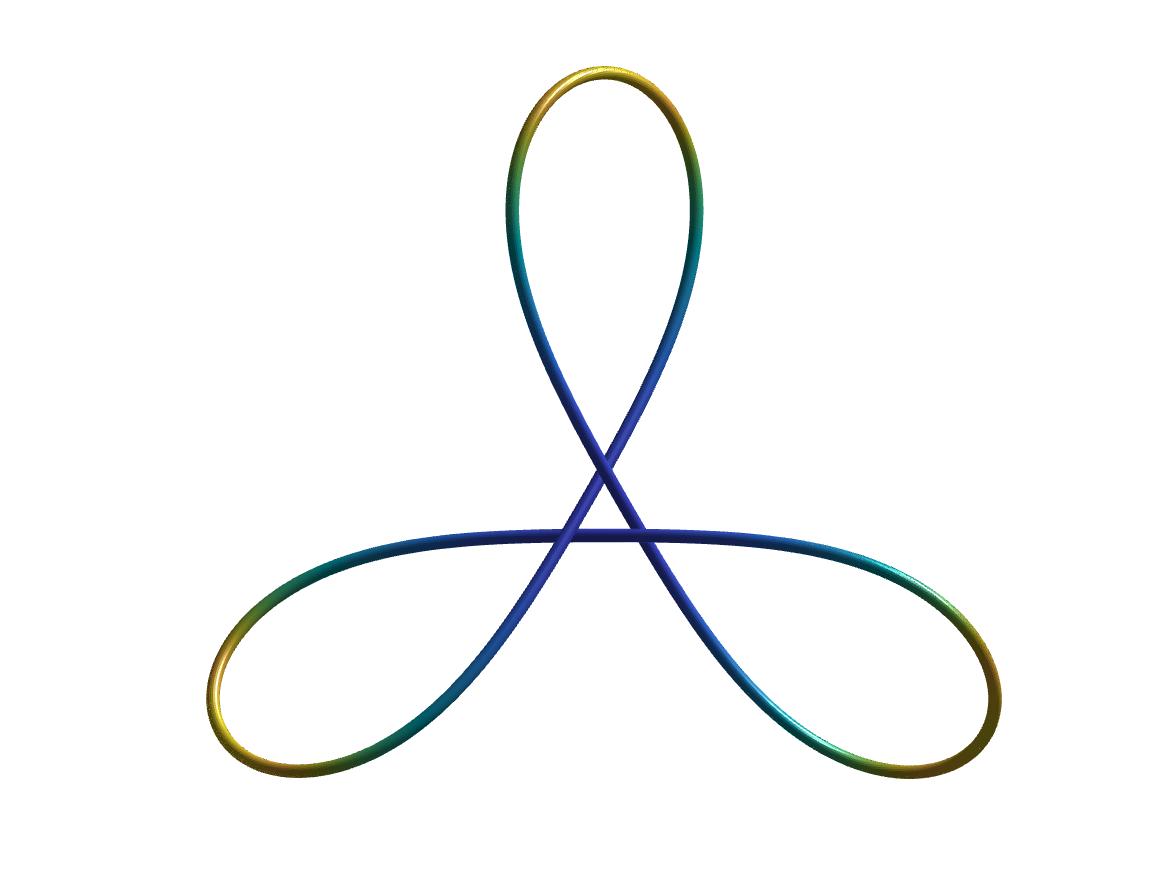}{600}
\bild{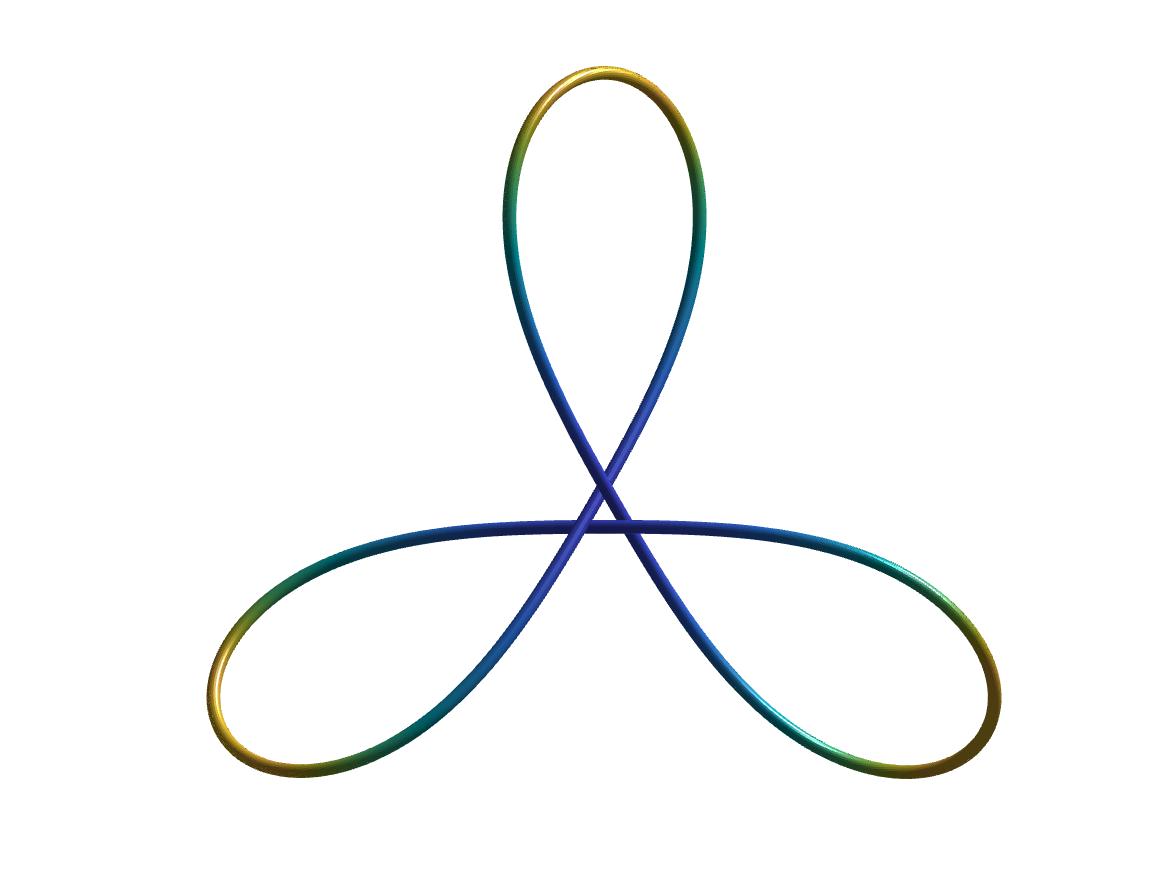}{700}

\bild{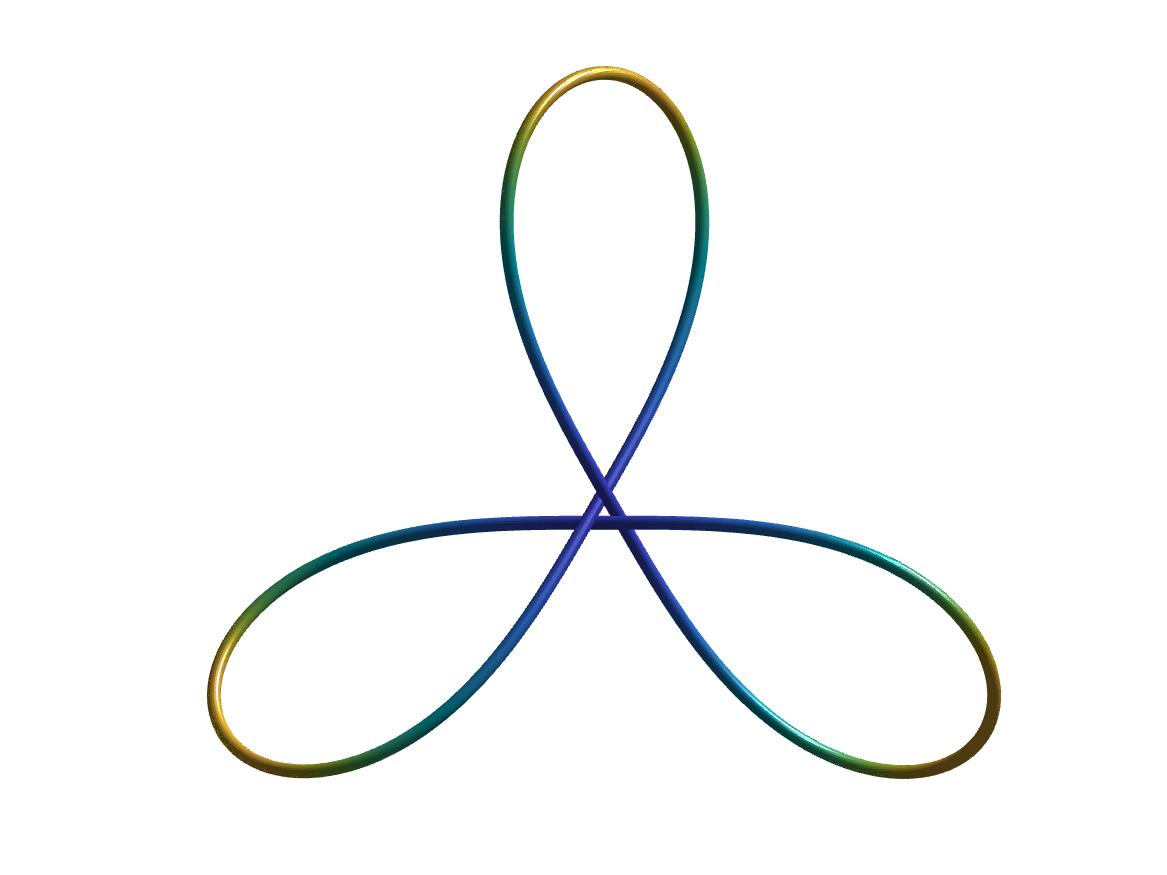}{800}
\bild{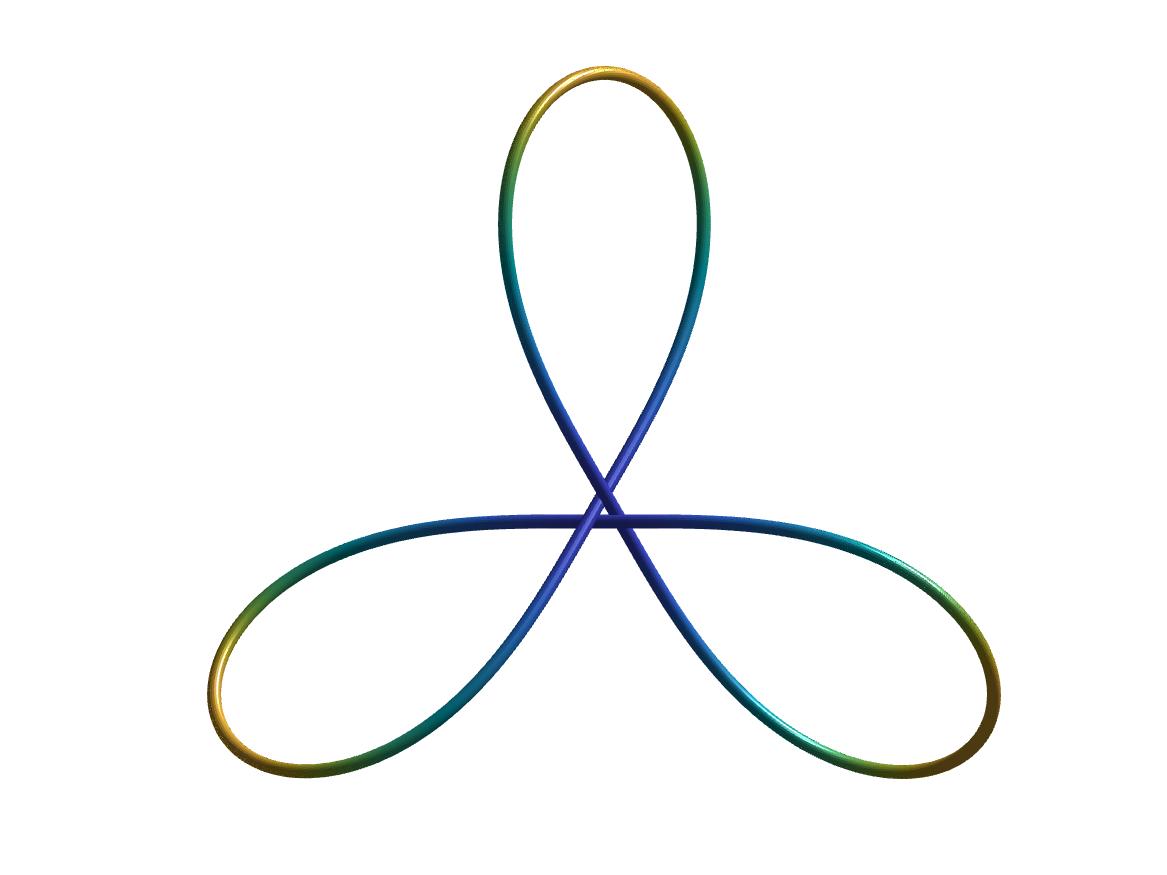}{900}
\bild{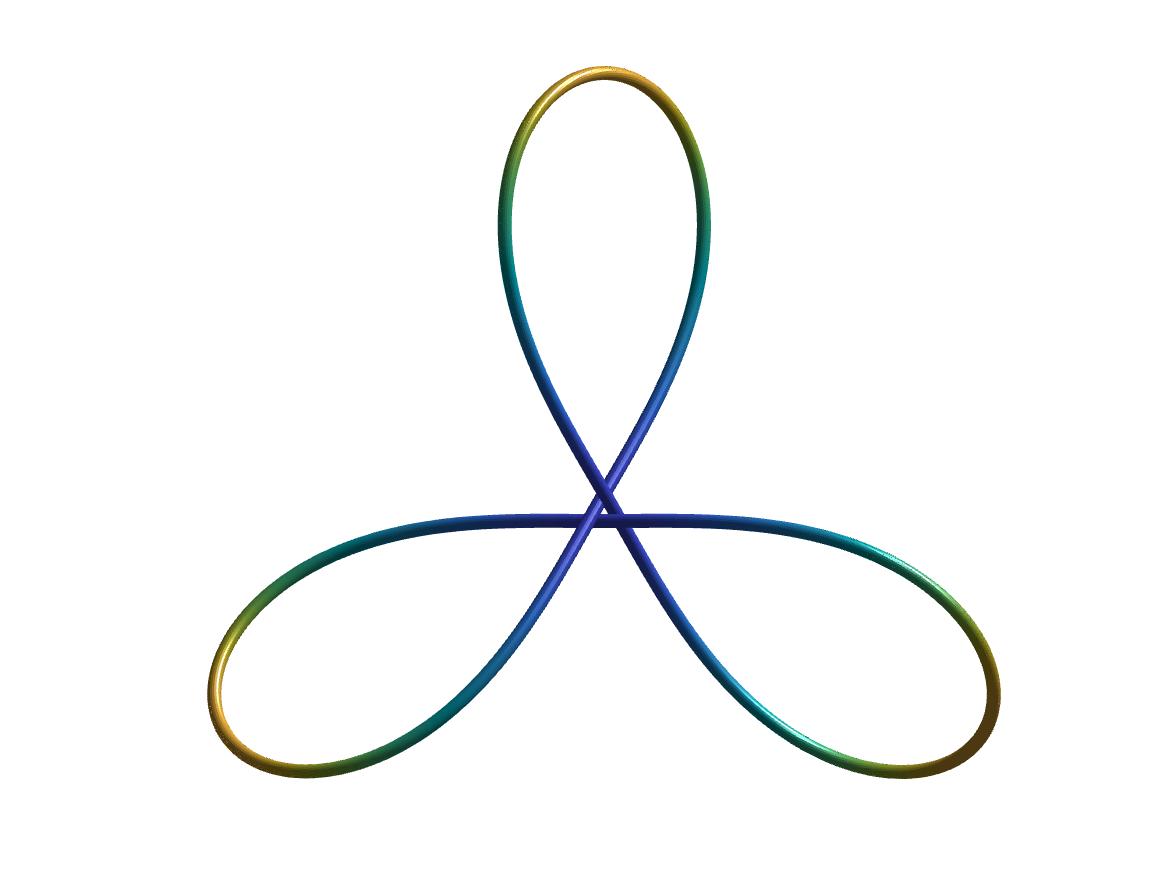}{1000}
\bild{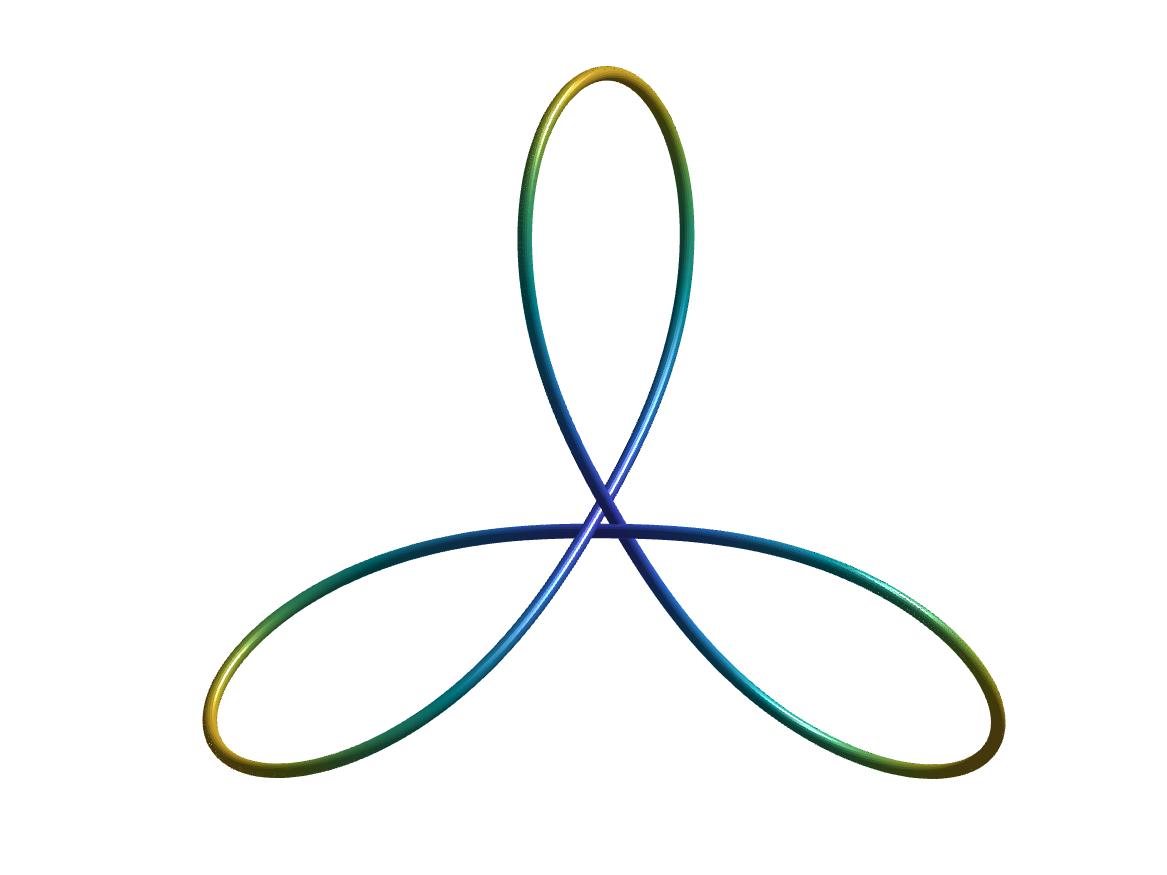}{23000}

\vspace{3mm}

\includegraphics[scale=.475]{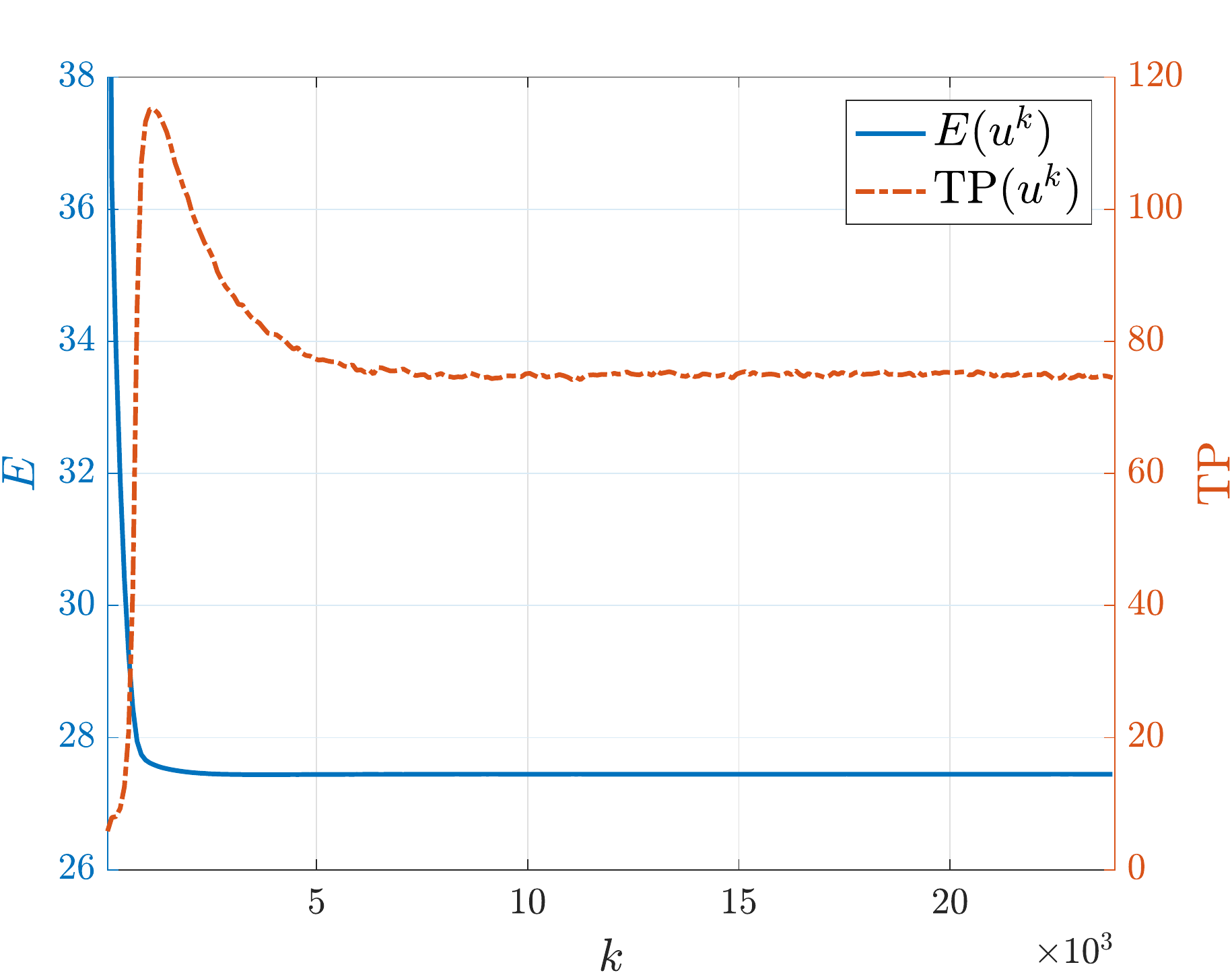}

\caption{
Numerical experiment with an unknot configuration defined in 
Section~\ref{sec:unkn}. Top: snapshots after $k=1,100,\dots,1000,23000$
iterations; bottom: monotonically decreasing energy curve.
The discrete gradient flow starting from an unknotted
symmetric polygon does not converge to the global
minimizer, i.e., the round circle, but to a configuration
which, being stable under small random perturbations, is likely 
to be a local minimizer.}\label{fig:unkn(a)}
\end{figure}

%

\subsection{Trefoil}\label{sec:tref}

In our second example we highlight the impact of symmetry to the evolution.
The fact that a curve belonging to the trefoil knot class
converges to a doubly covered circle (as predicted in~\cite{GRvdM})
has already been observed for the discretized $L^{2}$ gradient flow in~\cite{BRR}.

Here we start with an embedded curve belonging to the
trefoil knot class which is close to the three times covered circle.
The initial curve is obtained by discretizing and rescaling the curve
\begin{equation}\label{eq:trefoil}
x \mapsto 
\begin{pmatrix}
(2+0.1\cos(4\pi x))\cos(6\pi x) \\ (2+0.1\cos(4\pi x))\sin(6\pi x) 
\\ 0.1\sin(4\pi x)
\end{pmatrix}.
\end{equation}

The unperturbed discrete gradient flow (a) is depicted in the top part of
Figure~\ref{fig:tref(a)} and unfolds the curve but
seems to get stuck in a conformation approximating the shape
of two tangential circles tangentially meeting in an angle of $180$ degrees 
which might be a saddle point.

The perturbed discrete gradient flow (b) is depicted in the bottom part of
Figure~\ref{fig:tref(b)} and leaves that
state after some time and approaches the shape of the
elastic knot, i.e., the doubly covered (round) circle.
The final energy value is quite close to the analytically predicted
threshold of $\tfrac\kappa2\cdot50\cdot\br{\frac{4\pi}{50}}^{2} = \frac{(4\pi)^{2}}{10}\approx15.79$, cf.~\cite{GRvdM}.

The two evolutions are almost identical for the first 20,000 iteration
steps. Here the snapshots shown in the top and bottom parts of Figure~\ref{fig:tref(a)}
show essentially the same configurations from different viewpoints, and
the corresponding parts of the energy plot coincide
as shown in Figure~\ref{fig:tref}.


\renewcommand{\bild}[2]{\fbox{\includegraphics[scale=.1,trim=220 150 180 140,clip]{#1}}\makebox[0ex][r]{\tiny#2\ifnum#2>13000\hspace{19ex}\fi\ }\,\ignorespaces}

\begin{figure}

\bild{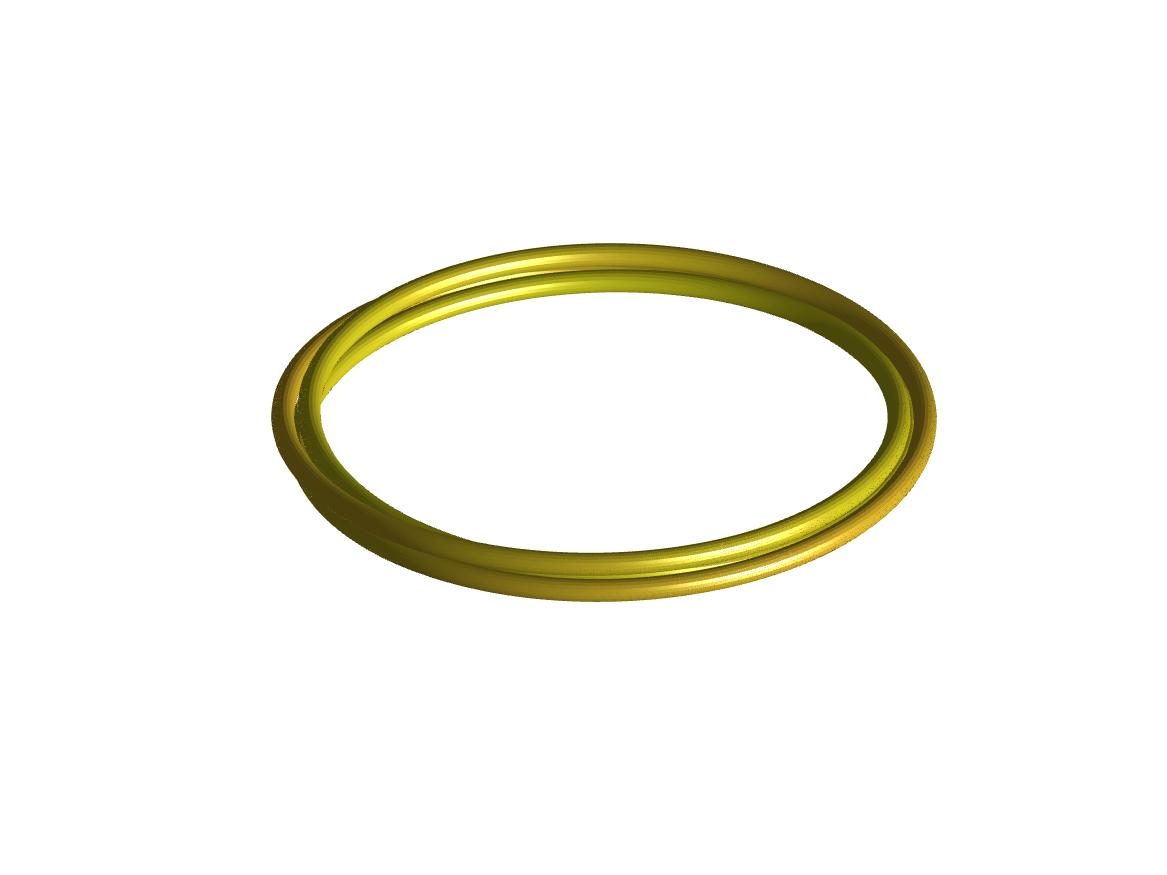}{1}
\bild{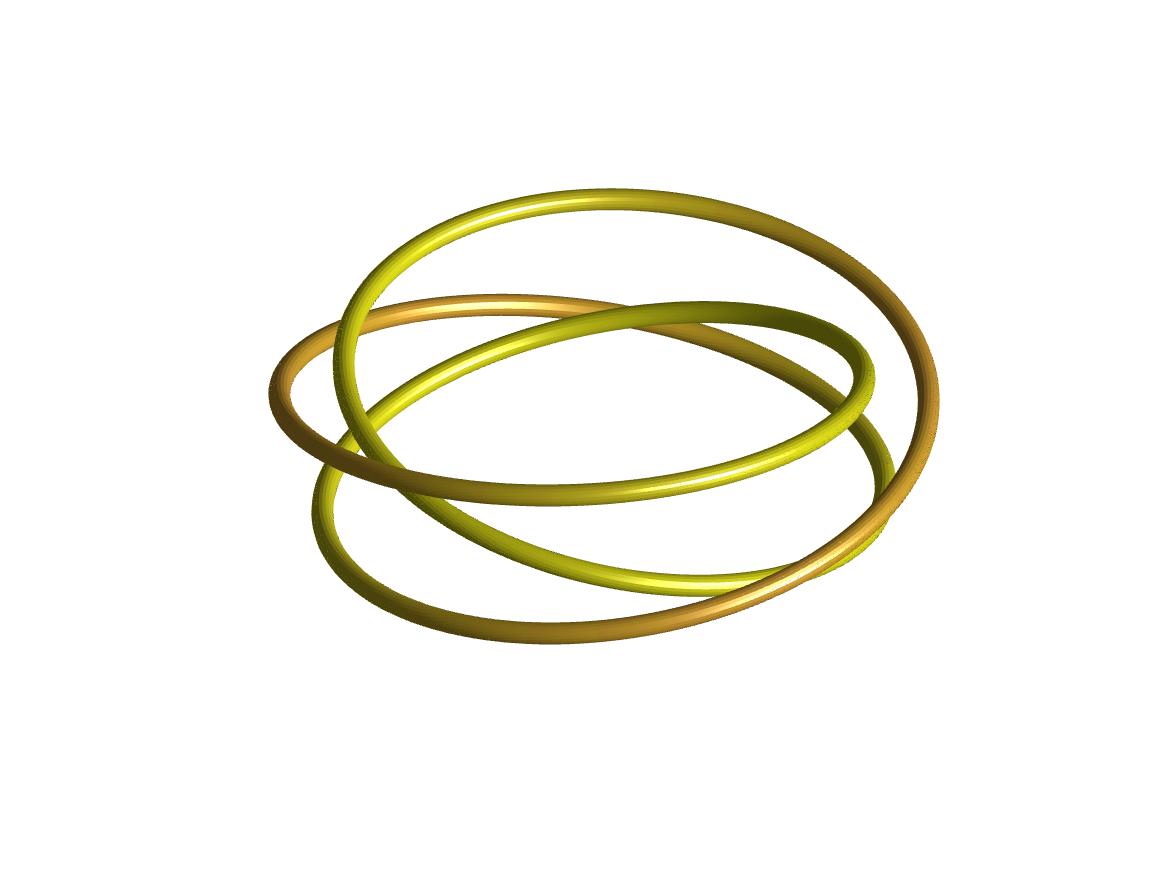}{1000}
\bild{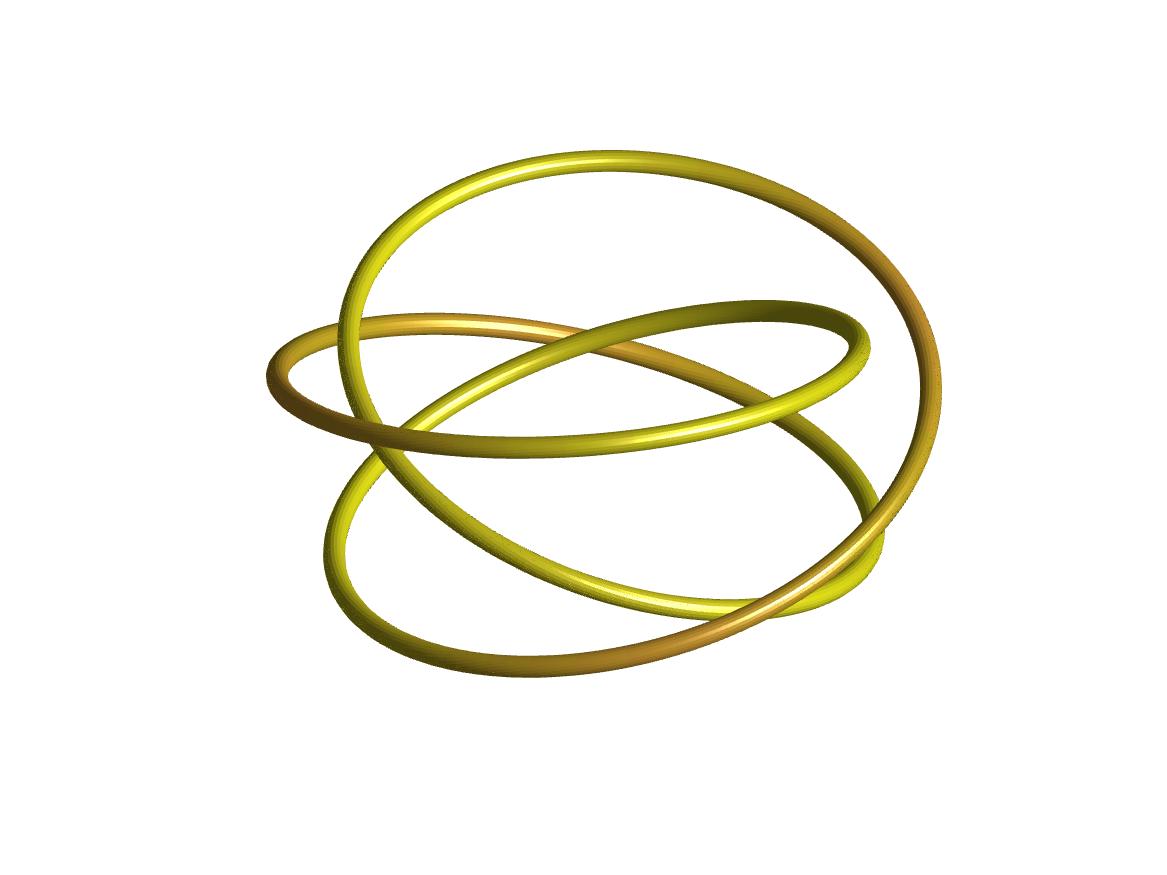}{2000}
\bild{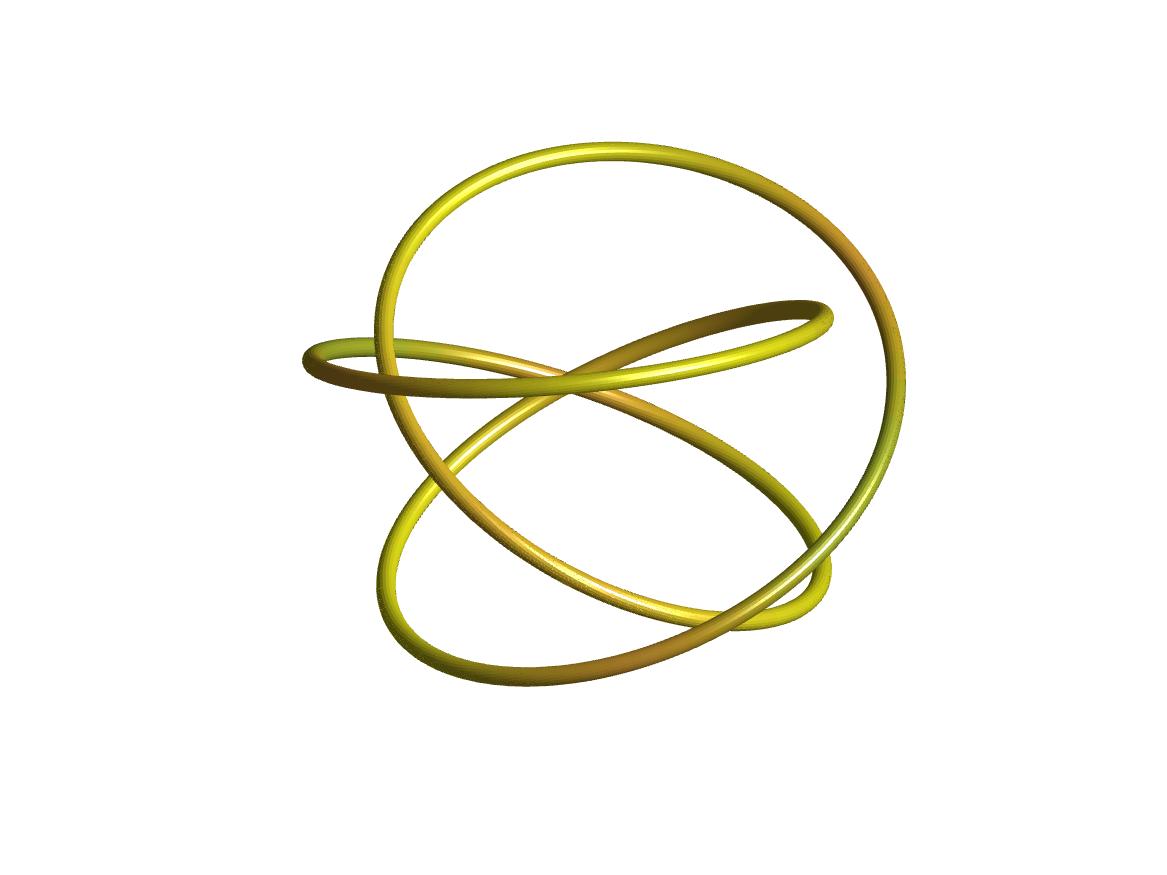}{3000}

\bild{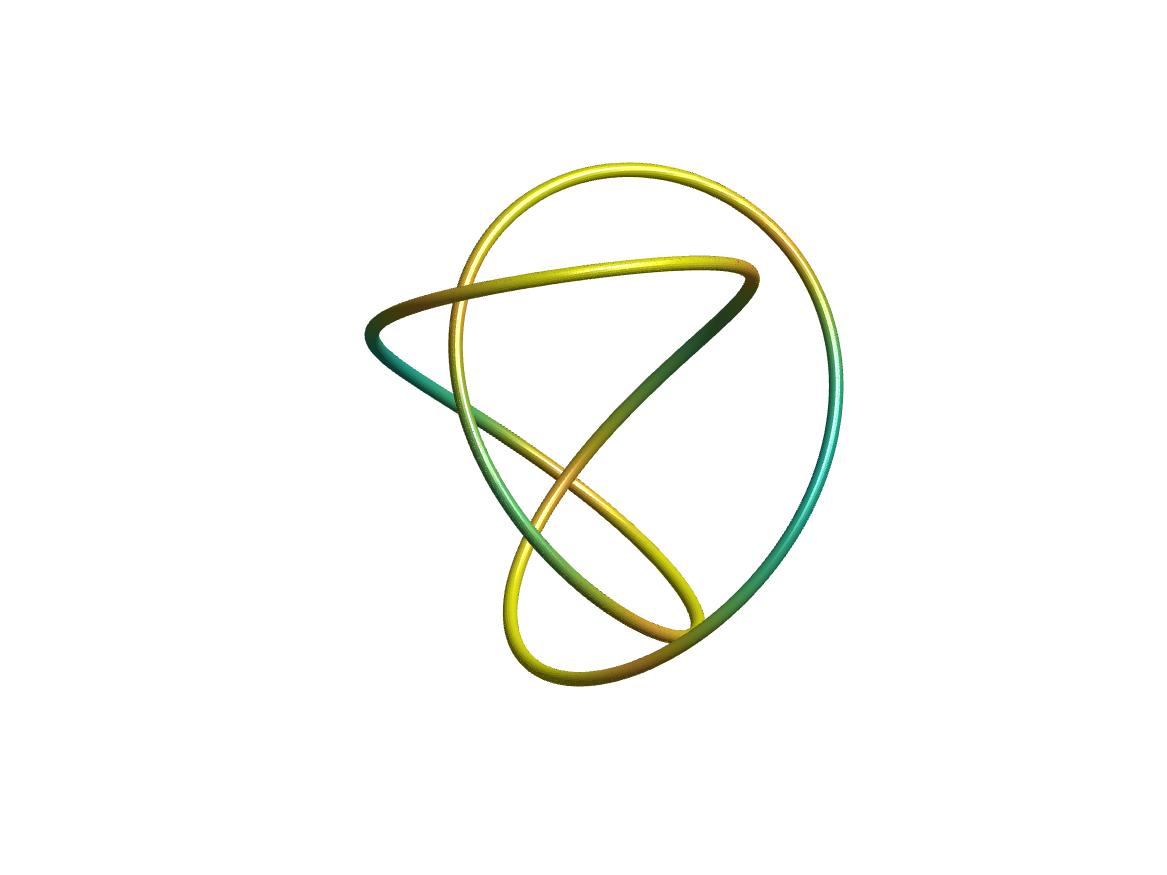}{5000}
\bild{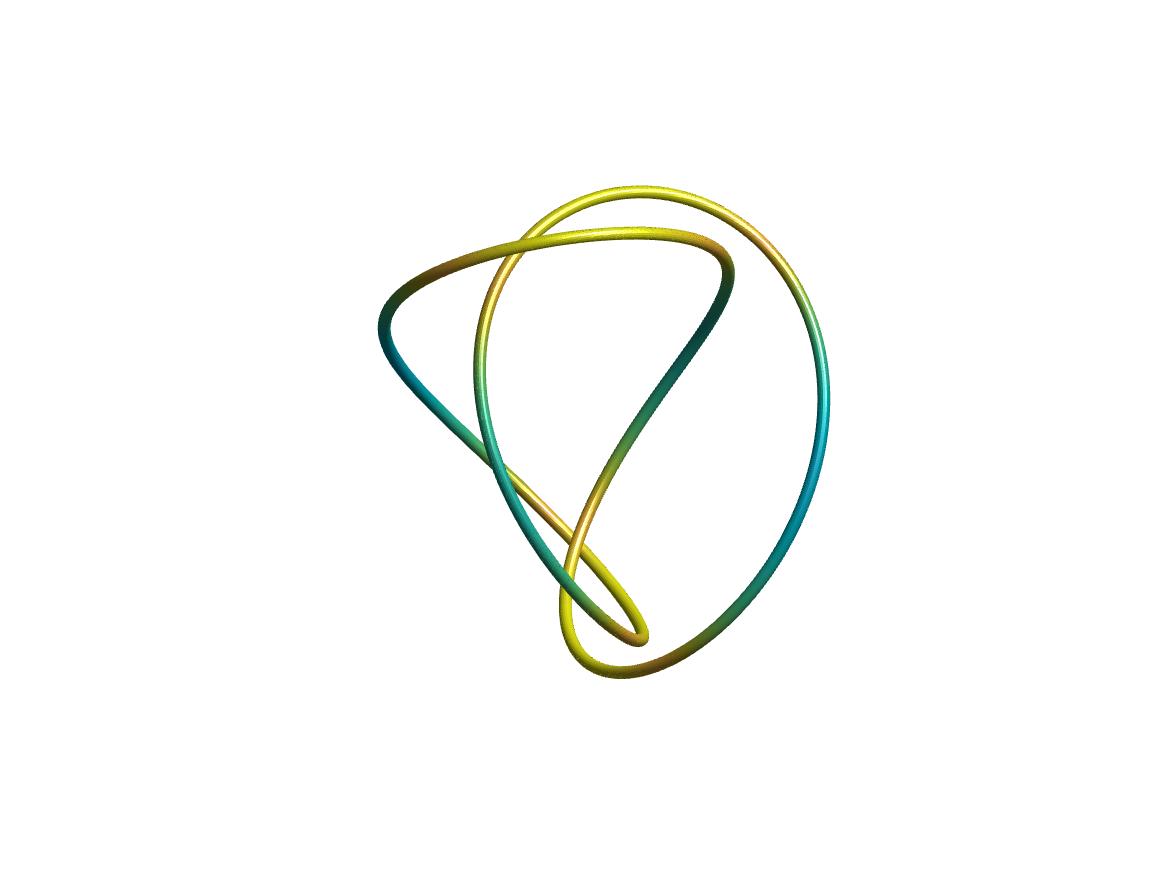}{6000}
\bild{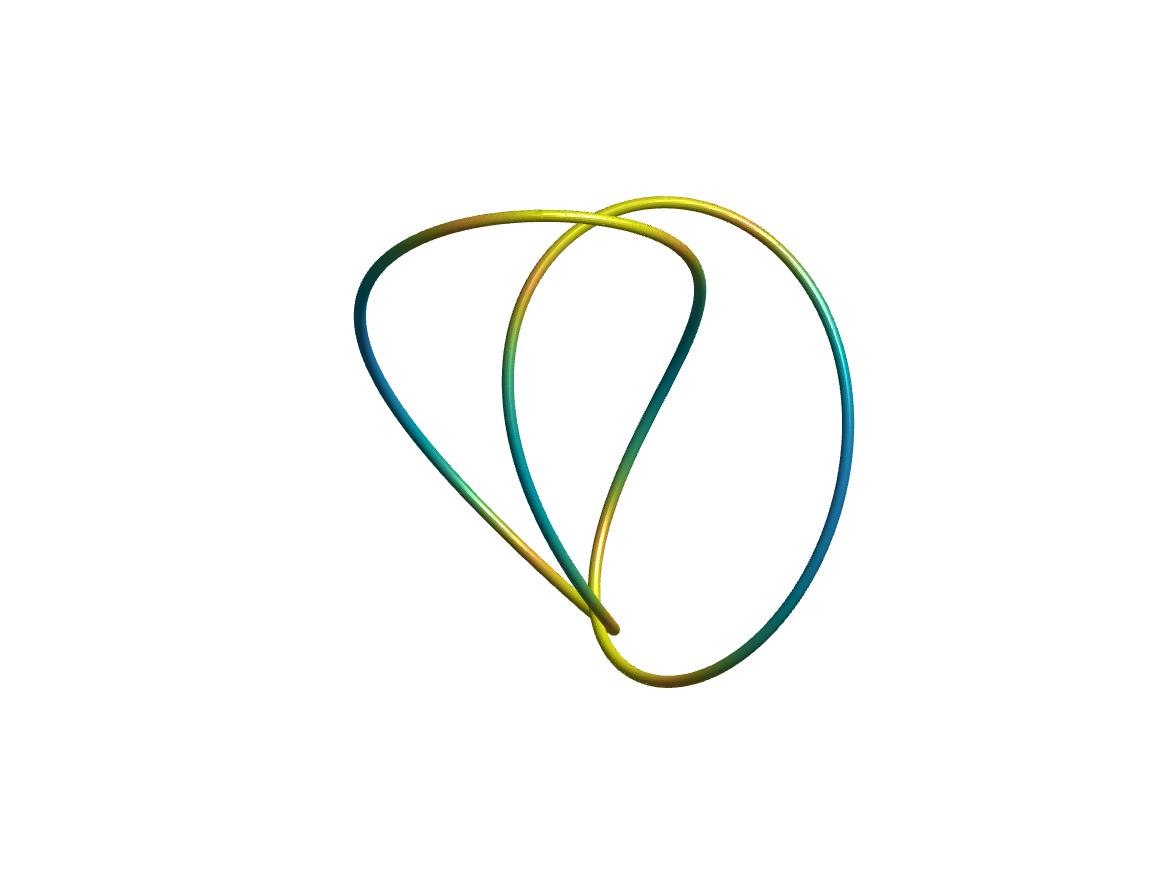}{10000}
\bild{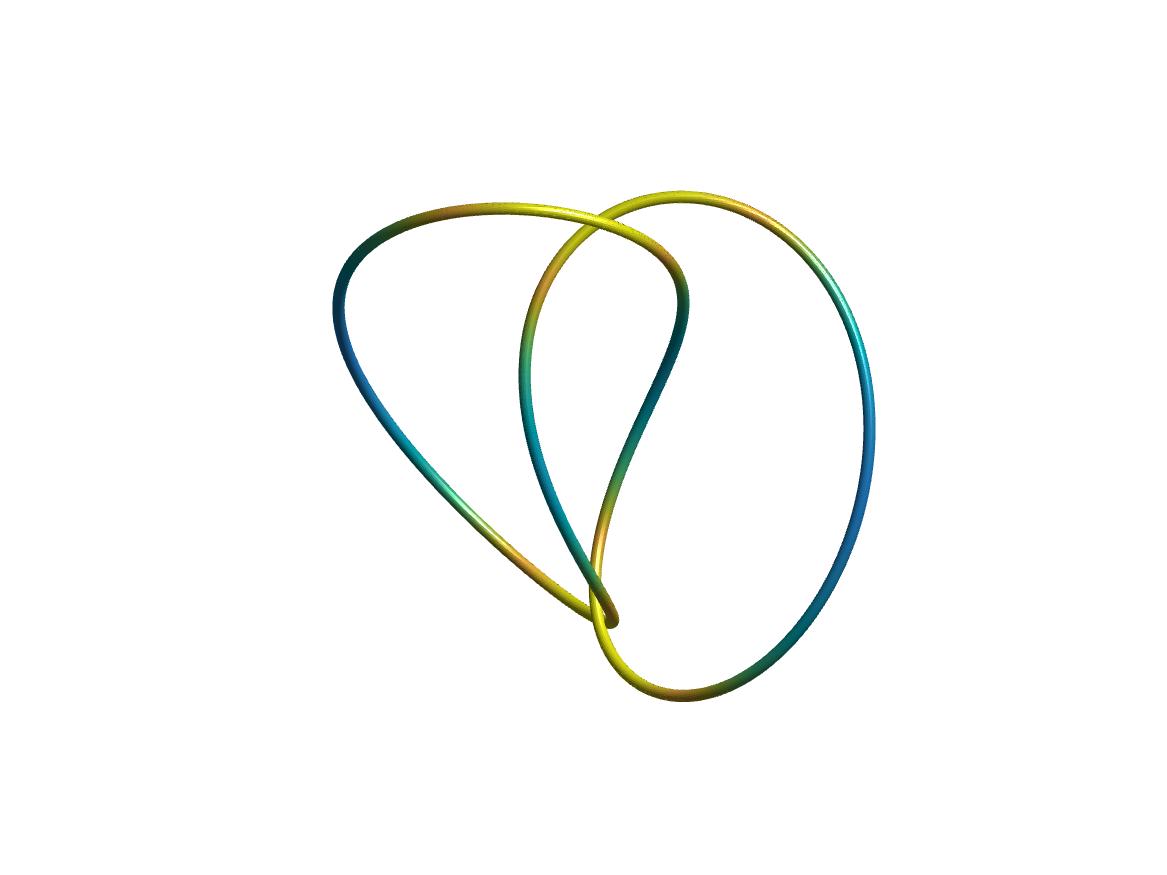}{12000}

\bild{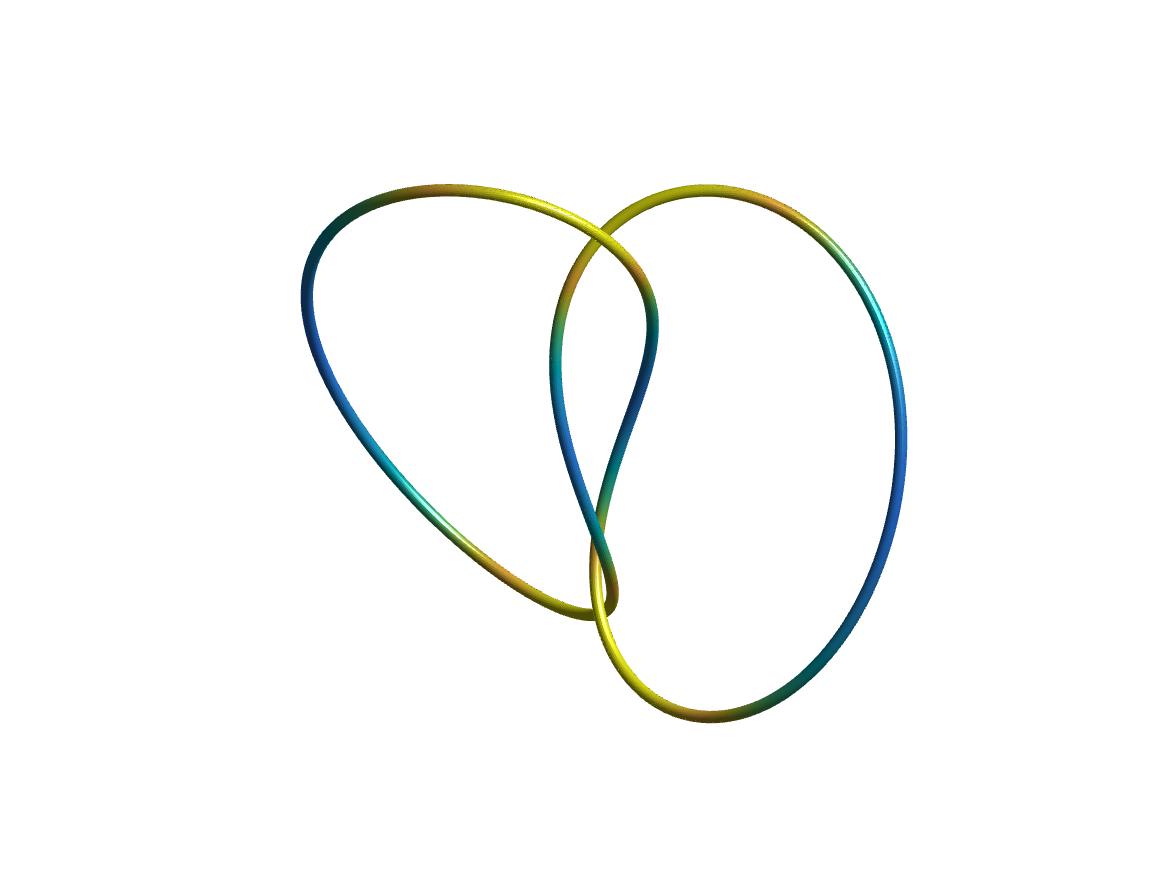}{14000}
\bild{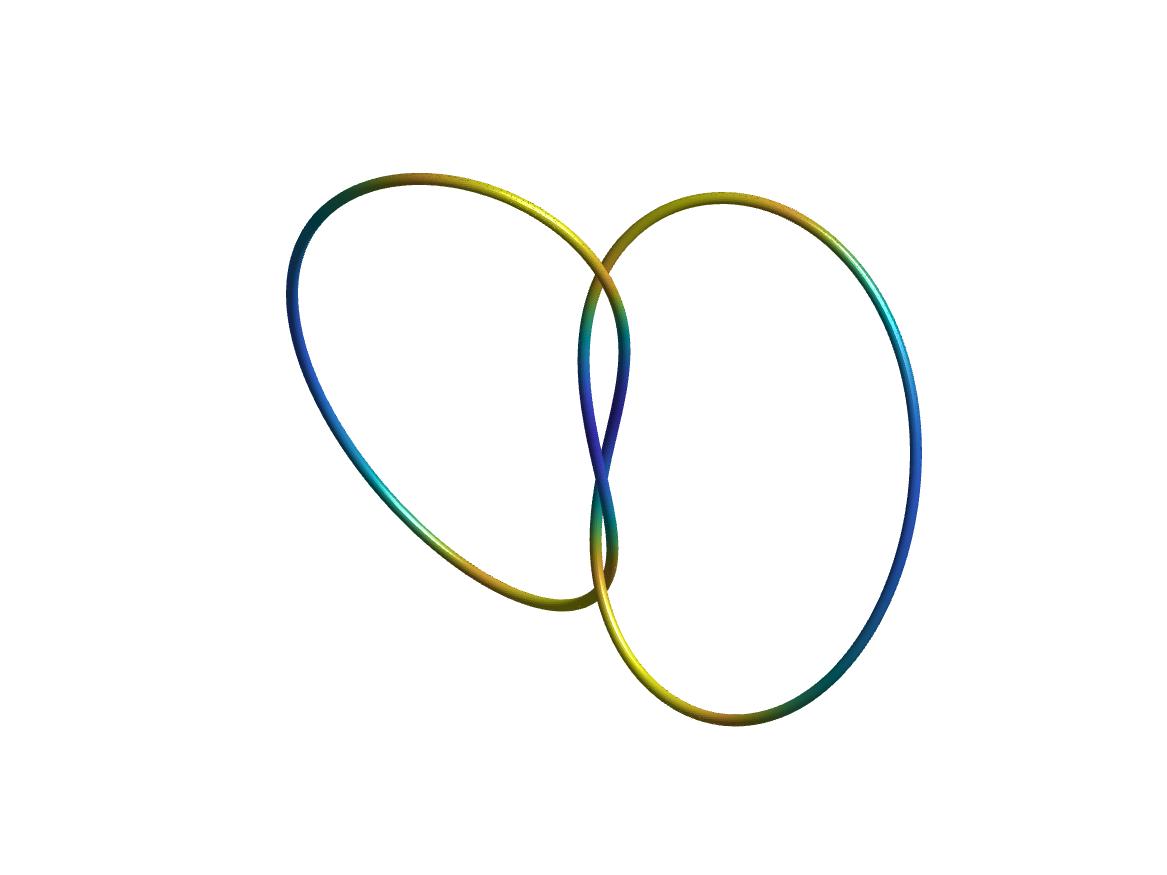}{15000}
\bild{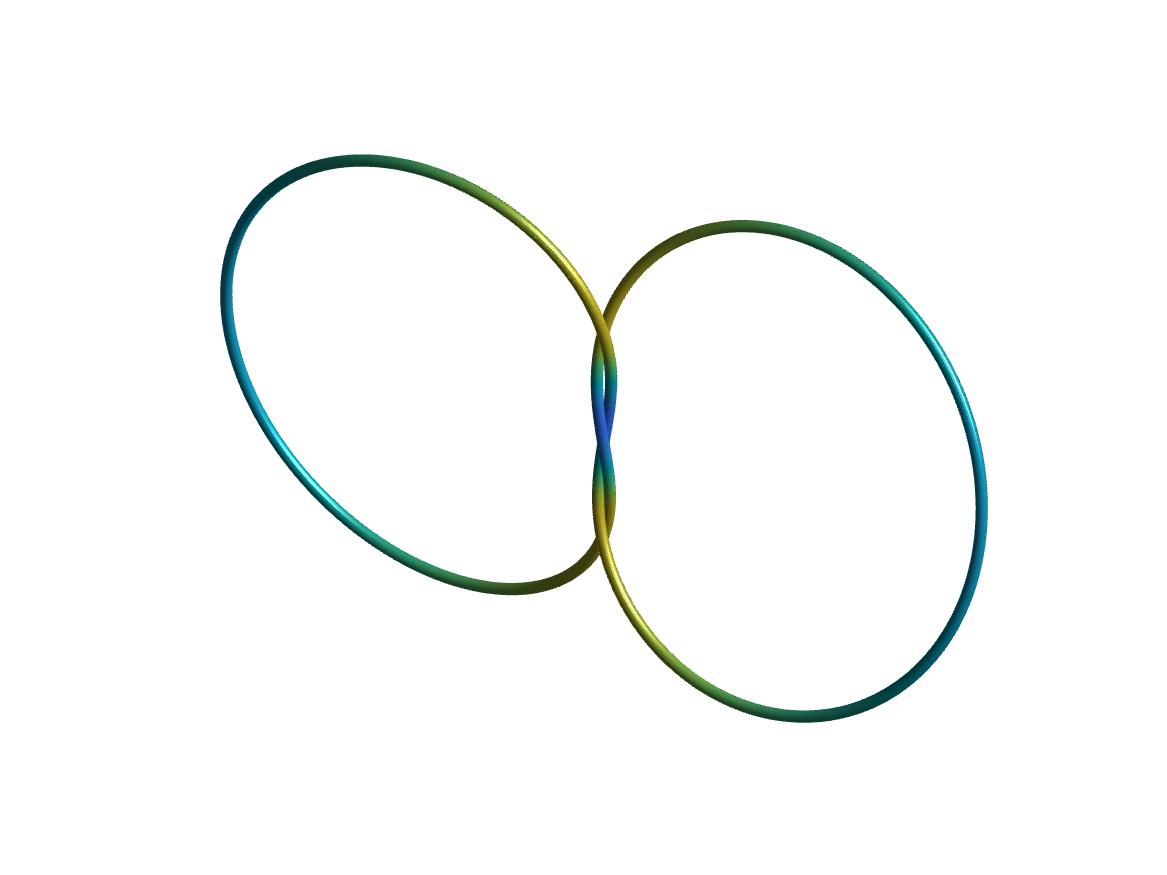}{20000}
\bild{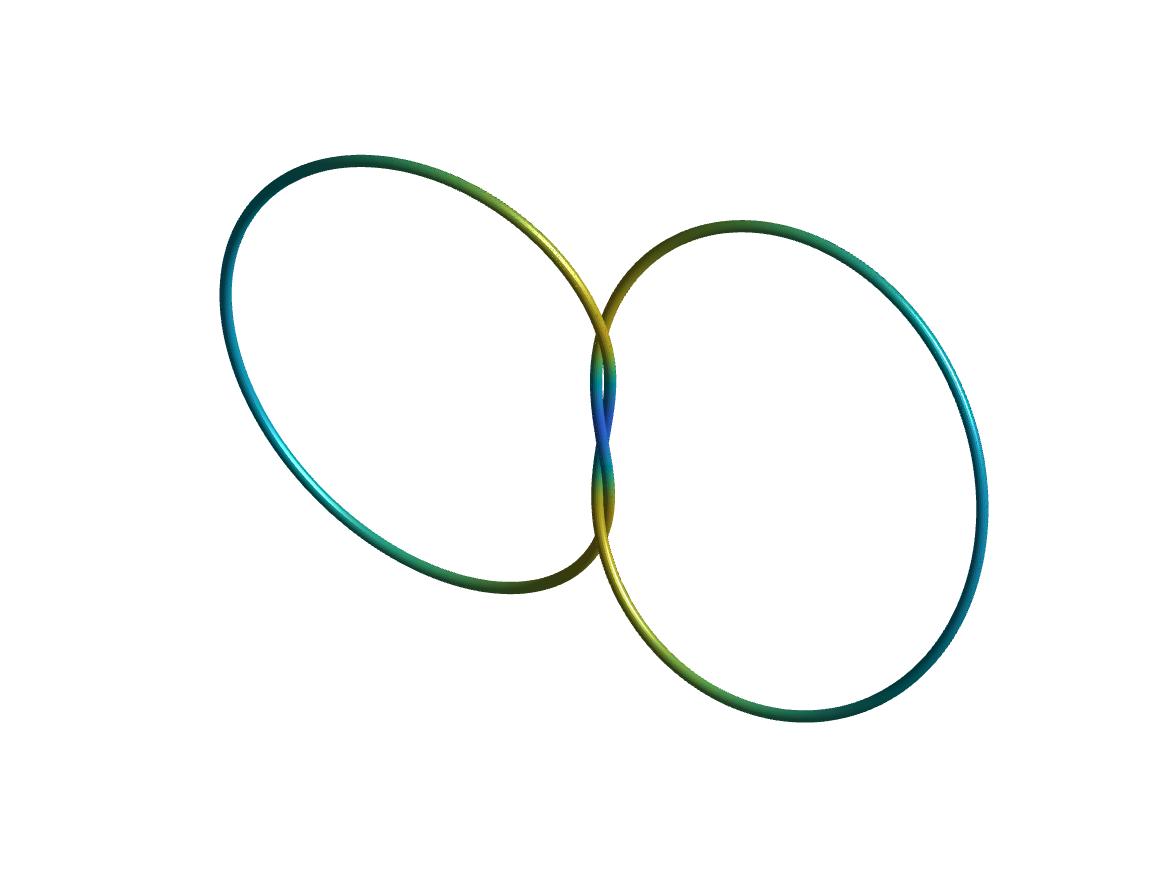}{80000}

%
%
\renewcommand{\bild}[2]{\fbox{\includegraphics[scale=.085,trim=150 100 110 70,clip]{#1}}\makebox[0ex][r]{\tiny#2\ }\,\ignorespaces}
%

\bigskip

\bild{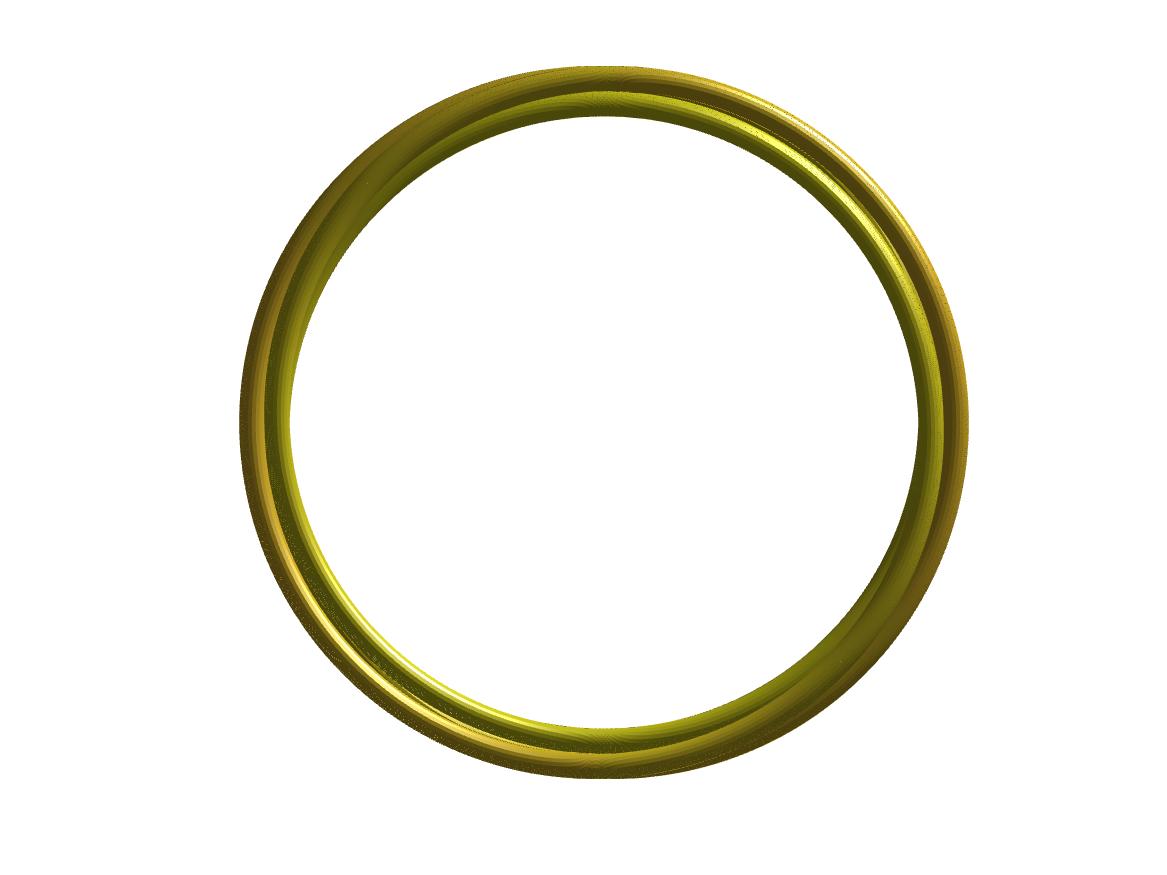}{1}
\bild{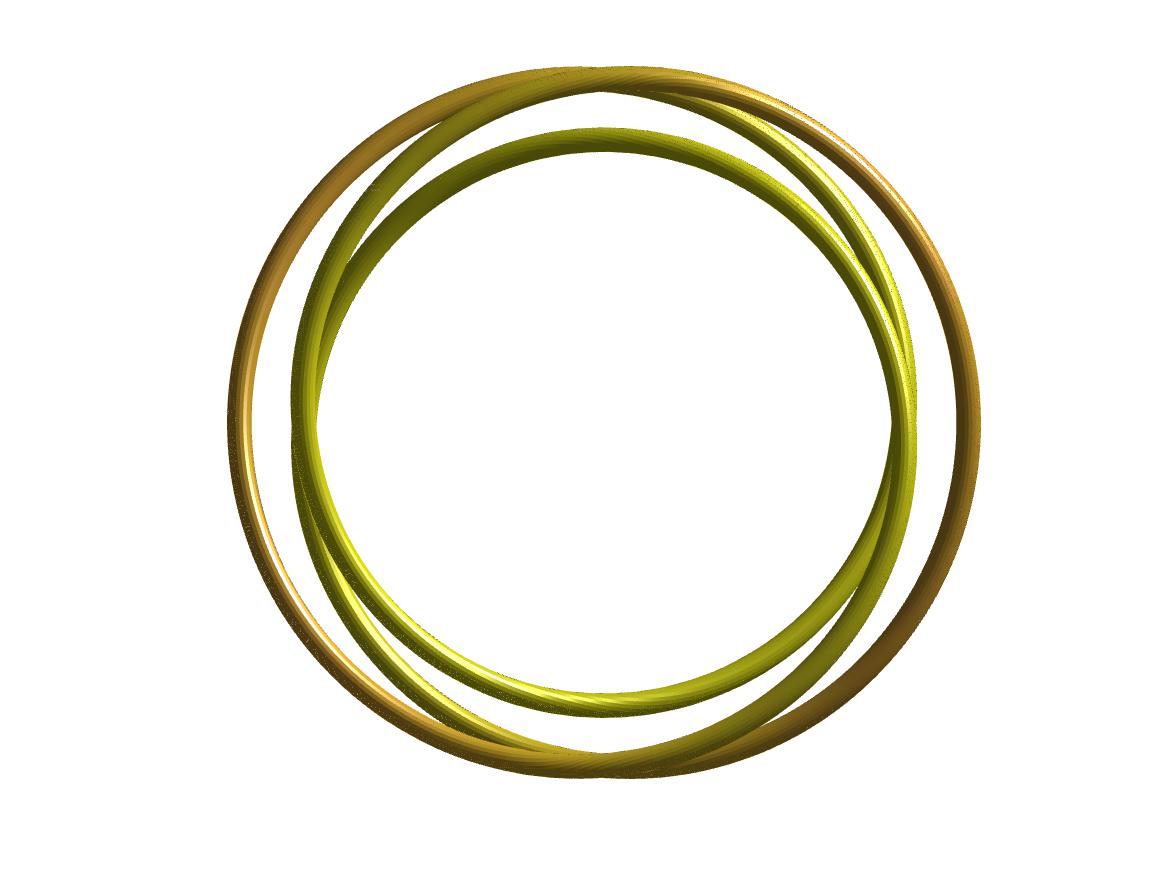}{1000}
\bild{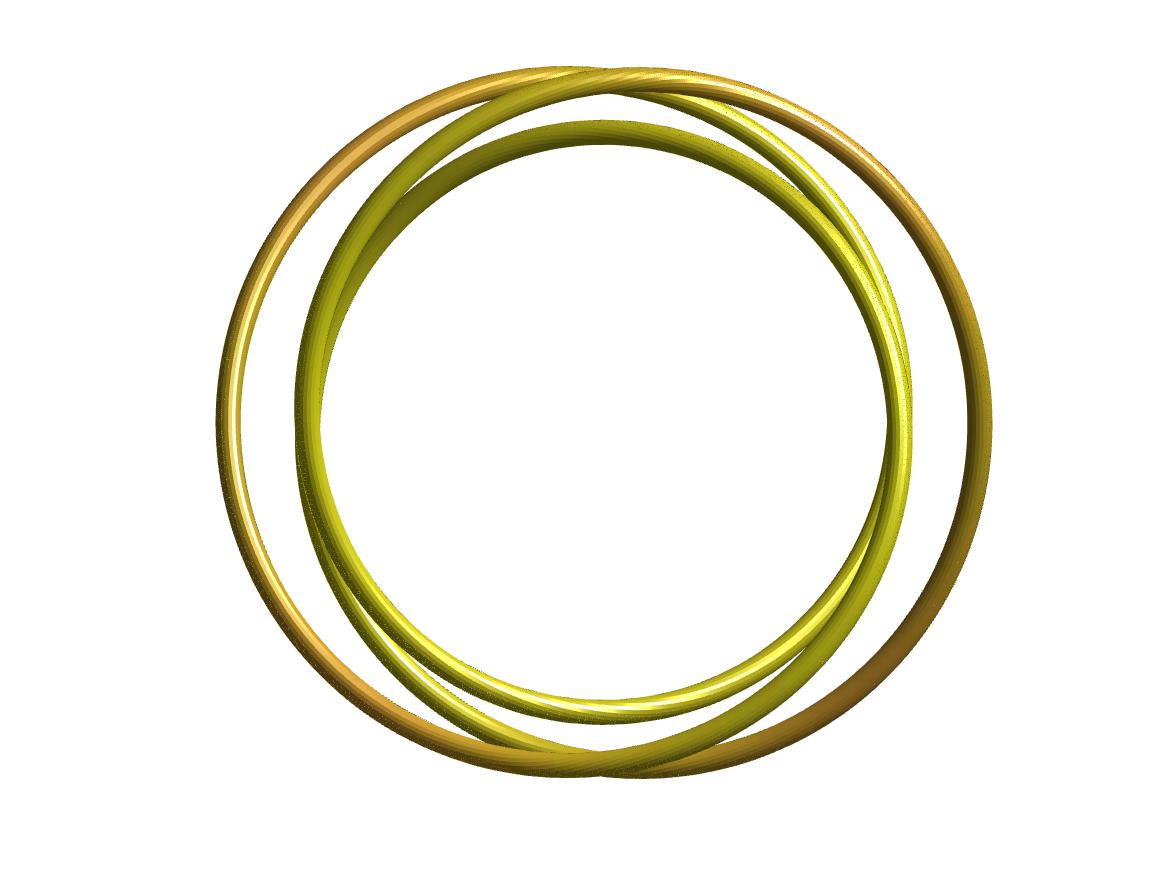}{2000}
\bild{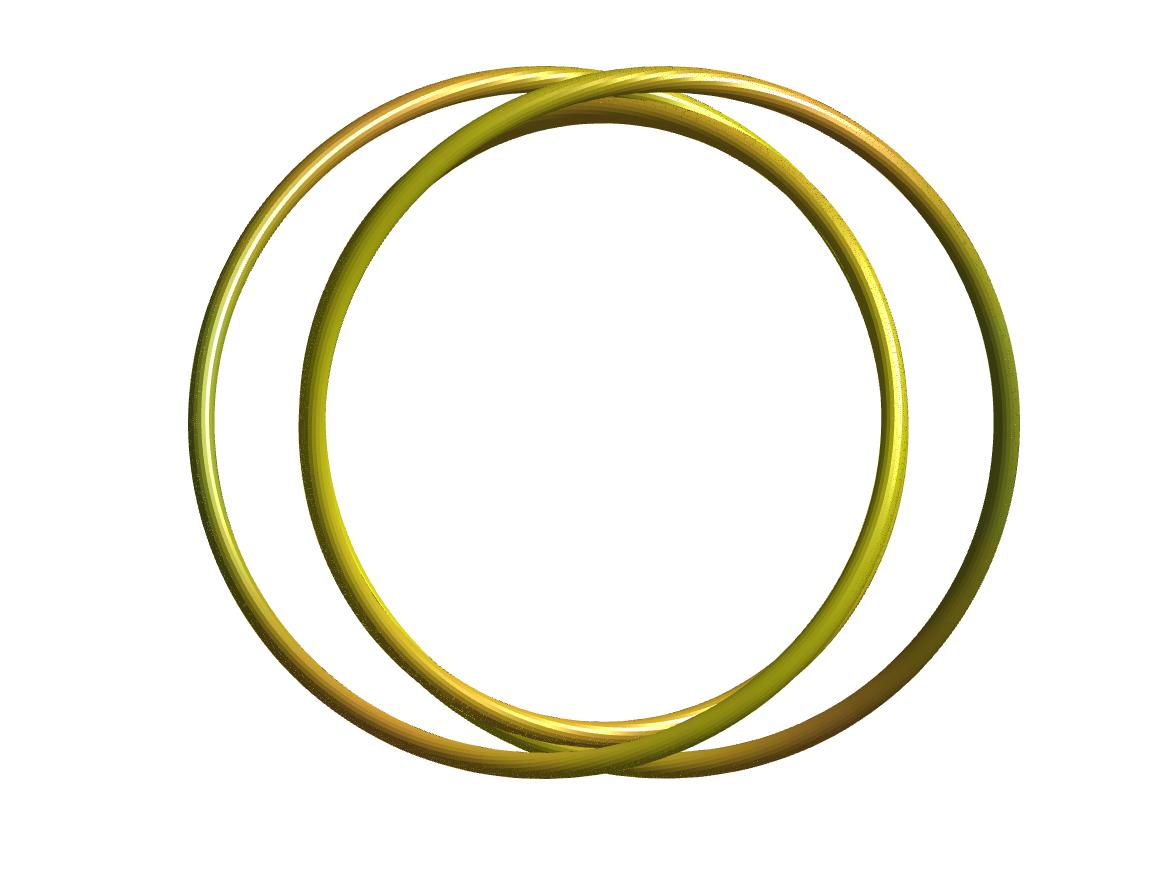}{3000}

\bild{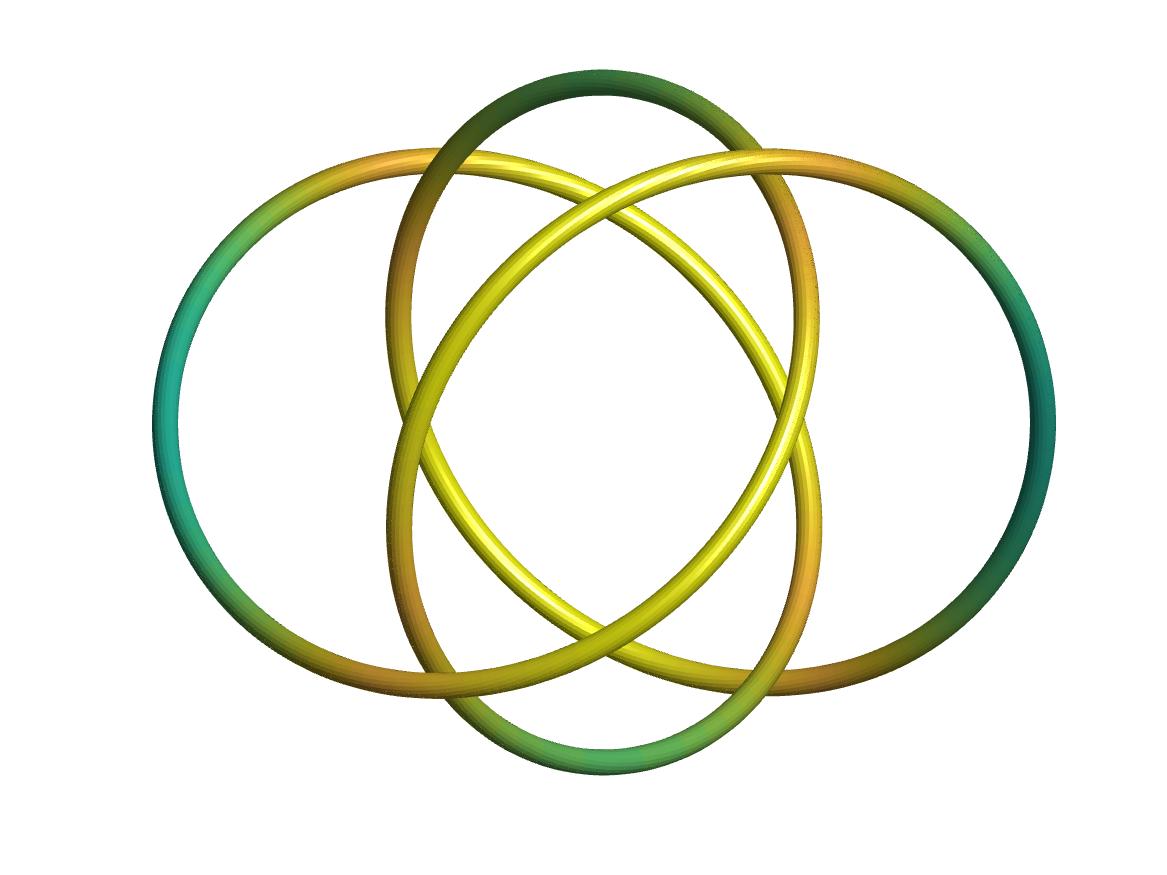}{5000}
\bild{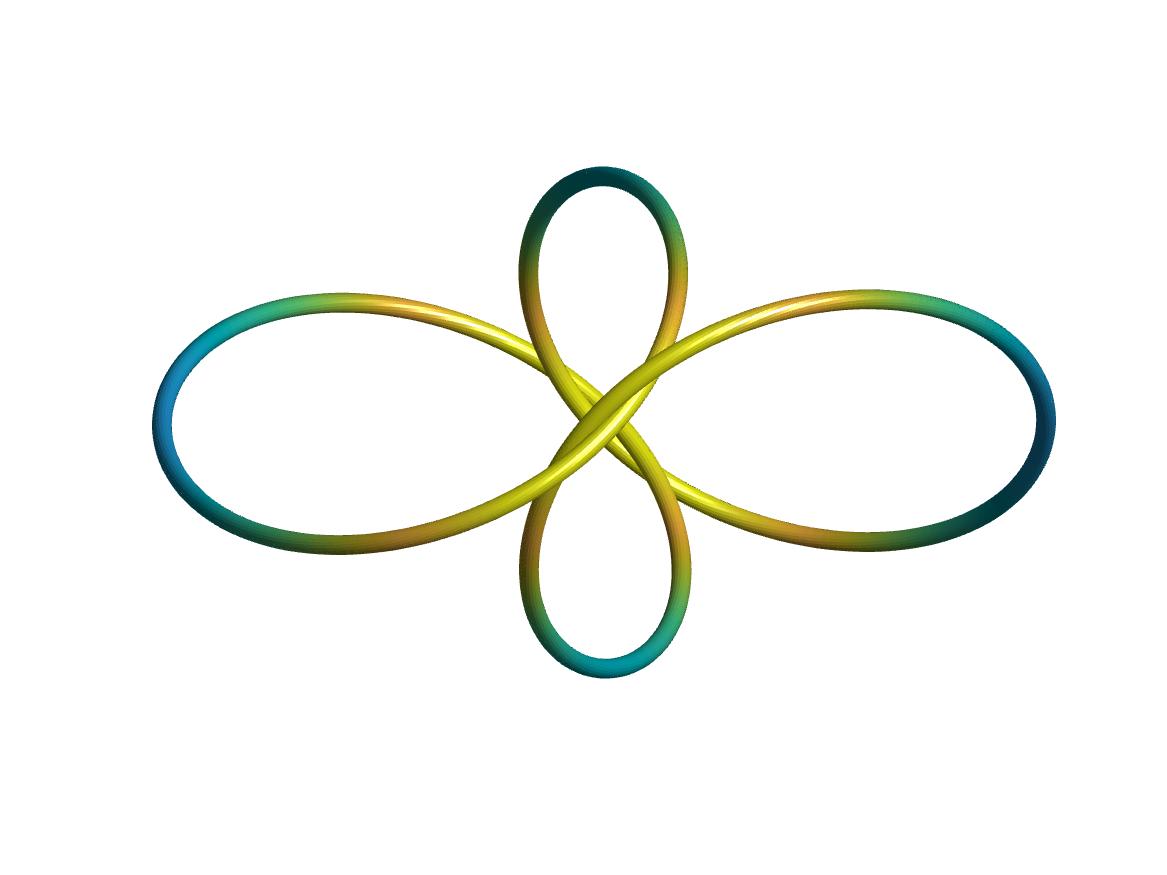}{10000}
\bild{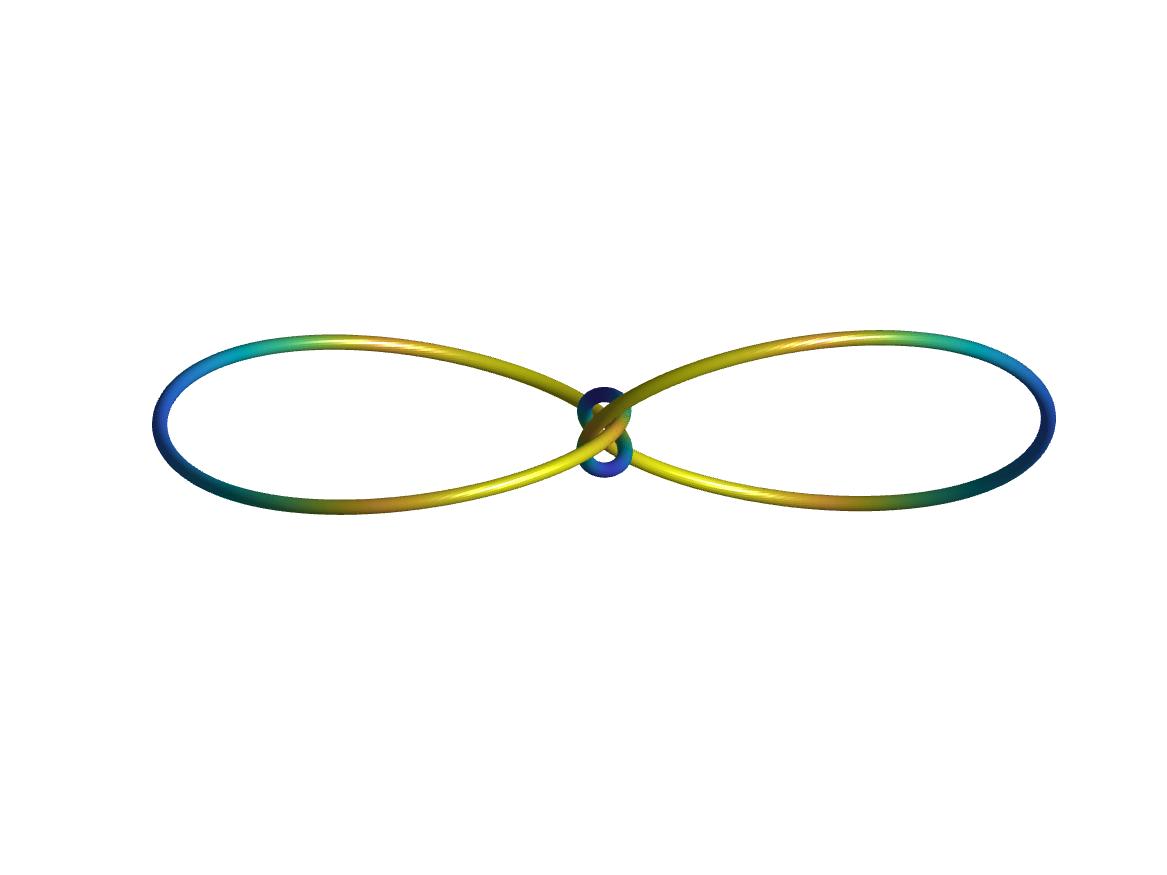}{15000}
\bild{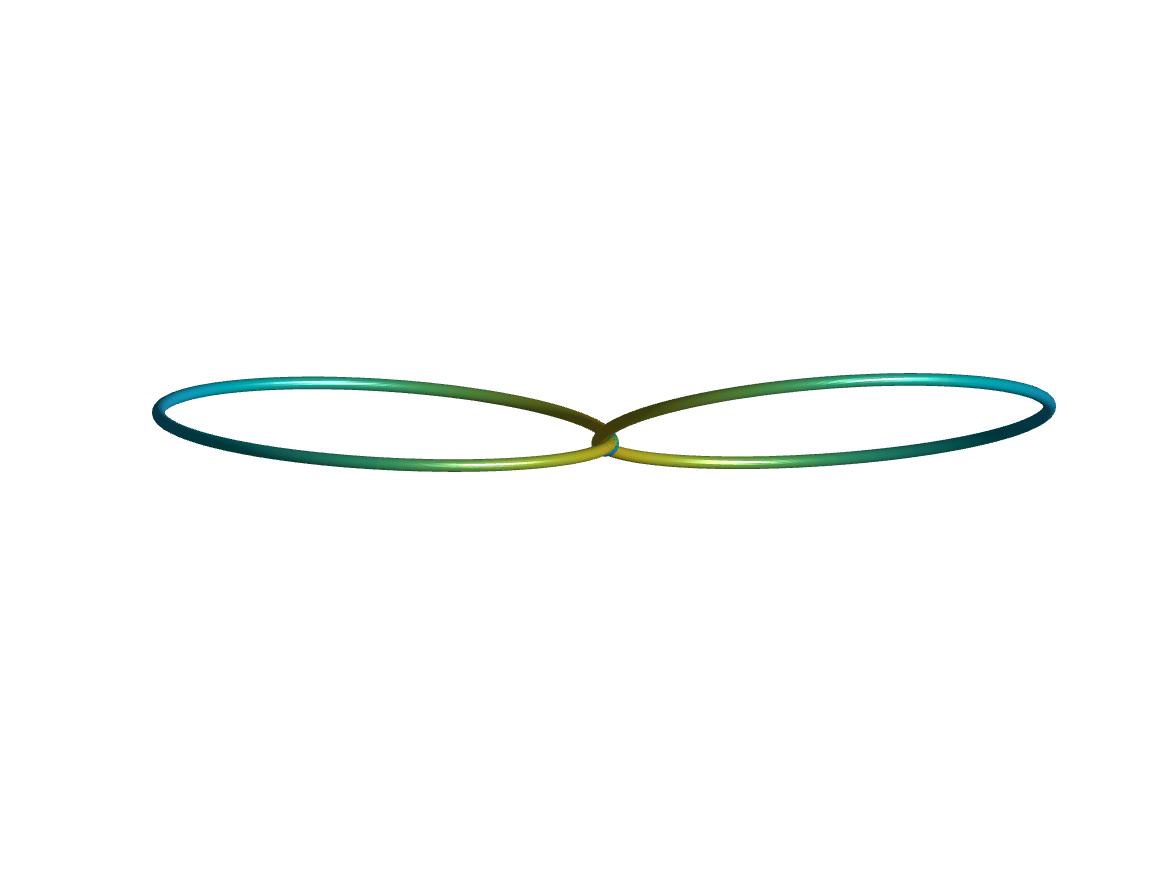}{20000}

\bild{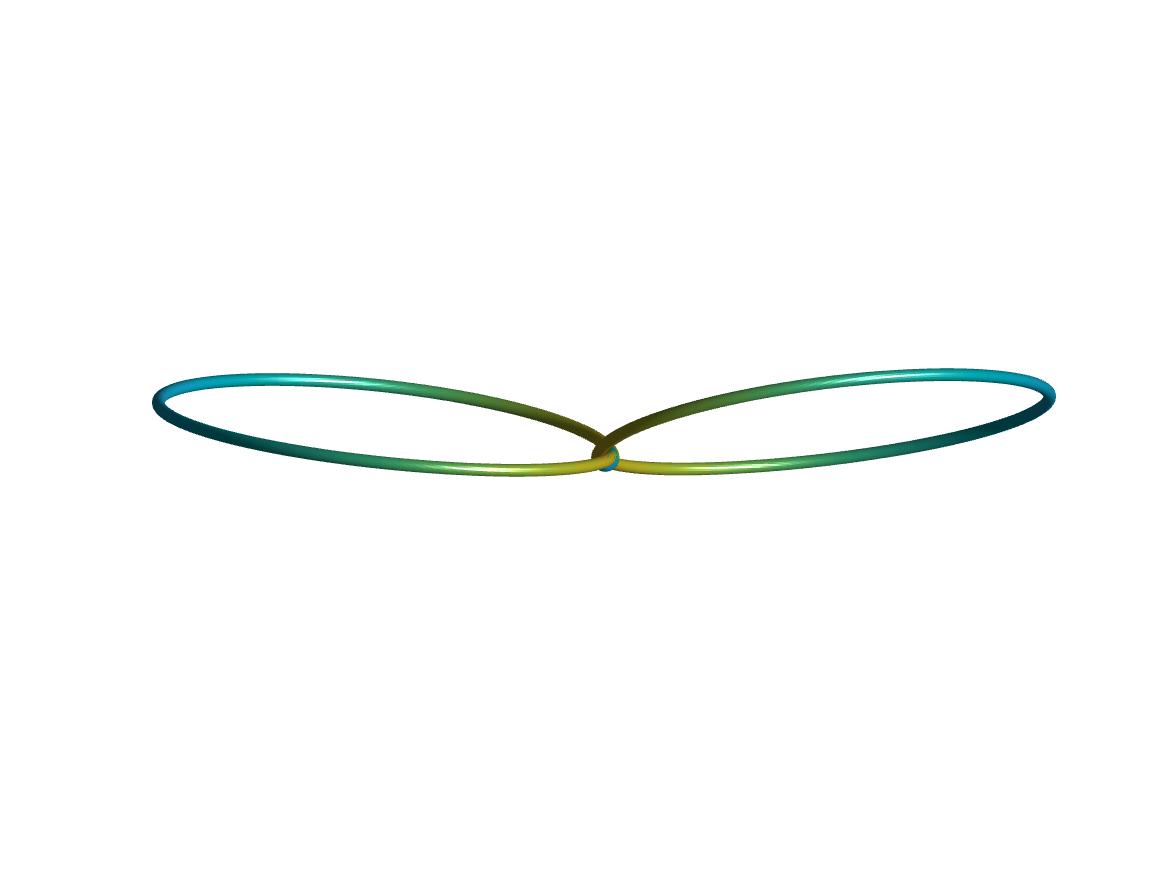}{25000}
\bild{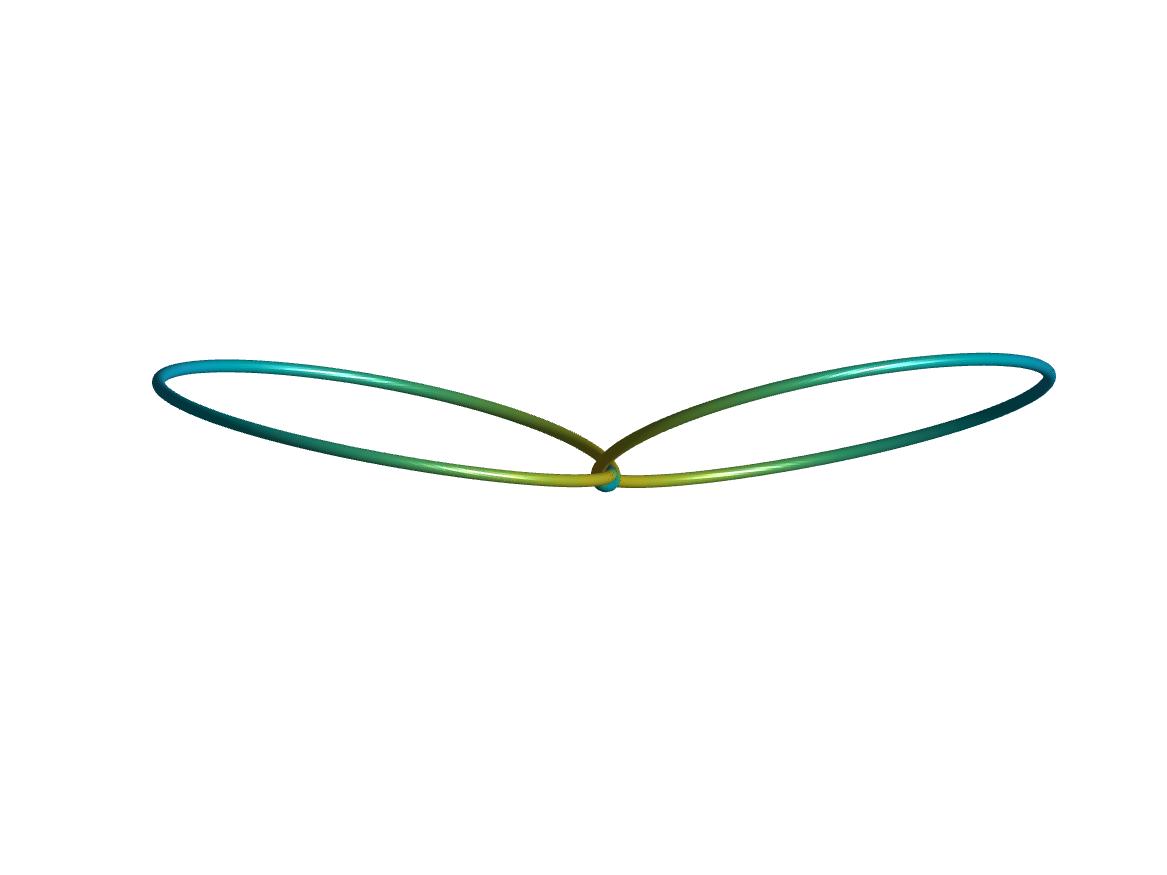}{30000}
\bild{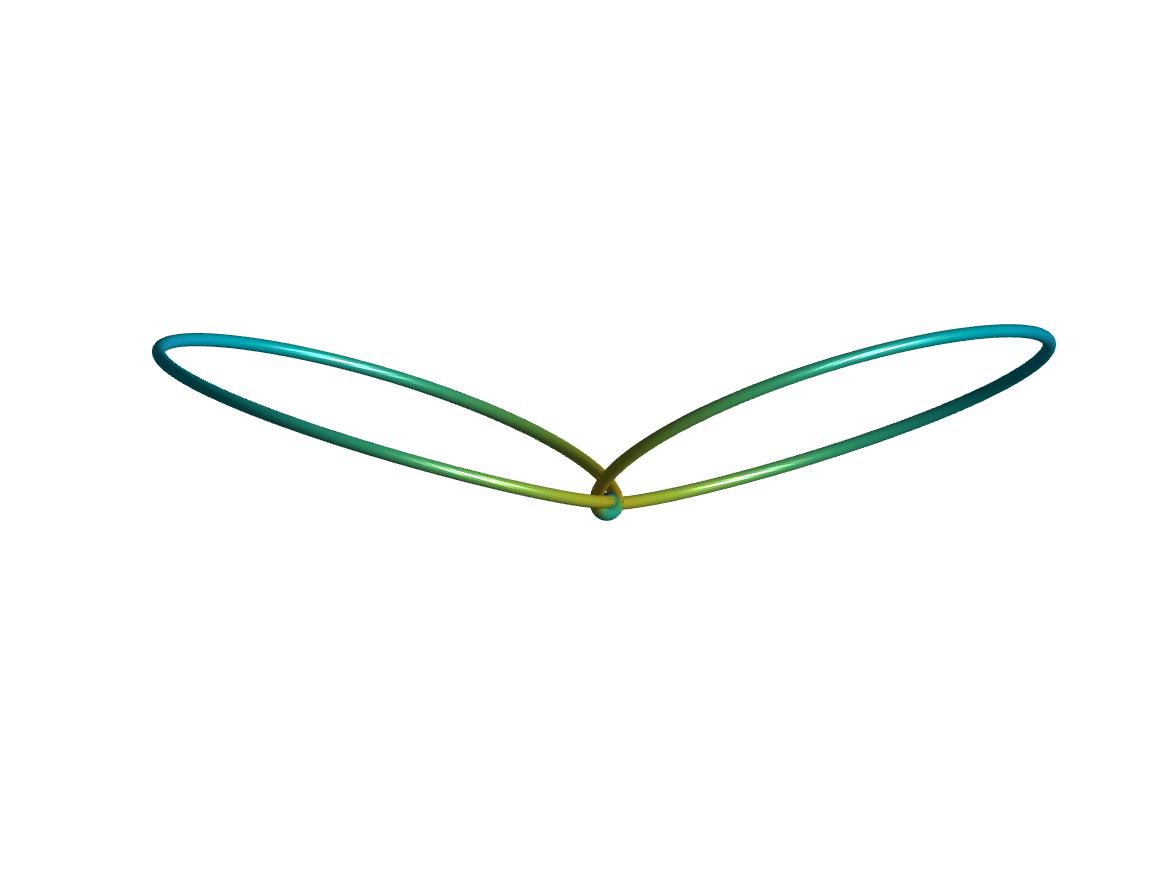}{35000}
\bild{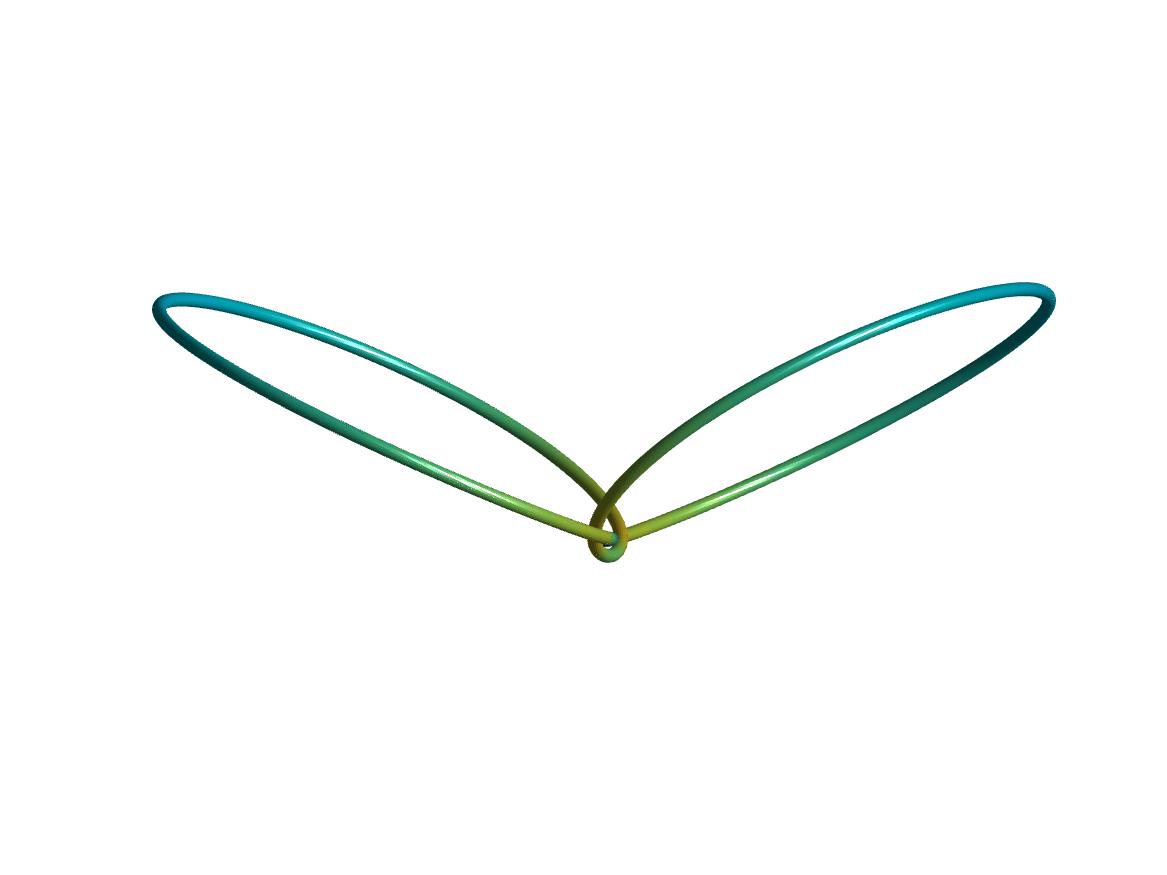}{40000}

\bild{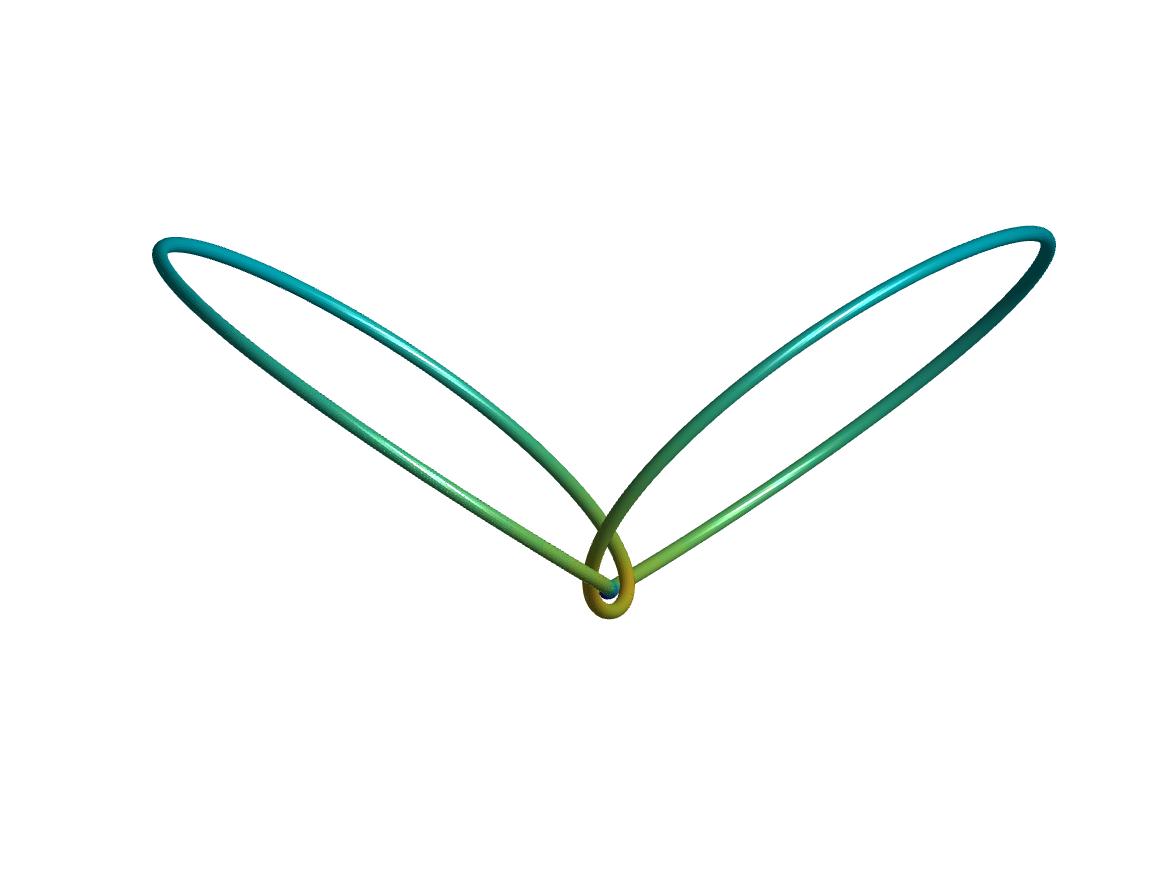}{45000}
\bild{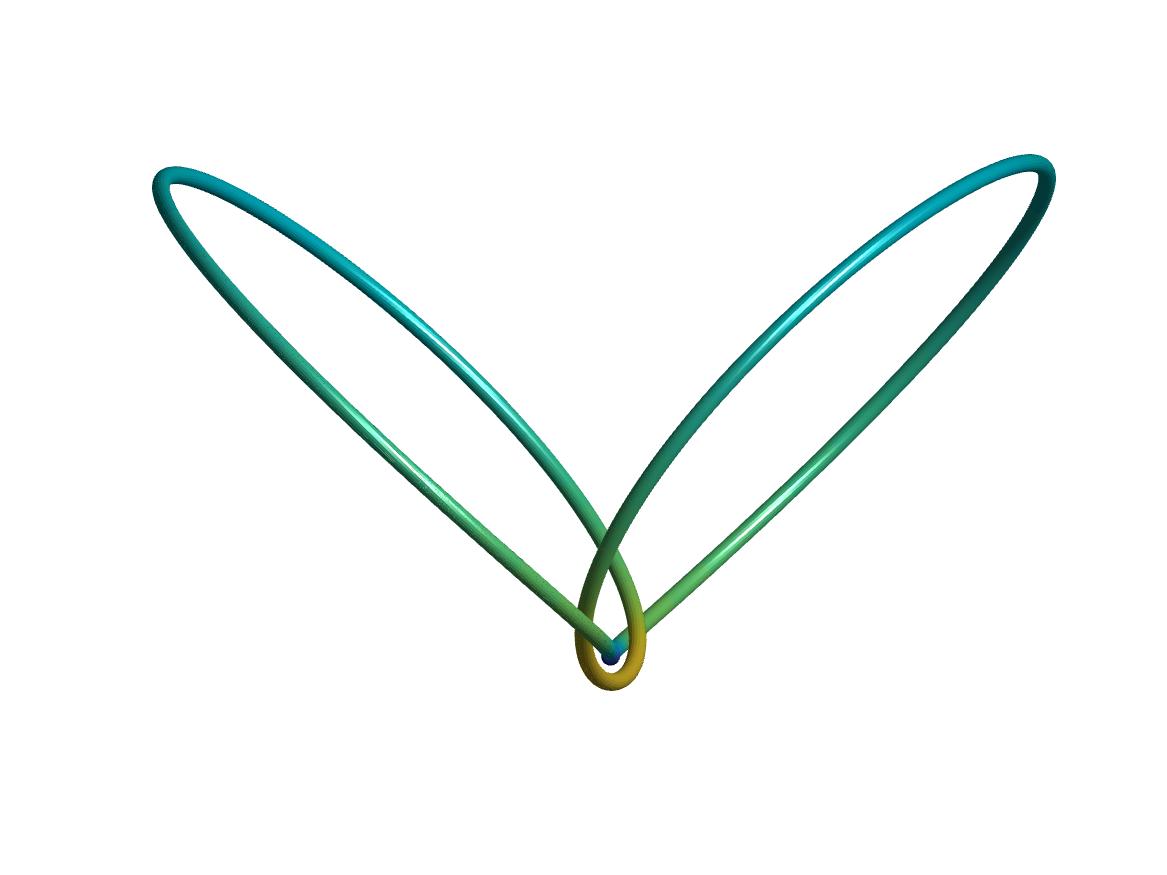}{50000}
\bild{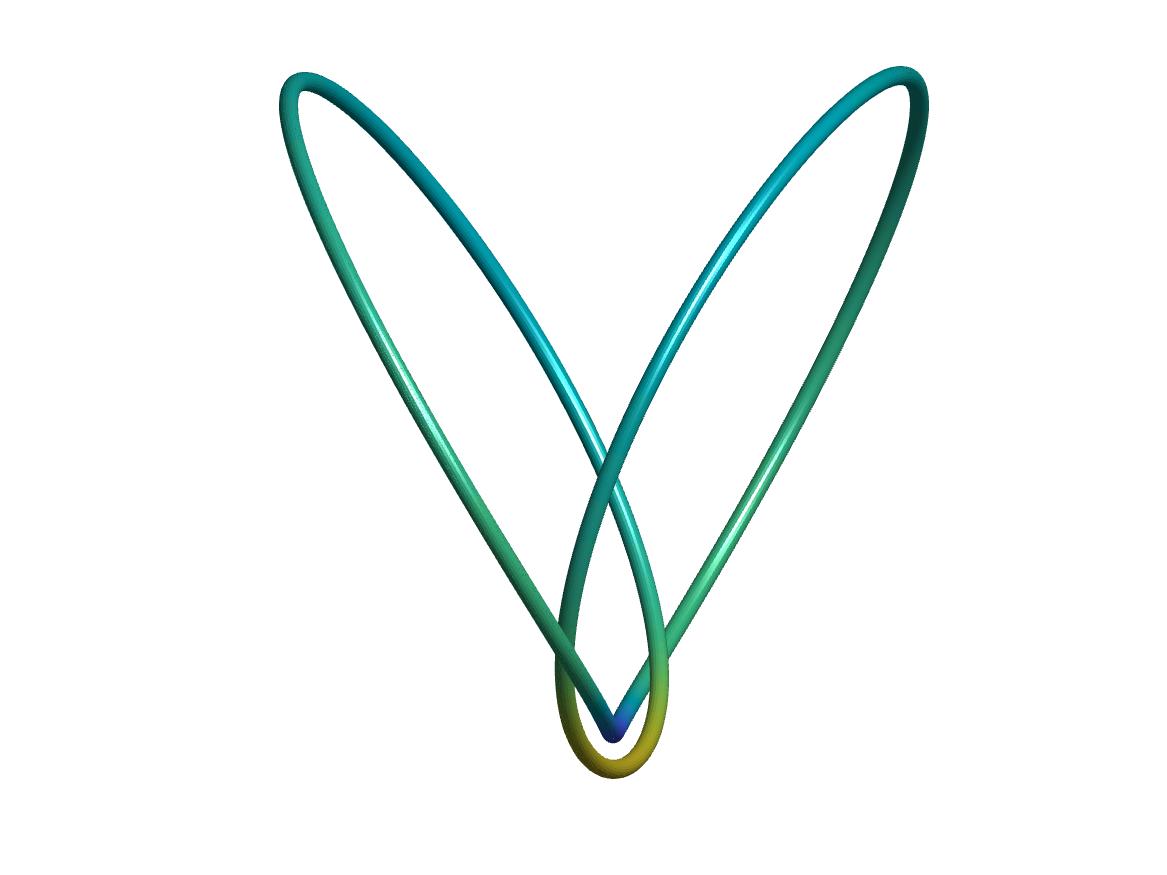}{60000}
\bild{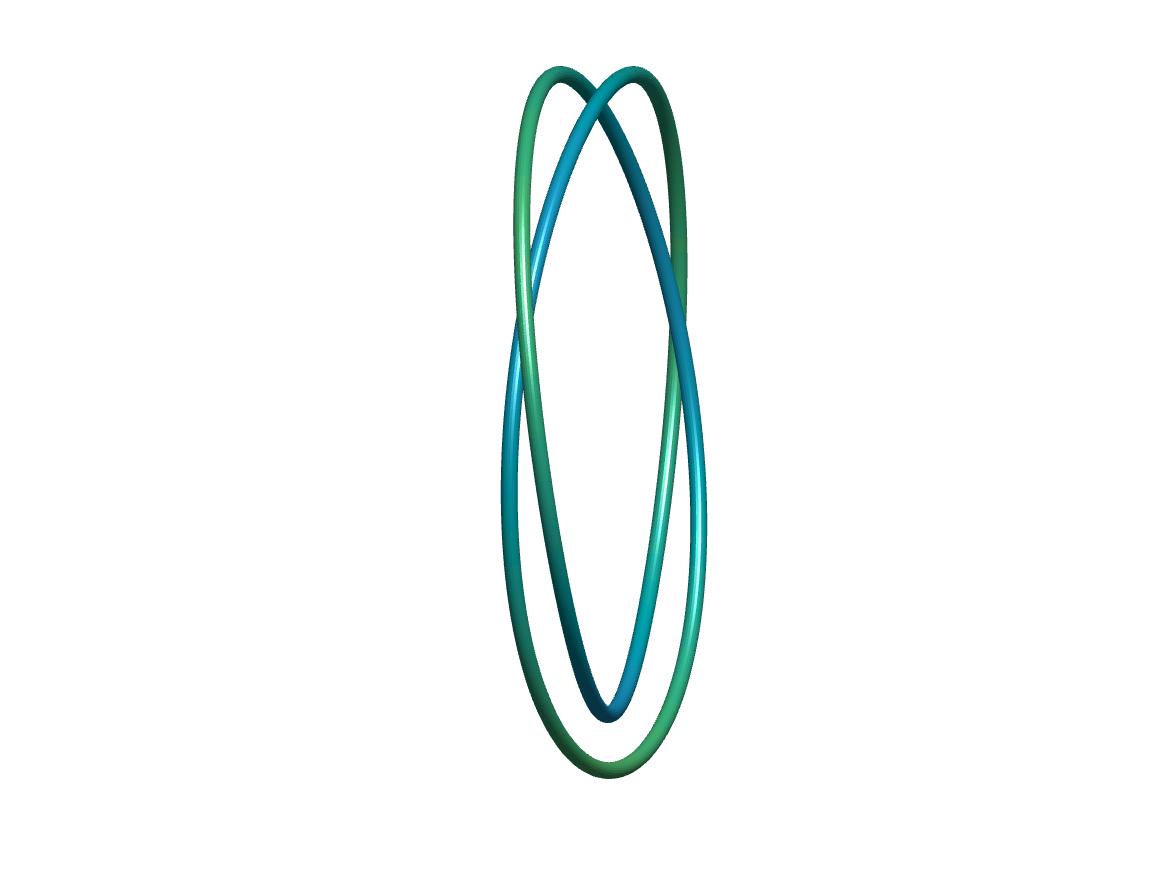}{80000}

\caption{\label{fig:tref(b)}\label{fig:tref(a)}
Snapshots of the unperturbed (top) and perturbed (bottom) discrete
evolutions for the initial trefoil configuration defined in Section~\ref{sec:tref}.
The unpertubed evolution becomes (experimentally) stationary in a symmetric
but not energy minimizing configuration. Small perturbations lead to 
a different evolution \sbbb{after approximately 25,000 steps}
with a stationary configuration approximating a doubly
covered circle that is globally energy minimizing.
}

\end{figure}

\begin{figure}
\includegraphics[scale=.475]{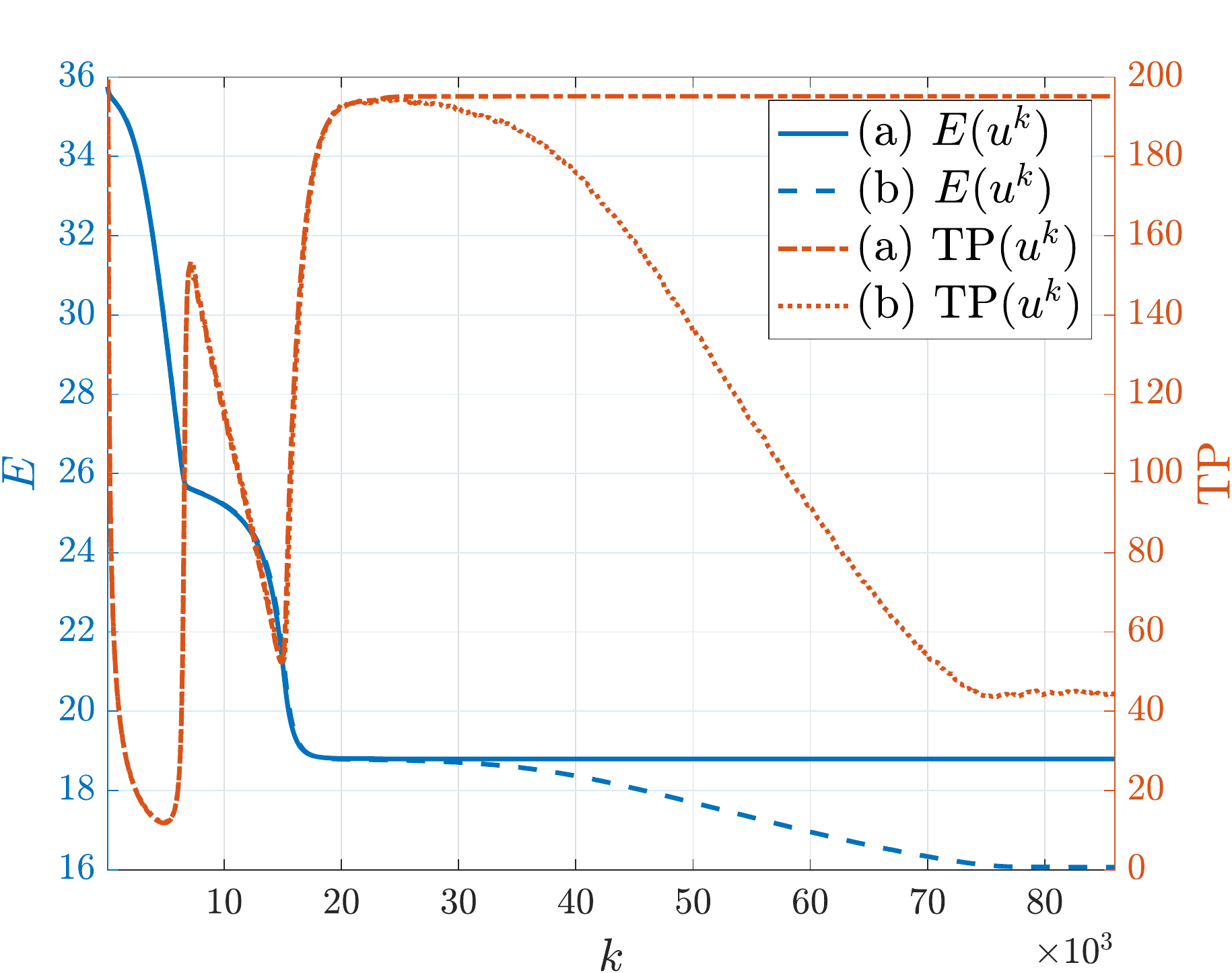}
\caption{Energy curves for the \sbbb{(a)~unperturbed and (b)~perturbed}
discrete evolutions for the initial configuration defined in  
Section~\ref{sec:tref}. The curves are nearly identical
up to step 25,000 when the unperturbed evolution becomes stationary
as depicted in Figure~\ref{fig:tref(a)}.}\label{fig:tref}
\end{figure}

\subsection{Figure-eight}\label{sec:f8}

Our third example illustrates how different initial configurations within
one knot class lead to different stationary configurations and thereby
show the limitations of gradient flows to determine particular representatives
of a given class. 


So far, there are no analytical results concerning the elastic figure-eight.
Numerical experiments carried out by Gallotti and Pierre-Louis~\cite{gallotti}
as well as by Gerlach et al.~\cite{GRvdM} led to a spherical configuration exhibiting 
a remarkable symmetry (as the final state shown in the top part of Figure~\ref{fig:f8(b)}).
Avvakumov and Sossinsky~\cite{sossinsky}
instead claim that a planar configuration (as in the bottom part of Figure~\ref{fig:f8(b)})
yields a lower energy value compared to the spherical curve
which they consider being merely a local minimizer.
Our experiments indicate some support for the latter observation.

We retrieved coordinates of knotted curves
\sbbb{from the website~\cite{knotserver},} namely
(a)~\href{http://www.colab.sfu.ca/KnotPlot/KnotServer/mseq-coord/3.html}{\tt{}mseq-coord/3.html} and
(b)~\href{http://www.colab.sfu.ca/KnotPlot/KnotServer/coord/3.html}{\tt{}coord/3.html}.
In order to produce suitable initial curves,
we added further nodes by cubic interpolation,
performed a few iteration steps including some small randomized perturbation
in order to allow for smoothing, and
then rescaled the curve to a length of $50$ units.

The resulting configuration was taken as the initial curve for the respective experiment
which again involved performing small randomized perturbations
each hundredth step in order to break symmetry.
The discrete evolutions are shown in Figures~\ref{fig:f8(b)}.
The corresponding energy plots can be found in Figure~\ref{fig:f8}.


\renewcommand{\bild}[2]{\fbox{\includegraphics[scale=.12,trim=290 210 230 160,clip]{#1}}\makebox[0ex][r]{\tiny#2\ }\,\ignorespaces}

\begin{figure}

\bild{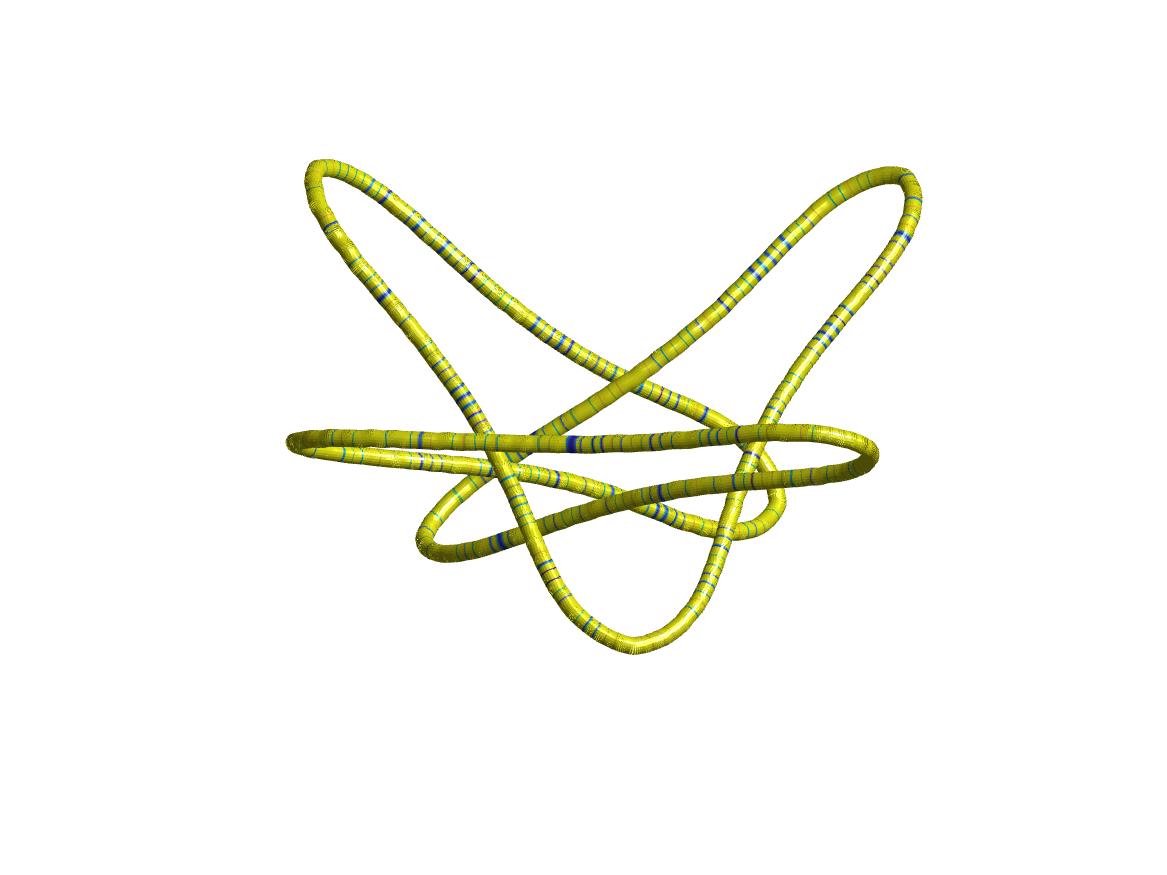}{1}
\bild{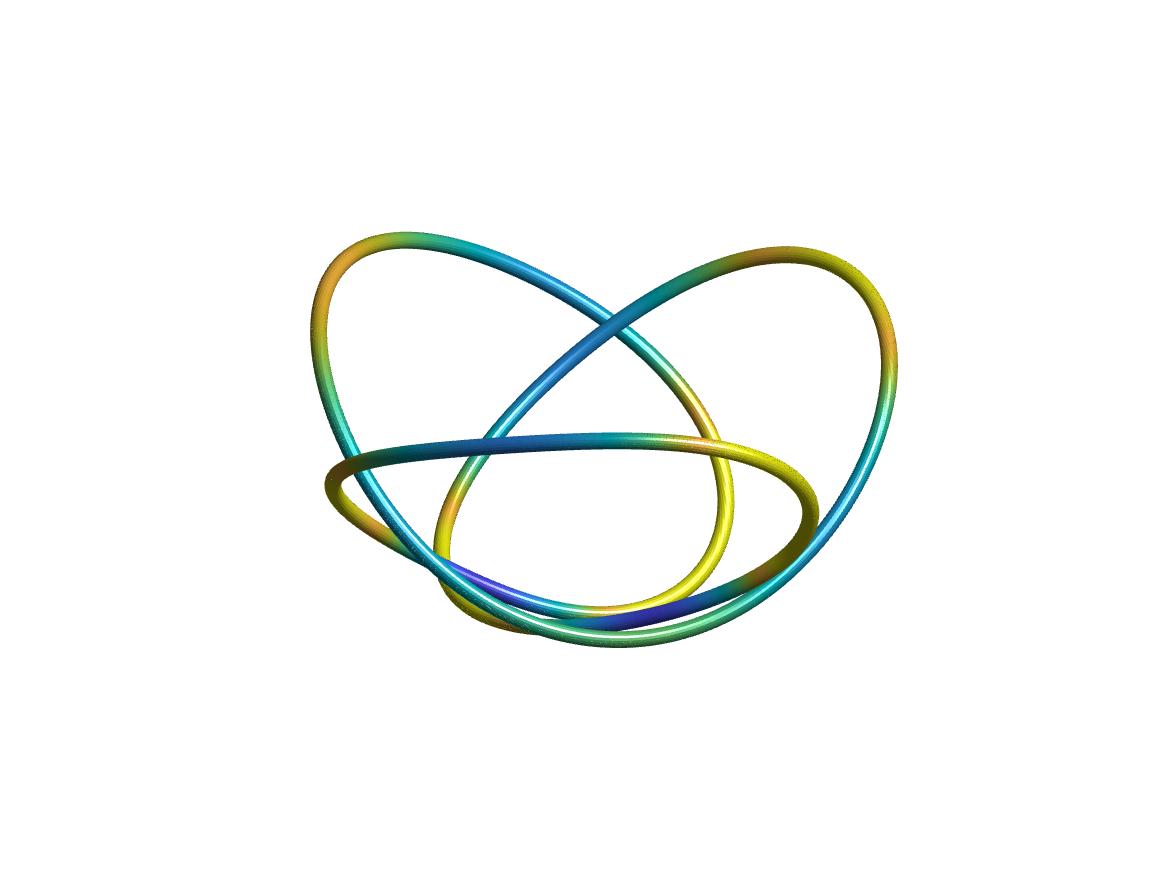}{100}
\bild{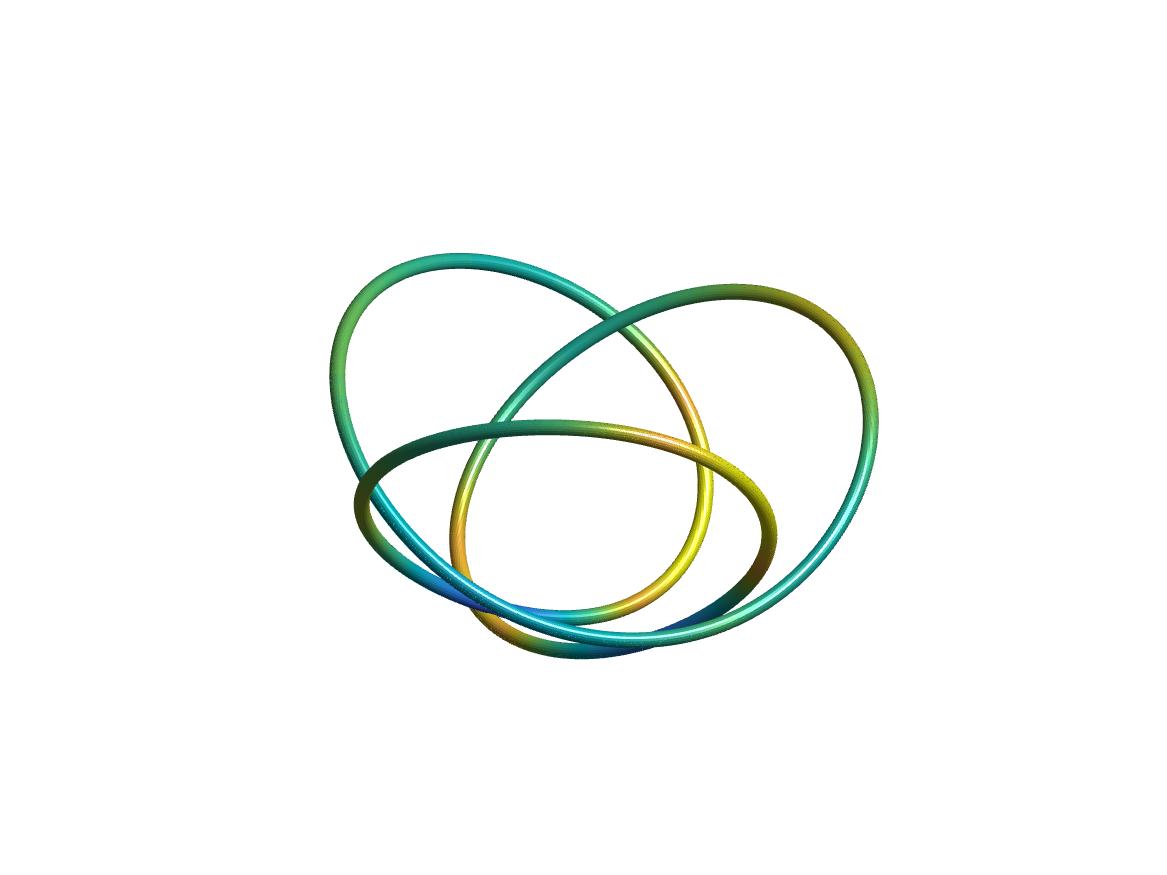}{200}
\bild{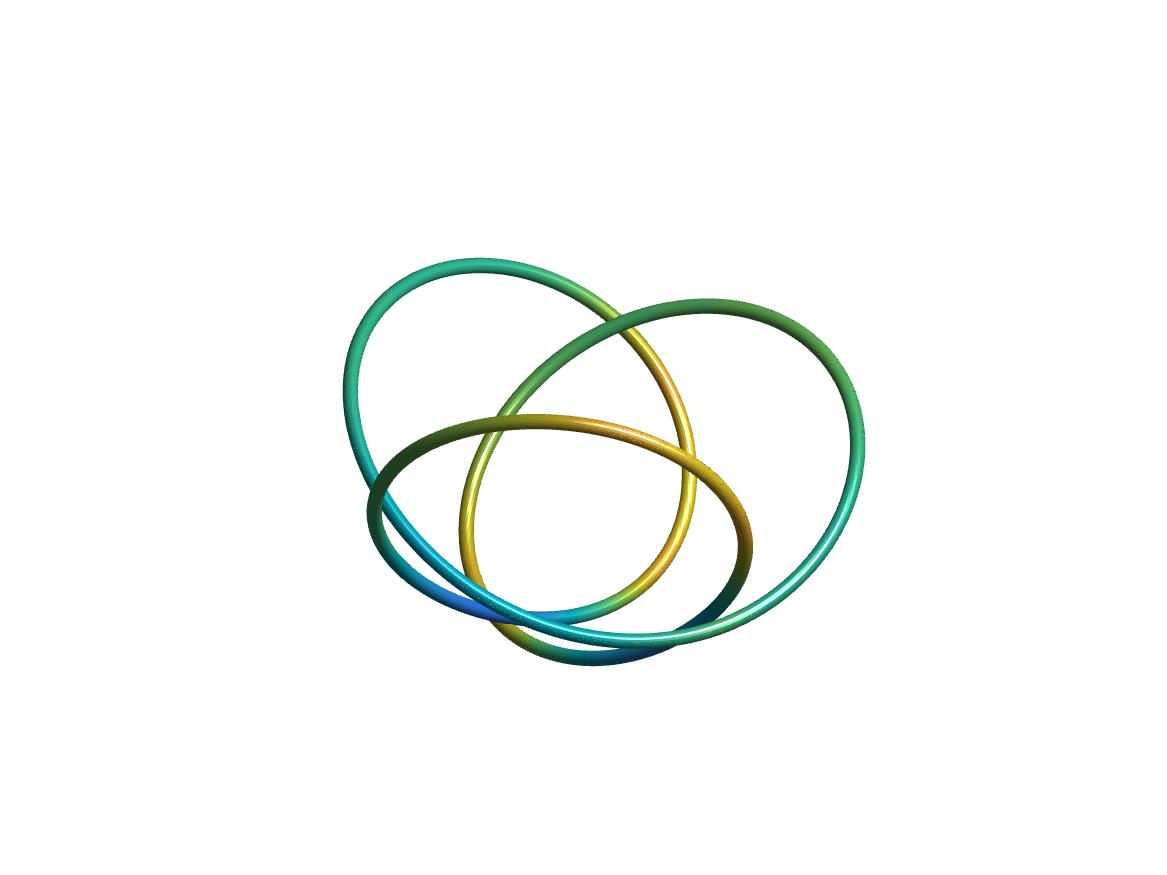}{300}

\bild{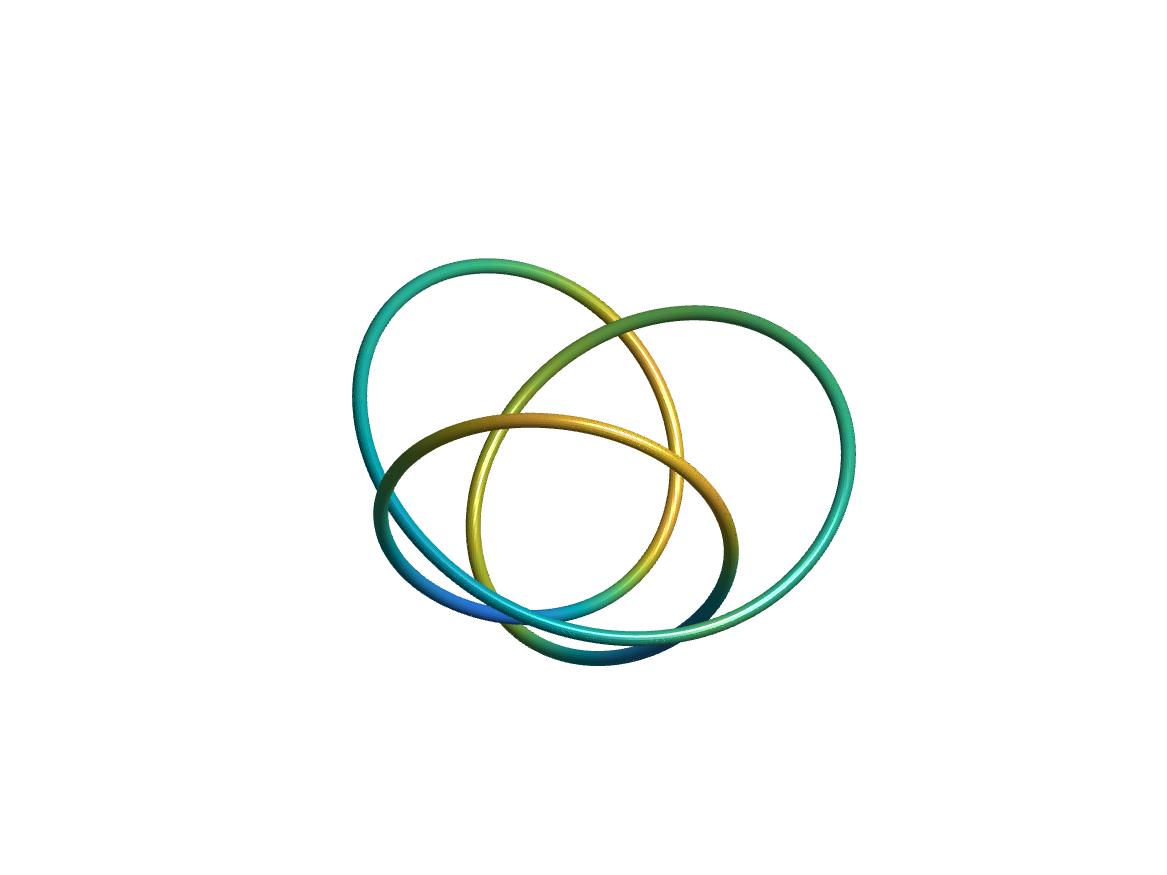}{400}
\bild{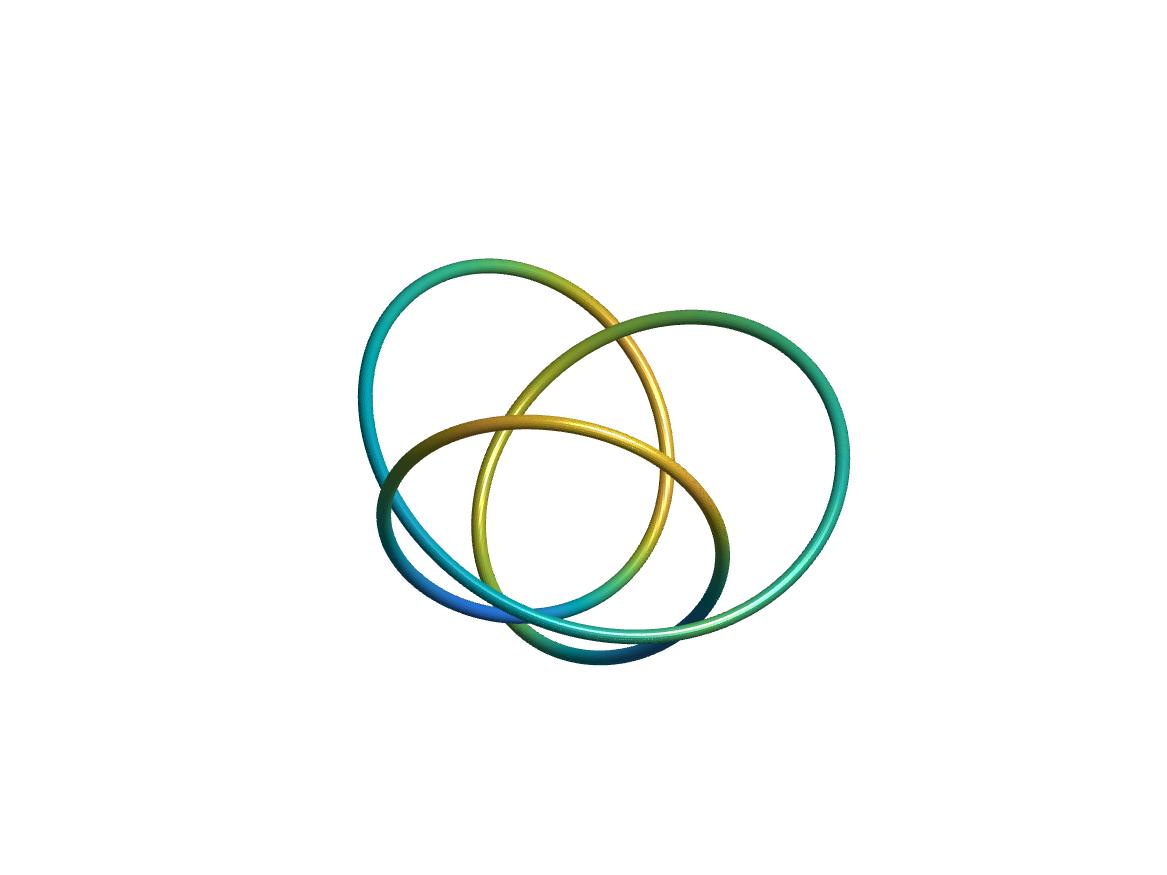}{500}
\bild{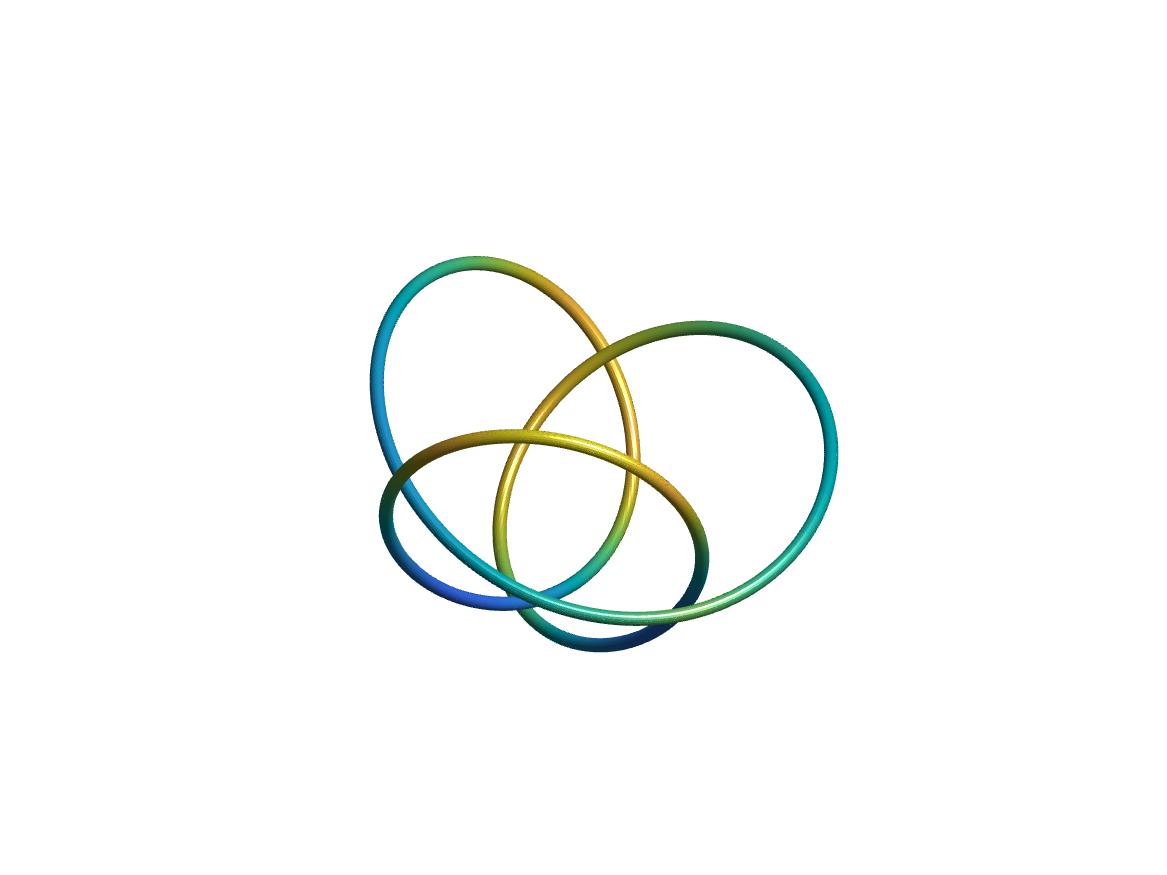}{1000}
\bild{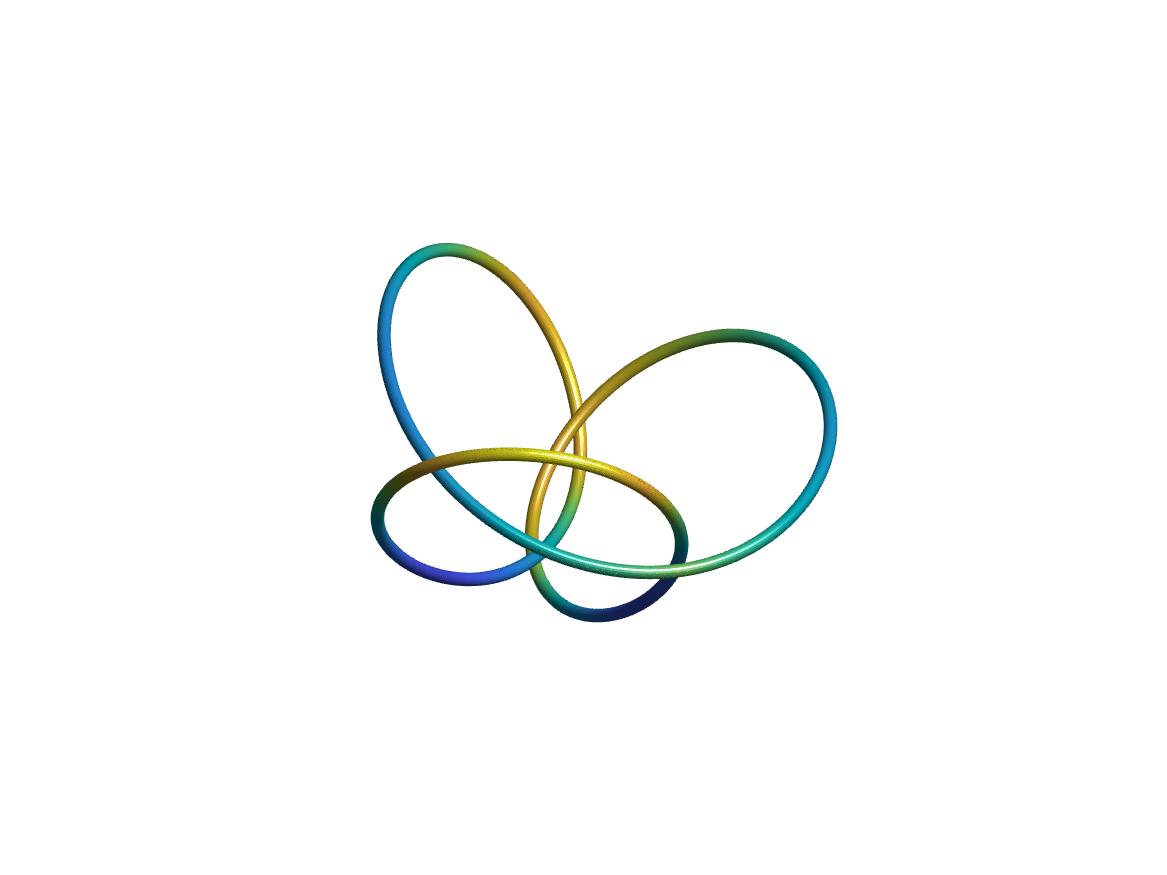}{2000}

\bild{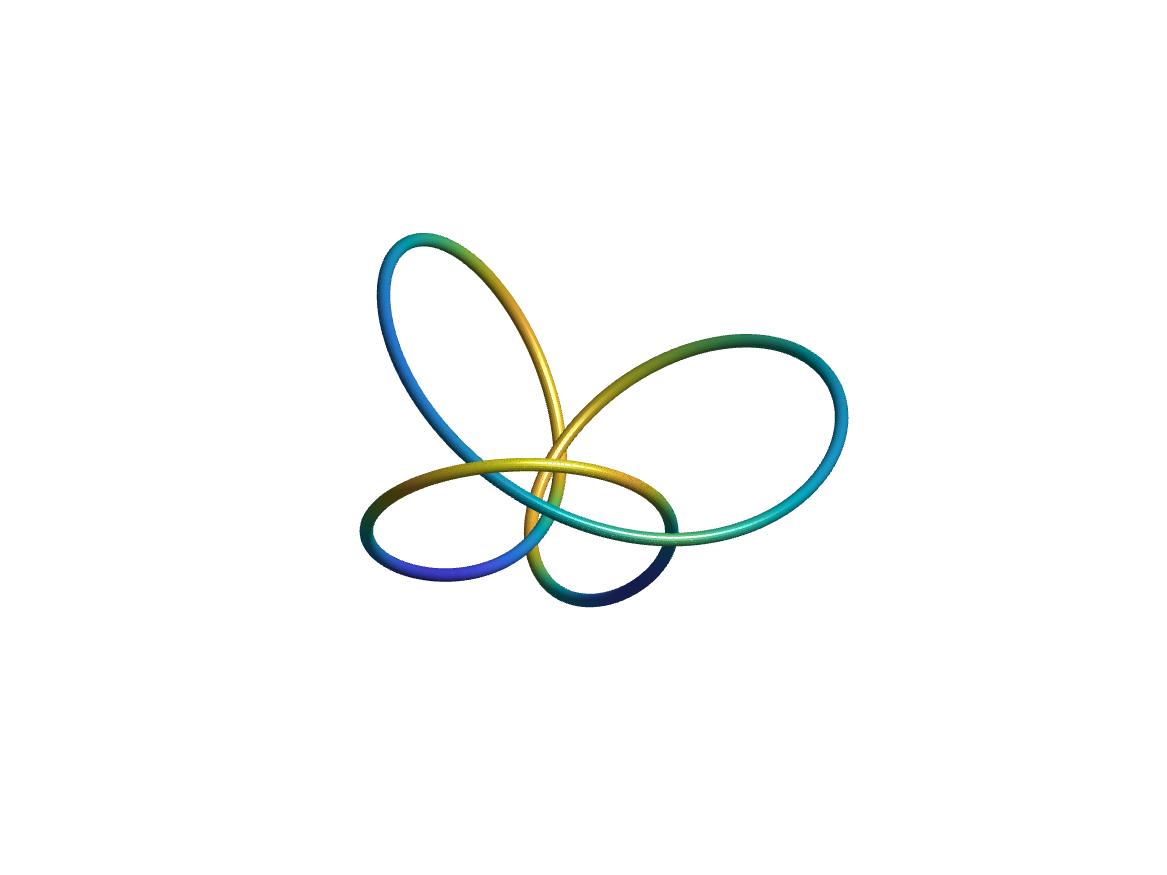}{3000}
\bild{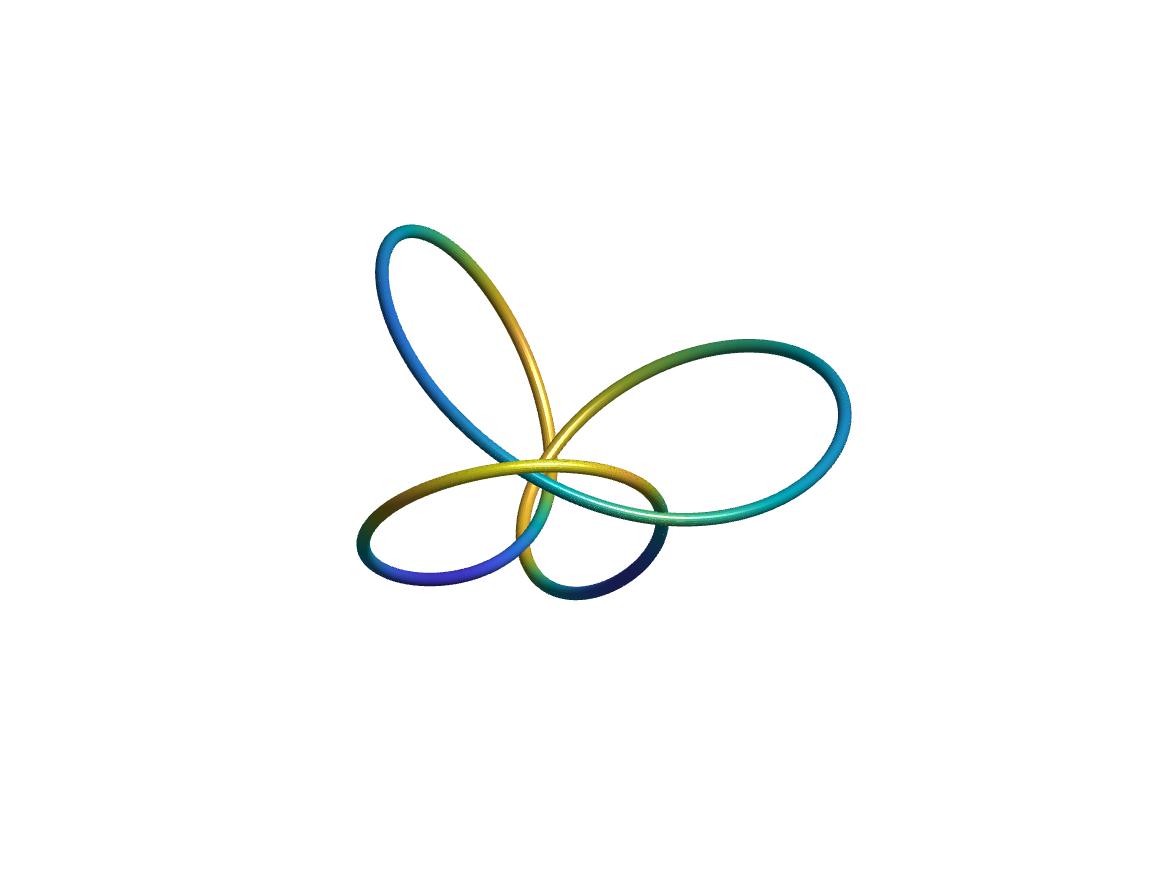}{4000}
\bild{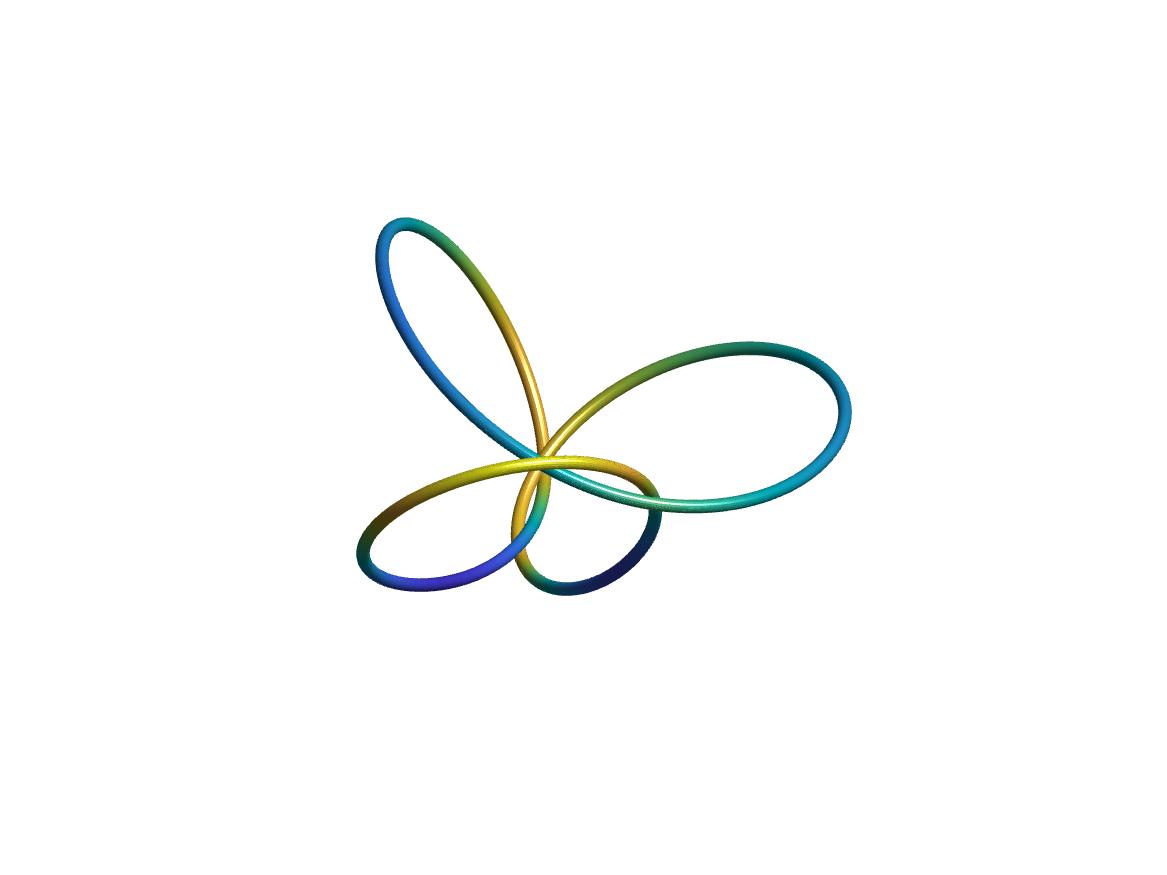}{5000}
\bild{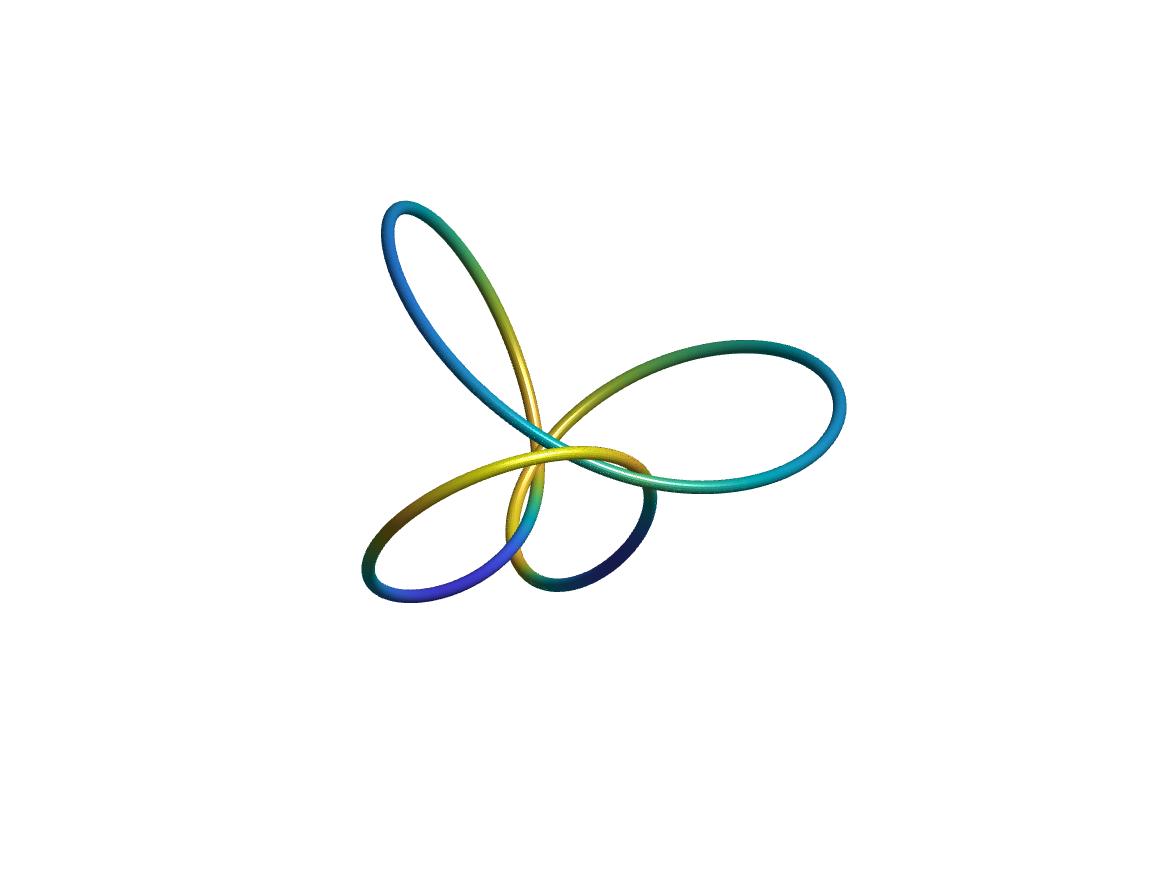}{10000}

\bild{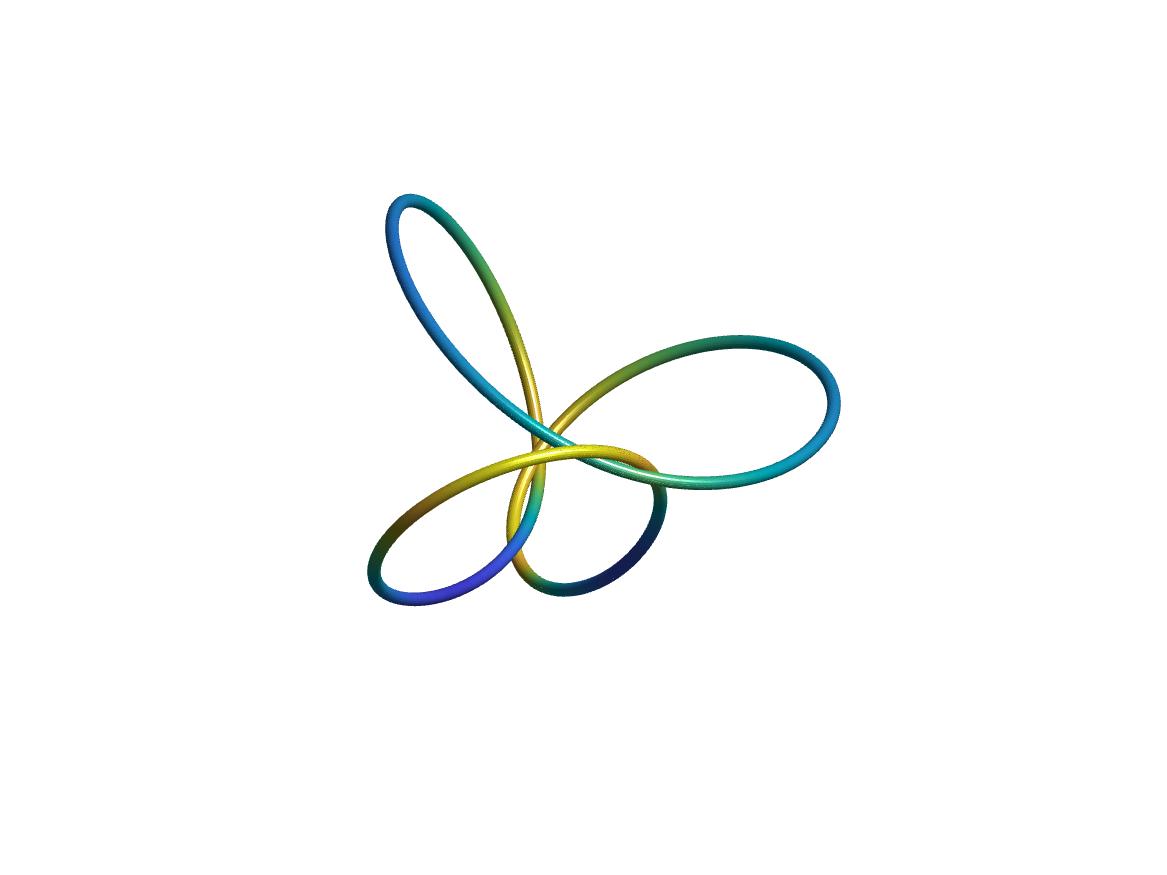}{15000}
\bild{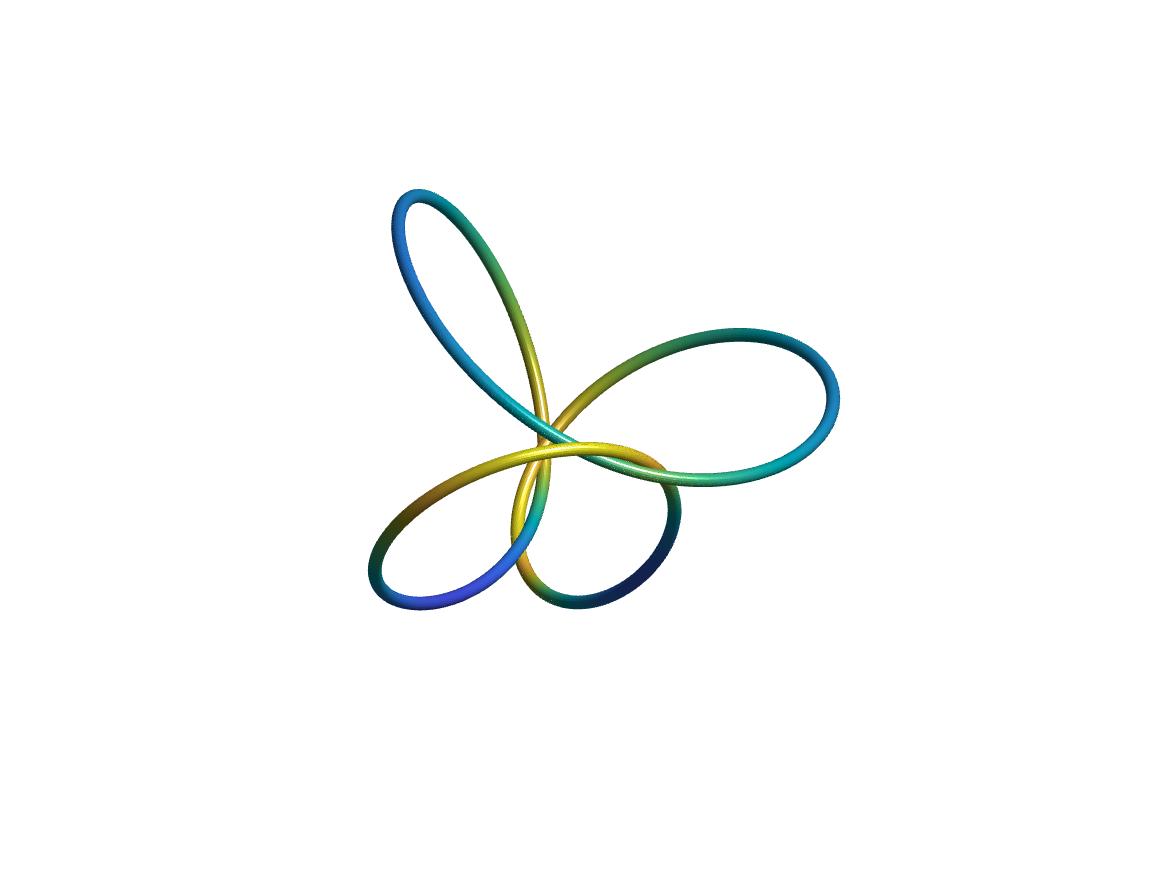}{20000}
\bild{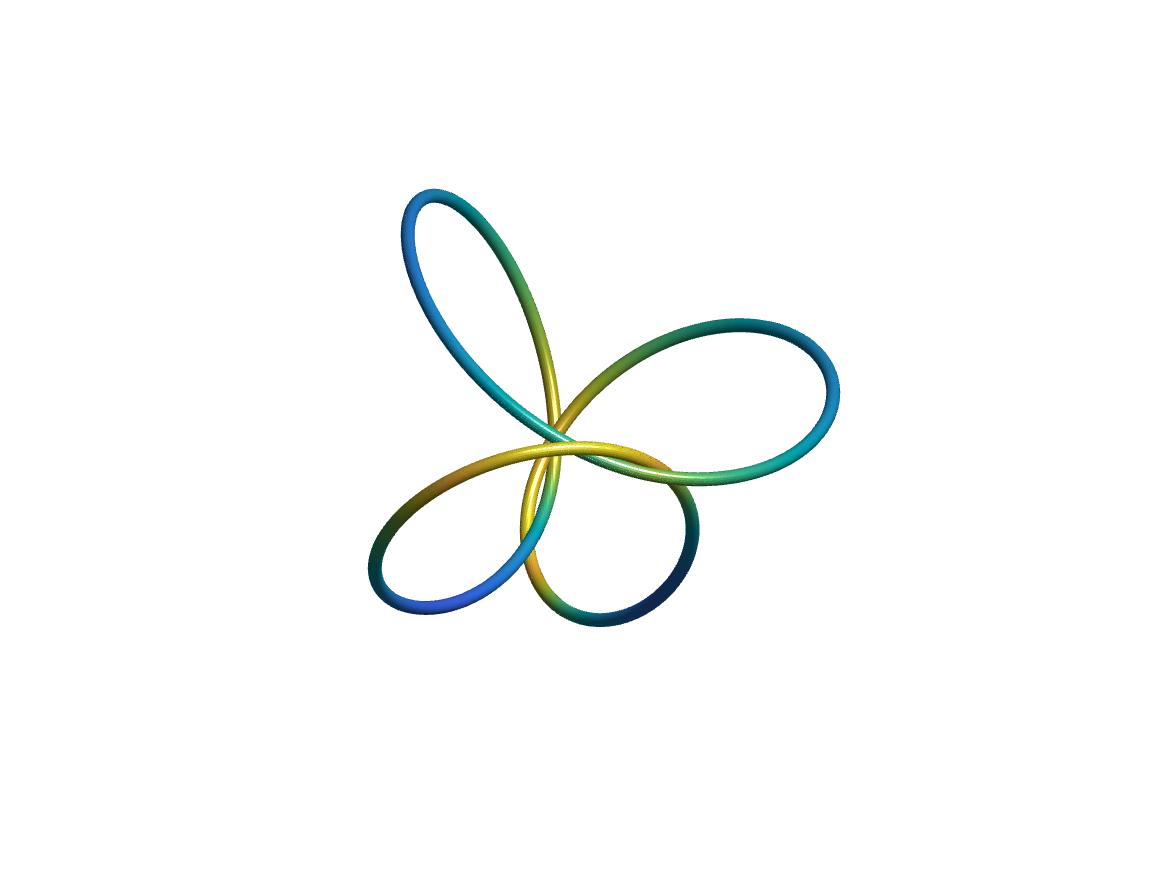}{25000}
\bild{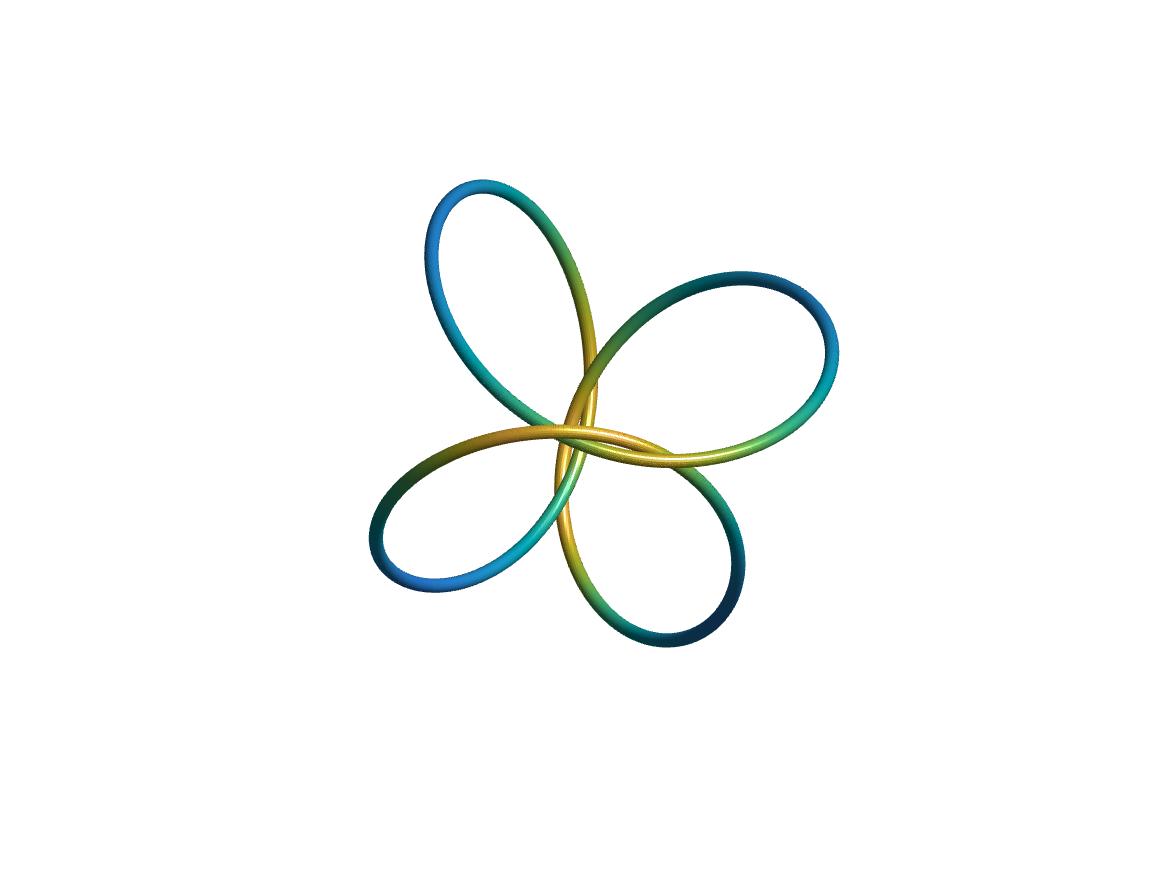}{40000}


\renewcommand{\bild}[2]{\fbox{\includegraphics[scale=.1185,trim=150 140 120 110,clip]{#1}}\makebox[0ex][r]{\tiny#2\ifnum#2>500\hspace{30ex}\fi\ }\,\ignorespaces}

\bigskip


\bild{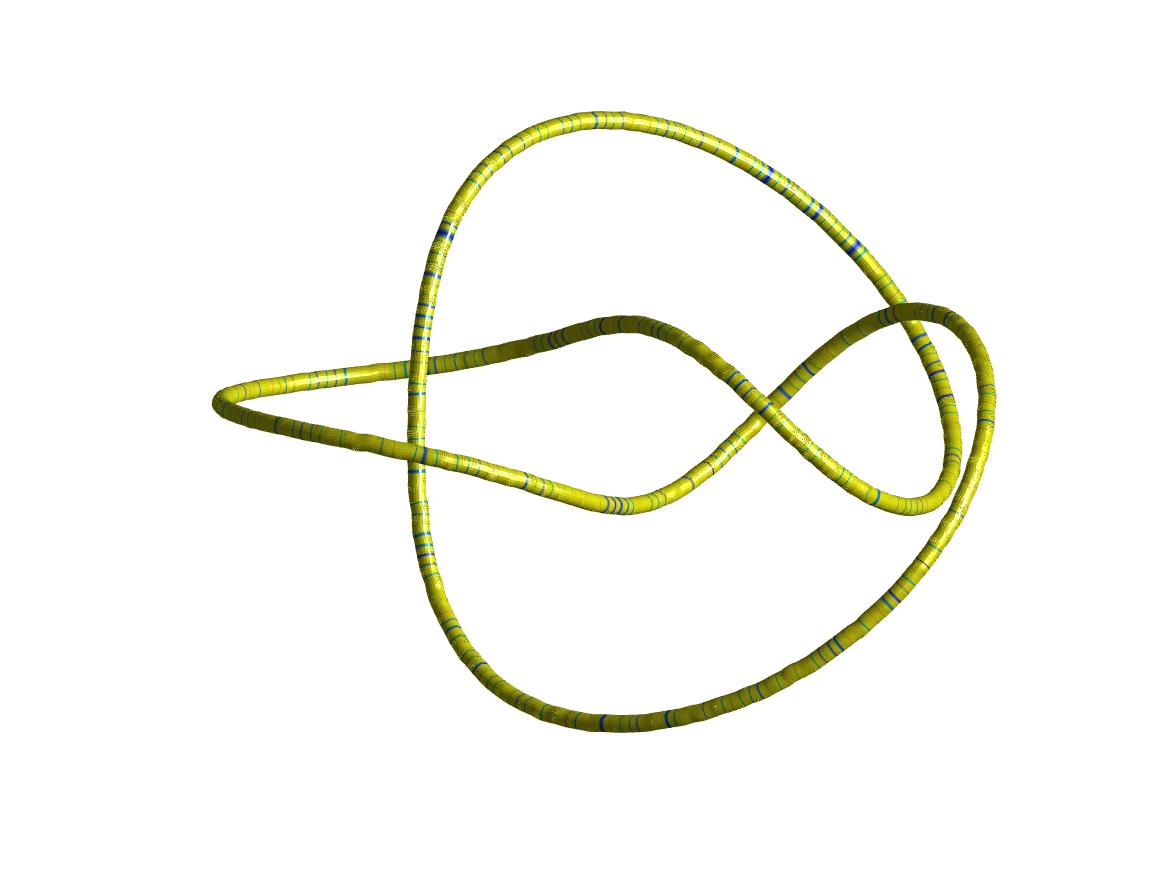}{1}
\bild{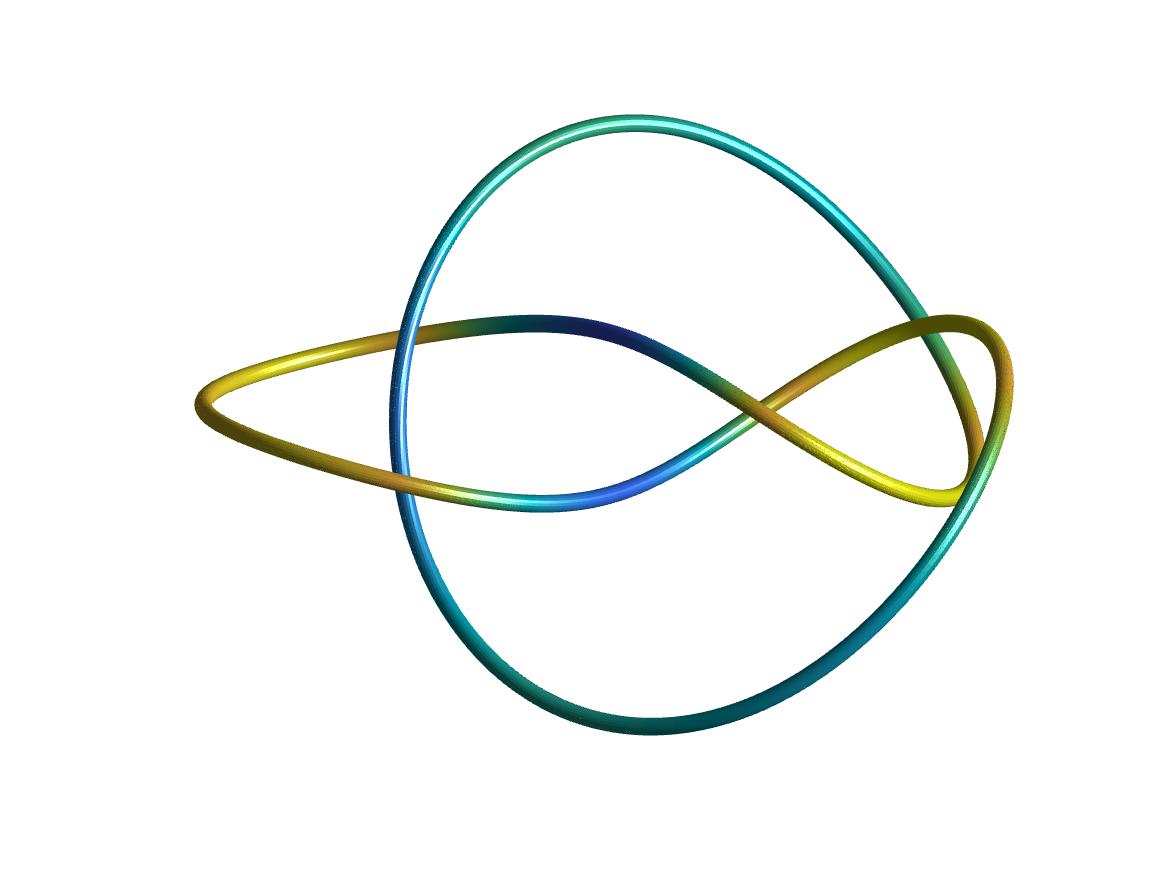}{100}
\bild{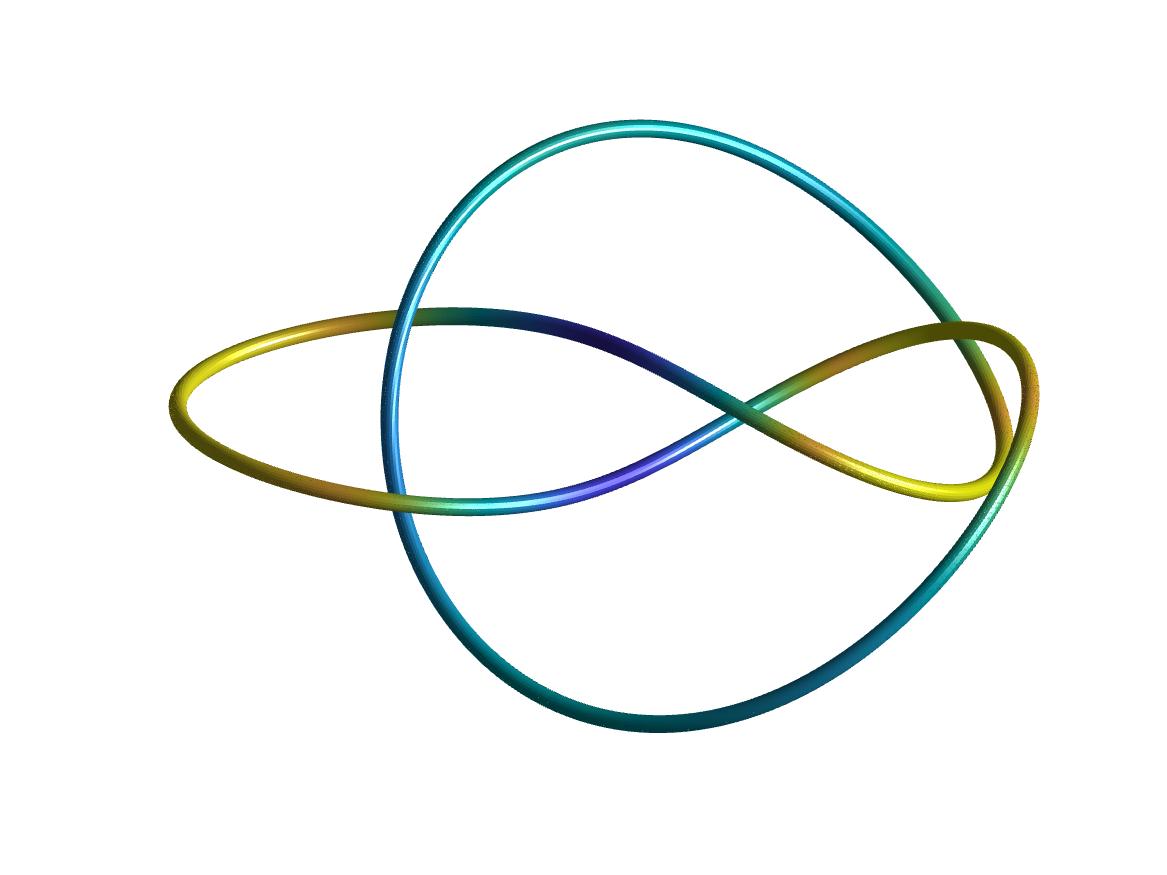}{200}

\bild{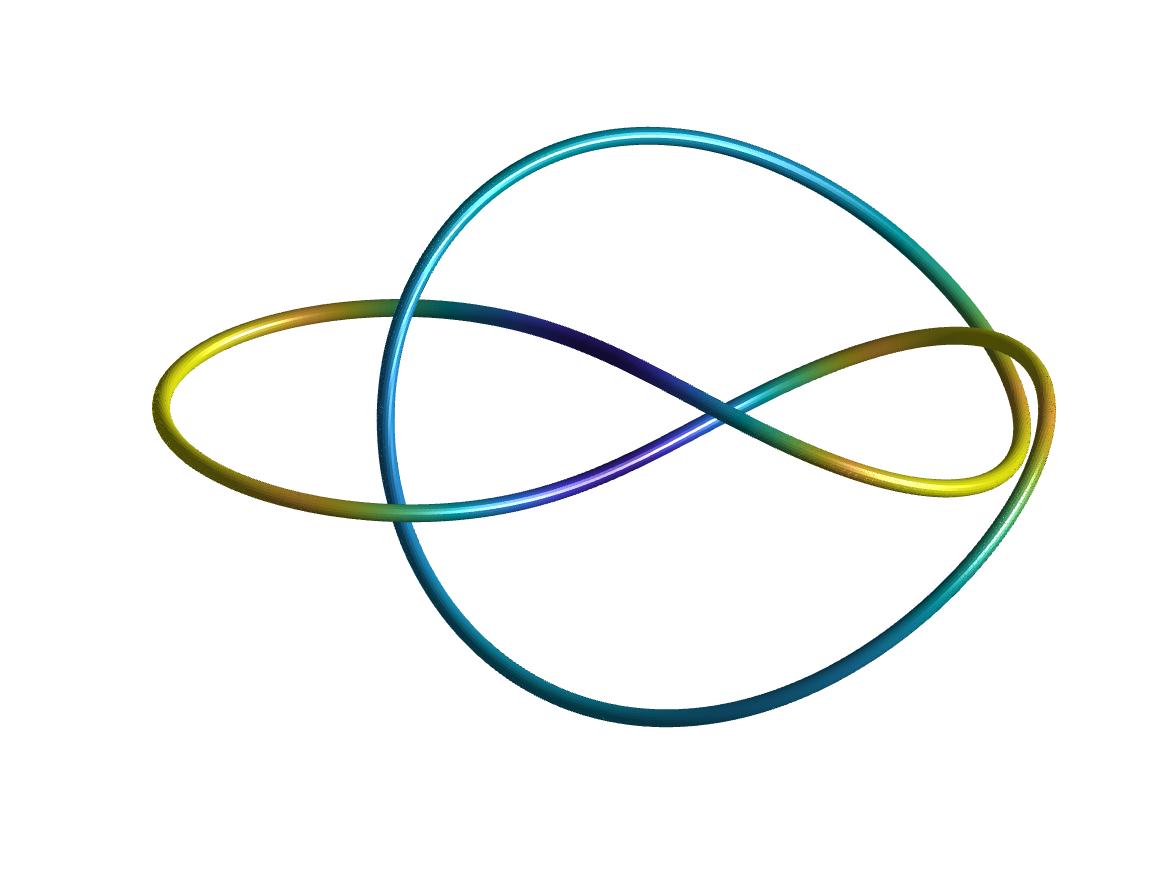}{300}
\bild{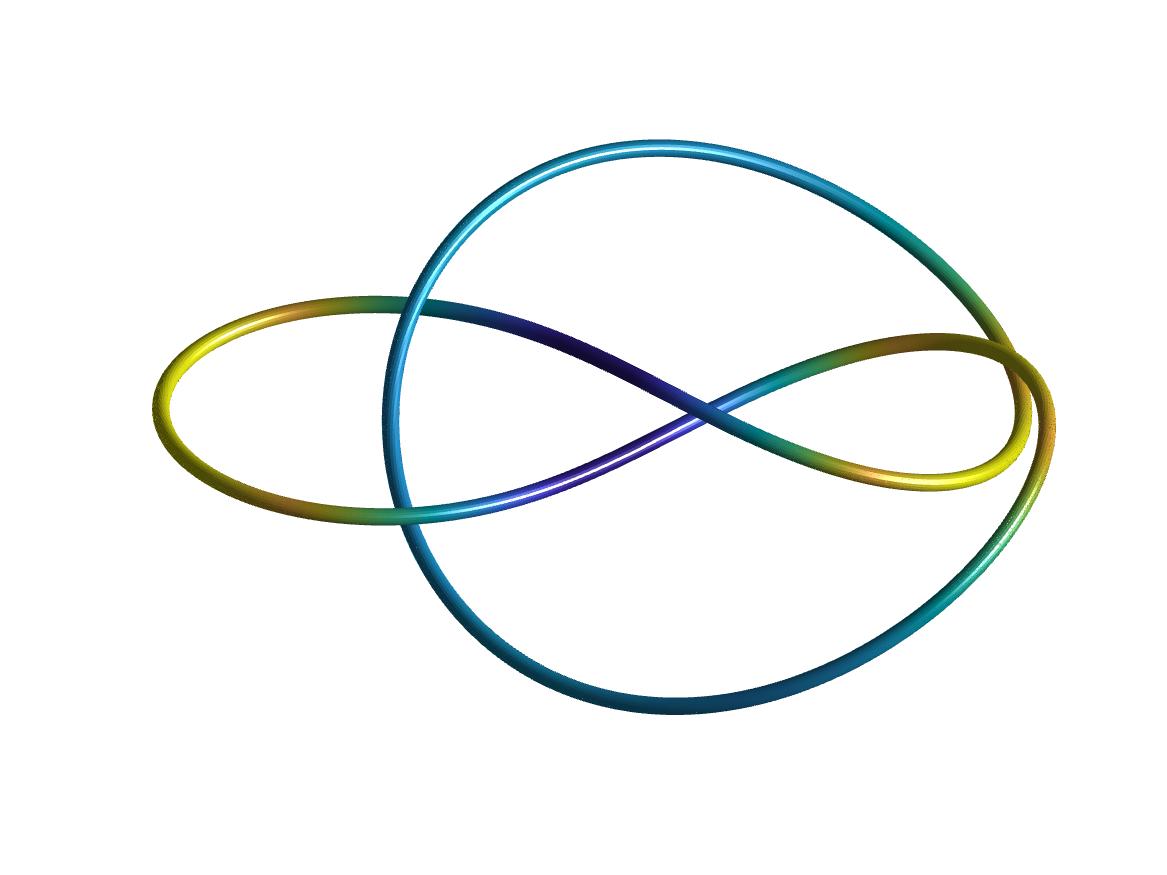}{400}
\bild{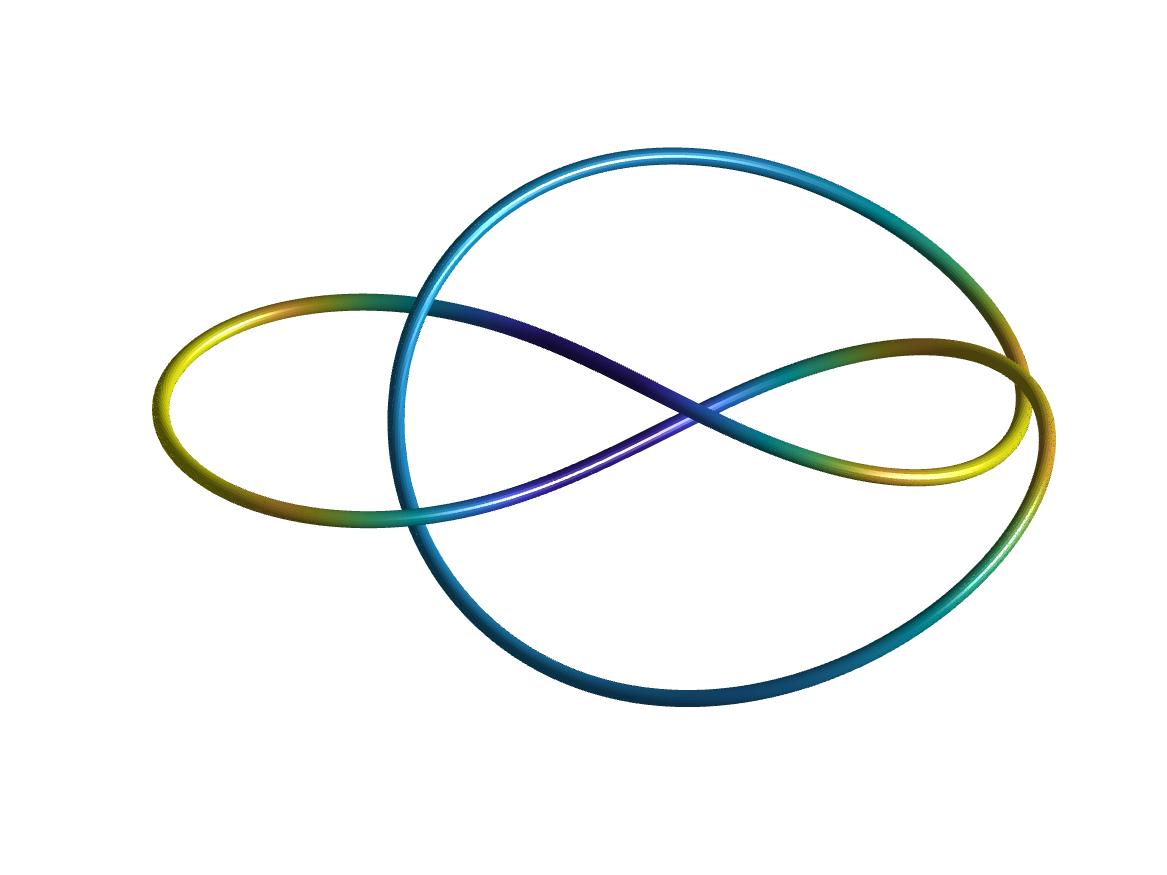}{500}

\bild{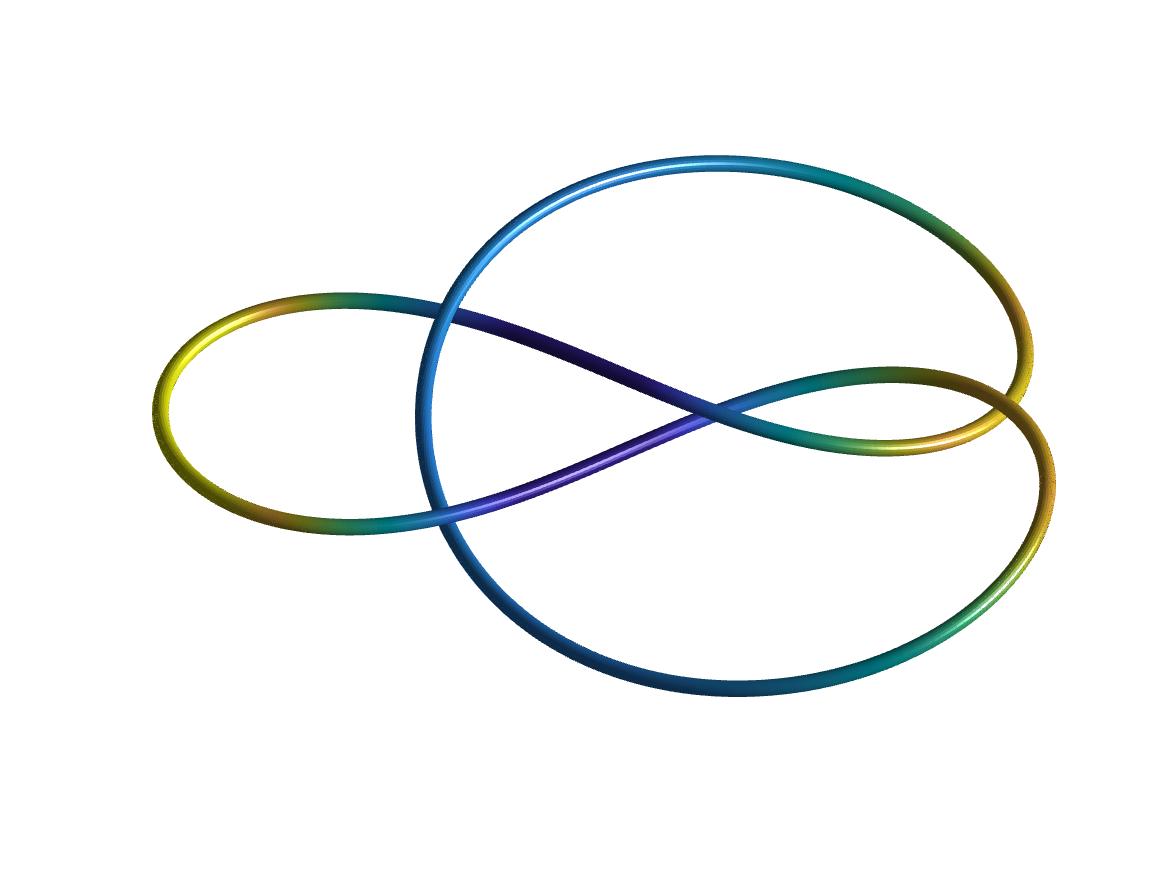}{1000}
\bild{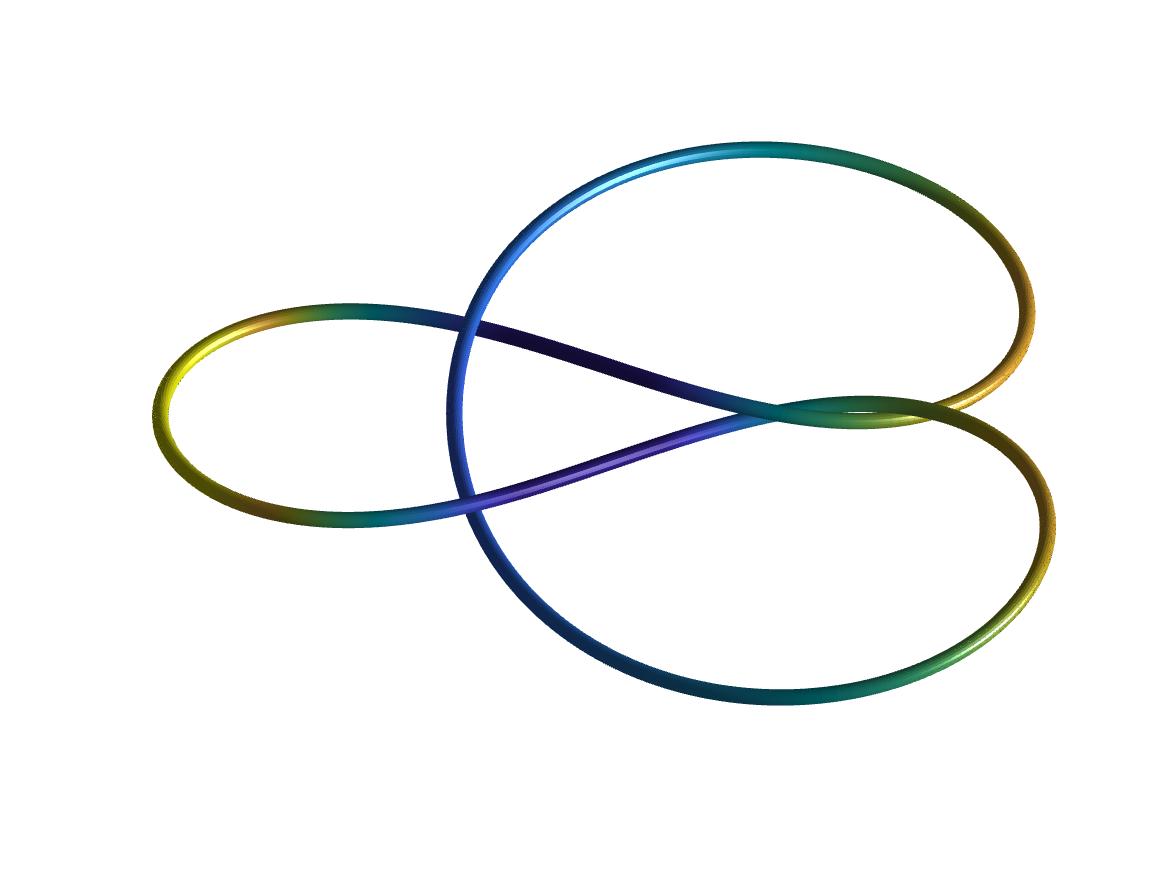}{2000}
\bild{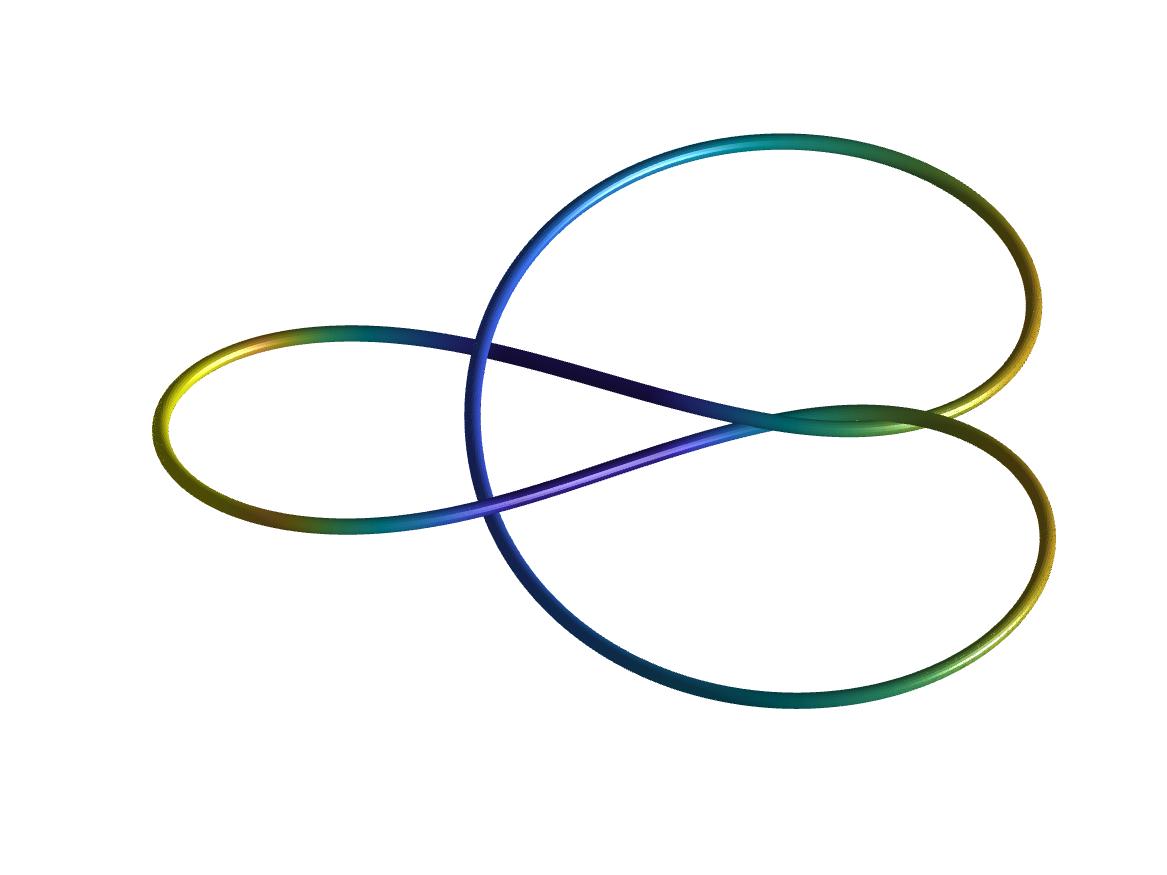}{40000}

\caption{\label{fig:f8(b)}
Snapshots of two different discrete evolutions defined by initial
configurations~(a) and~(b) specified in Section~\ref{sec:f8} 
within the figure-eight class leading to spherical (top) and planar (bottom) 
stationary configurations.} 
%
\end{figure}

\begin{figure}
\includegraphics[scale=.475]{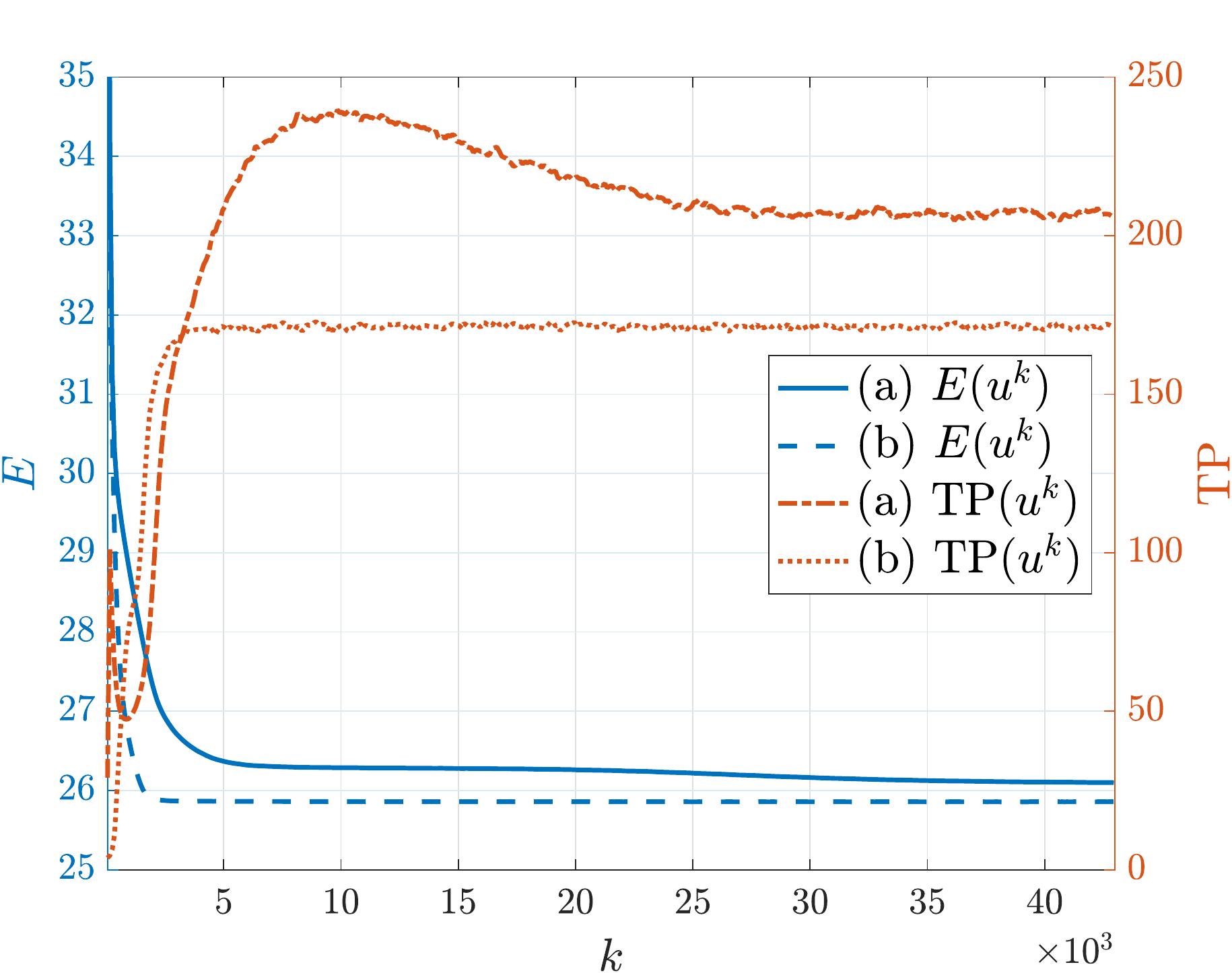}
\caption{The energy plots for the evolutions defined by the initial
configurations~(a) and~(b) in Section~\ref{sec:f8} (top and bottom parts
of Figure~\ref{fig:f8(b)}) give rise to assuming that the stationary
configuration obtained in~(b) could be the global minimizer in the 
figure-eight class.}\label{fig:f8}
\end{figure}
